\documentclass[11pt,a4paper]{amsart}

\usepackage[margin=2.5cm]{geometry}
\usepackage[utf8]{inputenc}
\usepackage[T1]{fontenc}
\usepackage{amsmath,amssymb,amsthm}
\usepackage{tikz}
\usepackage{booktabs,longtable,array}
\usepackage{graphicx}
\usepackage{seqsplit}
\usepackage{float}
\usepackage{placeins}
\usepackage[shortlabels]{enumitem}
\newcommand{\blocktables}{\setlength{\aboverulesep}{0pt}\setlength{\belowrulesep}{0pt}\renewcommand{\arraystretch}{1.15}}
\usepackage{xcolor}
\usepackage[colorlinks=true,linkcolor=blue!60!black,citecolor=blue!60!black,urlcolor=blue!60!black]{hyperref}

\newtheorem{theorem}{Theorem}[section]
\newtheorem{proposition}[theorem]{Proposition}
\newtheorem{lemma}[theorem]{Lemma}
\newtheorem{corollary}[theorem]{Corollary}
\theoremstyle{definition}
\newtheorem{definition}[theorem]{Definition}
\newtheorem{question}[theorem]{Question}
\newtheorem{problem}[theorem]{Problem}
\theoremstyle{remark}
\newtheorem{remark}[theorem]{Remark}

\newcommand{\vF}{\vec F}
\newcommand{\vFt}{\vec F_3}
\newcommand{\ind}{\operatorname{ind}}

\expandafter\def\csname witness@10-45\endcsname{\texttt{\seqsplit{0,10,110,11100,111010,111011,111100,111101000,111101001,111101010,111101011,1111011,11111}}\ $|$\ \texttt{\seqsplit{00000,00001,000100,000101,00011,0010,0011,010,0110,0111,10,110,111}}}
\expandafter\def\csname witness@10-163\endcsname{\texttt{\seqsplit{0,10000,1000100,1000101,100011,10010,10011,1010,1011,11}}\ $|$\ \texttt{\seqsplit{000,001,01,1000,10010,10011,1010,1011,110,111}}}
\expandafter\def\csname witness@10-132\endcsname{\texttt{\seqsplit{0,100,101,1100,11010,110110,1101110,1101111,111}}\ $|$\ \texttt{\seqsplit{00,01,100000,100001,10001,1001,101,110,111}}}
\expandafter\def\csname witness@10-156\endcsname{\texttt{\seqsplit{0,1000,100100,100101,10011,101,1100,11010,11011,111}}\ $|$\ \texttt{\seqsplit{000,001,01,1000,100100,1001010,1001011,10011,101,11}}}
\expandafter\def\csname witness@8-17-F\endcsname{\texttt{\seqsplit{0,100,101,1100,1101,11100,111010,111011,1111}}\ $|$\ \texttt{\seqsplit{000,001,0100,0101,01100,01101,0111,10,11}}}
\expandafter\def\csname witness@8-17-F3\endcsname{\texttt{\seqsplit{0,1,20,21,220,221,22200,22201,22202,22210,22211,22212,2222}}\ $|$\ \texttt{\seqsplit{0,100,101,10200,10201,102020,102021,102022,1021,1022,11,12,2}}}
\expandafter\def\csname witness@8-17-Fvec\endcsname{\texttt{\seqsplit{00,010,011,100,10100,1010100,10101010,10101011,101011,10110,10111,110,111}}\ $|$\ \texttt{\seqsplit{0000,000100,0001010,0001011,00011,0010,0011,01000,010010,010011,0101,011,1}}}
\expandafter\def\csname witness@8-17-F3vec\endcsname{\texttt{\seqsplit{0,1,200,201,202,210,211,212,220,2210,2211,2212,222}}\ $|$\ \texttt{\seqsplit{00,0100,0101,0102,011,012,0200,0201,0202,021,022,1,2}}}
\expandafter\def\csname witness@8-17-Fvec-K\endcsname{\texttt{\seqsplit{0,1000,10010,10011,1010,1011,11000,110010,110011,11010,110110,110111,111}}\ $|$\ \texttt{\seqsplit{0000,00010,00011,001000,00100100,001001010,001001011,0010011,001010,001011,0011,01,1}}}
\expandafter\def\csname witness@8-17-Fvec-Kr\endcsname{\texttt{\seqsplit{0000,000100,0001010,0001011,00011,0010,0011,01000,010010,010011,0101,011,1}}\ $|$\ \texttt{\seqsplit{00,010,011,100,10100,1010100,10101010,10101011,101011,10110,10111,110,111}}}
\expandafter\def\csname witness@8-17-F3vec-K\endcsname{\texttt{\seqsplit{00,0100,0101,0102,011,012,0200,0201,0202,021,022,1,2}}\ $|$\ \texttt{\seqsplit{0,1,200,201,202,210,211,212,220,2210,2211,2212,222}}}
\expandafter\def\csname witness@8-17-F3vec-Kr\endcsname{\texttt{\seqsplit{00,0100,0101,0102,011,01200,01201,01202,0121,0122,020,021,022,1,20,21,22}}\ $|$\ \texttt{\seqsplit{0,1,20,21,220,221,22200,2220100,2220101,2220102,222011,222012,22202,22210,22211,22212,2222}}}
\expandafter\def\csname witness@11n-42-Fvec\endcsname{\texttt{\seqsplit{000,0010,0011,0100,0101000,01010010,01010011,0101010,0101011,01011,011,100,101,110,111}}\ $|$\ \texttt{\seqsplit{00,0100,010100,0101010,0101011,01011000,010110010,010110011,010110100,0101101010,0101101011,01011011,010111,011,1}}}

\DeclareRobustCommand{\witness}[1]{{\footnotesize\csname witness@#1\endcsname}}

\title{Thompson indices of knots and the invertibility problem}
\author{Valeriano Aiello}
\date{\today}

\begin{document}

\begin{abstract}
Jones showed that Thompson's group $F$ produces all unoriented links, and introduced its oriented subgroup
$\vF$, which produces all oriented links; the same holds for the Brown--Thompson group $F_3$ and its
oriented subgroup $\vFt$. The corresponding Thompson indices measure the complexity of a knot in these
groups, much as the braid index does for braid groups. 
In this article we prove that the span of the Jones polynomial plus $1$ is a lower bound for the $F$- and
$\vF$-index, and we give lower bounds for the $\vF$- and $\vFt$-index in terms of the degrees of the HOMFLY-PT
polynomial. We then compute the four indices --- the $F$-, $F_3$-, $\vF$- and $\vFt$-index --- for the prime
knots with at most 13 crossings, by an exhaustive enumeration of tree diagrams (up to 13 leaves for $F$, 17 for
$\vF$, with the HOMFLY-PT polynomials of all the diagrams with 19 leaves, and 9 and 11 for $F_3$ and $\vFt$) completed by targeted searches; every value is rigorously certified, by exact computations
and by verified computations of hyperbolic structures in interval arithmetic. The $F$-index is determined for all knots with at most 11 crossings, the $\vF$-index for all knots
with at most 8 crossings, and the $F_3$- and $\vFt$-index for all knots with at most 12 crossings. The main outcome of these computations is that the
oriented indices are sensitive to the orientation: we find hundreds of non-invertible knots $K$ whose reverse
$K^r$ has a different $\vF$- or $\vFt$-index, the smallest being $8_{17}$ (with $\vFt$-index $7$ and $9$). The
$\vF$-index also distinguishes the Conway knot from its mutant, the Kinoshita--Terasaka knot.
\end{abstract}

\maketitle
\tableofcontents

\section{Introduction}\label{sec:intro}

A knot $K\subset S^3$ is \emph{invertible} (or reversible) if it is isotopic to its reverse $K^r$, the same
knot with the opposite orientation. We write $K^m$ for the mirror image of $K$ and $K^{mr}=(K^m)^r=(K^r)^m$; for a
link $L$, $L^m$ is its mirror image and, if $L$ is oriented, $L^r$ is obtained by reversing the orientation of
all its components. The first examples of non-invertible knots are due to the work of Trotter \cite{Trotter1963}. Later, Kawauchi showed that $8_{17}$, the first
candidate in the tables, is non-invertible~\cite{Kawauchi1979}; see also Hartley~\cite{Hartley1983}. Today
  invertible knots appear to be the exception: in KnotInfo~\cite{KnotInfo}, 8\,286 of the 12\,965 prime
knots with at most 13 crossings are not invertible, and as the crossing number goes to infinity the
proportion of invertible knots is expected to go to zero.

Yet non-invertibility is notoriously hard to detect. The classical polynomial invariants (Alexander, Jones,
HOMFLY-PT, Kauffman) do not change under reversal, and neither do most quantum invariants. For hyperbolic
knots the question can be decided from the symmetry group of the complement~\cite{HenryWeeks1992,
Kuperberg1996}, which is how SnapPy~\cite{SnapPy} answers it in practice; 
algebraically, a particularly effective tool is given by Eisermann's knot colouring
polynomials~\cite{Eisermann2007}, which record the images of the longitude under the representations of the knot group
into a finite group $G$ sending the meridian to a fixed element; they detect non-invertibility, as for $8_{17}$
\cite[\S2.4]{Eisermann2007}, provided one finds a finite group with suitable representations (note that
in~\cite{Eisermann2007}, following Conway, a knot with $K\neq K^r$ is called non-reversible). The problem of detecting
non-invertibility is thus interesting, important and subtle.

In a series of works on unitary representations of Thompson's groups started with~\cite{Jones2017}, Jones introduced a
general machinery producing actions of groups of fractions of categories, and observed that, applied to the
category of binary forests and planar tangles, it gives a way to construct a knot or link from every element of
Thompson's group $F$, in a way that is parallel to the closure of braids~\cite{Jones2017,Jones2019} (see~\cite{Aiello2025} for a survey of the field). An element
$g\in F$ is represented by a binary tree diagram (i.e., a pair of binary trees $(T_+,T_-)$) and yields an
unoriented link $L(g)$. Jones proved that every unoriented link arises in this way, the analogue of Alexander's
theorem for braids. These links do not come with a natural orientation, however the elements of Jones'
\emph{oriented subgroup} $\vF\le F$ produce oriented links. Jones originally proved that every oriented link can
be obtained from $\vF$ up to disjoint union with unknots~\cite{Jones2017}, and later the author showed that every
oriented link arises from $\vF$~\cite{Aiello2020}.  In this new framework, the role of positive braids is taken by the positive Thompson
monoids, which were studied by Baader and the author~\cite{AielloBaader2022,AielloBaader2024}.

In~\cite{Jones2019} the unoriented construction is extended to the Brown--Thompson group $F_3$, where it becomes
particularly transparent (each internal vertex of a ternary tree becomes a crossing), and the oriented one to its
oriented subgroup $\vFt$, defined by the condition that the checkerboard shading gives a Seifert surface (see
also~\cite{AielloNagnibeda2022} for a study of this subgroup). Hence a knot $K$ has a \emph{Thompson index} in each
of these groups: the minimal size of a tree diagram realizing it. Jones proposed to reboot the theory of braids
and links with Thompson groups in place of the braid groups, and listed the Thompson indices among the first
questions, next to the analogues of Markov's and Birman's theorems~\cite[Section~7]{Jones2019}. The very first bound for one of these indices, in particular $F$-index (Theorem~\ref{thm:unoriented}\ref{it:golan-sapir}), was discovered by Golan and Sapir. More recently, Luo and Wan related the Thompson
indices to grid diagrams, the braid index and the Seifert genus~\cite{LuoWan2024}. 

Having new constructors of oriented links, we have a potentially new way to approach the invertibility problem. This motivates the following problem.

\begin{problem}[J. Birman, personal communication]\label{prob:invariants}
Find a complete invariant, or partial invariants, of the oriented knots $\vec K$ and $\vec K^r$ (two
orientations of the same knot) in terms of the elements of the oriented Thompson group $\vF$ that realize them.
In particular, find computable invariants of elements of $\vF$ that detect non-invertibility.
\end{problem}

The $\vF$-index and the $\vFt$-index are the simplest such invariants: they are invariants of oriented knots,
and computing them also helps to understand Jones' construction better. This paper shows that they detect
non-invertibility in many cases. Another possible choice would be to count the number of tree diagrams that realize the same index, but it would require more computations.

The results of this paper are of two kinds: lower bounds for the indices, proved in Section~\ref{sec:known}, and the
computation of the four indices for the prime knots with at most 13 crossings (Section~\ref{sec:computation}):
all the tree diagrams of $F$ with at most 13 leaves, of $\vF$ with at most 17 (and the HOMFLY-PT polynomials of
all those with 19), of $F_3$ with $n\le 9$ and of $\vFt$ with $n\le 11$ are enumerated, the knots are identified with KnotInfo and SnapPy, and single open cases
beyond the enumeration are settled by targeted searches (Section~\ref{sec:targeted}). The main findings are the
following.
\begin{itemize}
\item \textbf{Lower bounds from polynomial invariants.} For every link $L$,
  $\ind_{\vF}(L)\ge\ind_F(L)\ge\operatorname{span}V_L+1$, where $V_L$ is the Jones polynomial
  (Theorem~\ref{thm:span}); hence $\ind_F(L)\ge c(L)+1$ for a non-split alternating link, and
  $\ind_{\vF}(K)\ge c(K)+2$ for an alternating knot with an odd  crossing number
  (Corollary~\ref{cor:alternating}). The bound of Theorem~\ref{thm:span} is attained by 1\,440 of the knots
  with known $F$-index, and in $\vF$ (raised to the next odd integer) by 8 knots, 7 of them alternating with
  $\ind_{\vF}(K)=c(K)+2$; for these knots the index is proved by a single tree diagram, without any
  enumeration. The degrees of the HOMFLY-PT polynomial give
  lower bounds for the oriented indices (Proposition~\ref{prop:homfly}).
\item \textbf{The oriented indices detect non-invertibility.} Among the non-invertible knots with at most
  13 crossings, 311 have orientations with different $\vF$-index and 946 have orientations with different
  $\vFt$-index (Section~\ref{sec:oriented} and Tables~\ref{tab:oriented-Fvec}, \ref{tab:oriented-F3vec} in Appendix~\ref{app:tables}).
  The smallest example is the first non-invertible knot, $8_{17}$: its $\vFt$-index is $7$ for $K=K^{mr}$ and
  $9$ for $K^r=K^m$ (Figure~\ref{fig:817-F3vec}). All these results are rigorously certified
  (Section~\ref{sec:certified}). For $\vF$ the smallest examples have 10 crossings
  ($10_{88}$: $15$ and $13$). This answers Question~1 of~\cite{Aiello2026}; all the differences found are equal
  to $2$, in accordance with Theorem~\ref{thm:oriented}\ref{it:parity} and~\ref{it:reverse}.
  Since the oriented indices of $K$ and $K^{mr}$ coincide (Theorem~\ref{thm:oriented}\ref{it:inverse}), each of
  these knots also has $\ind(K^m)\ne\ind(K)$: the oriented indices detect chirality as well. This is a
  difference with the braid index, which sees neither the reverse nor the mirror image: if $K$ is the closure
  of a braid $\beta$ on $b$ strands, then $K^r$ is the closure of the braid obtained by reading $\beta$
  backwards and $K^m$ is the closure of the braid obtained by replacing each $\sigma_i$ with $\sigma_i^{-1}$,
  both on $b$ strands.
\item \textbf{The $\vF$-index distinguishes mutants.} The Conway knot $11n_{34}$ has $\vF$-index $\ge 19$
  in each orientation, while every orientation of its mutant, the Kinoshita--Terasaka knot $11n_{42}$, has
  $\vF$-index $15$ (Theorem~\ref{thm:mutants}). Of the
  1\,125 pairs of mutant knots with at most 13 crossings~\cite{Stoimenow2010}, 458 are distinguished by one of
  the indices, among them 12 of the 16 pairs with 11 crossings (Table~\ref{tab:mutants} in Appendix~\ref{app:tables}).
\item \textbf{Tables.} The $F$-index is known for all the knots with at most 11 crossings, the $F_3$-index for all the knots with at most 12 crossings; the
  $\vFt$-index for all the knots with at most 12 crossings, in both orientations; the $\vF$-index for all
  the knots with at most 8 crossings, for 36 of the 49 knots with 9 crossings and for 128 of the 165 with 10
  (Table~\ref{tab:summary} and Table~\ref{tab:small} in Appendix~\ref{app:tables}). All the values are rigorously certified by computer
  (Section~\ref{sec:certified-tables}): a diagram is shown to be a knot $K$, or not to be $K$, by exact invariants
  (HOMFLY-PT polynomial, knot Floer homology, crossing number after simplification) or by hyperbolic structures,
  volumes and canonical triangulations computed in interval arithmetic, so that no step relies on floating-point
  approximations.
\item We answer a question of Jones~\cite[Question 7]{Jones2019}: we compute the $F_3$-index of the Borromean rings (it turns out to be $6$, see Theorem~\ref{thm:borromean}). We also confirm that the  $F$-index is $9$ (a result originally found by
  Sweeney).
\item We compare the indices with the bounds of Golan--Sapir~\cite{GolanSapir2017} and
  Luo--Wan~\cite{LuoWan2024} (Section~\ref{sec:bounds}): the braid index is always strictly smaller than
  each index; the lower bound $c(K)\le 2(\ind(K)-1)$ is attained by 105 knots for $F_3$; the genus bound
  $g(K)\le(\ind(K)-1)/2$ is attained by 2\,571 knots for $\vFt$ and by none for $\vF$.
\item We correct two entries of an earlier table of $F$-indices~\cite{Sweeney2021} and complete it with two
  values that were missing (Appendix~\ref{app:sweeney}).
\end{itemize}

All the code and the tables are available in the repository~\cite{Repository}, and a selection of the data (all
the HOMFLY-PT classes of knot diagrams of every size, with one representative each) is archived on
Zenodo~\cite{Data}.

\subsection*{Acknowledgements}
The author thanks Joan Birman for conversations following his talk at the Geometric Topology Seminar of
Columbia University (online, November 2021), which further motivated his interest in the oriented indices,
and for suggesting Problem~\ref{prob:invariants}.

The author is partially supported by Sapienza Università di Roma (Progetto di Ateneo Dipartimentale 2024,
``New research trends in Mathematics at Castelnuovo''). This work was also supported by the Gruppo Nazionale
per l'Analisi Matematica, la Probabilità e le loro Applicazioni (GNAMPA--INdAM) through the program
``Partecipazione a Convegni, Scuole, Workshop e Cicli di Seminari'' and by the INdAM/GNAMPA project
``Simmetrie distribuzionali per processi stocastici quantistici'' CUP E53C25002010001. 

The first computations of the 
indices were carried out between October 2022 and
April 2023, during an internship of the author at CERN; the author thanks CERN, and in
particular Maurizio Pierini, for providing the computing resources. These computations were later resumed with the help of Claude Code, which made it possible to extend them 
to larger diagrams, and to make them considerably
more efficient. The author acknowledges the Claude Team plan for scientists and thanks the Facoltà di Scienze
Matematiche, Fisiche e Naturali of Sapienza Università di Roma for the subscription to Claude. Part of the targeted searches of
Section~\ref{sec:targeted} ran on the server Lagrange of the Dipartimento di Matematica Guido Castelnuovo,
Sapienza Università di Roma, which the author thanks for the computing resources. 

\section{Knots from Thompson groups}\label{sec:construction}

\subsection{Tree diagrams}
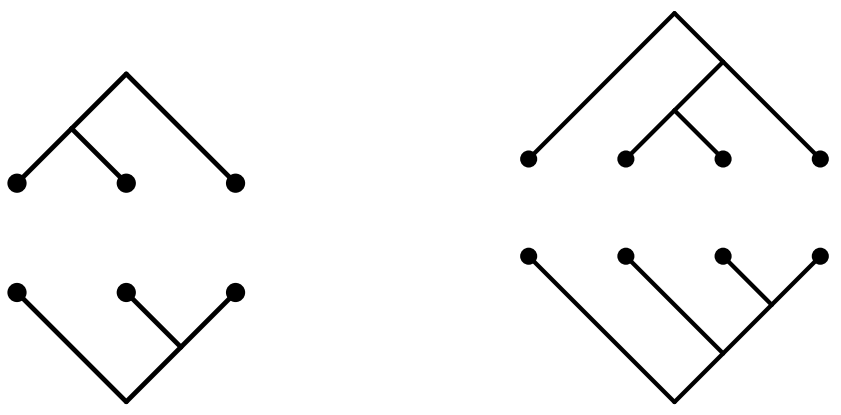
\begin{figure}[ht]
\centering
\begin{tabular}{c@{\qquad\qquad}c}
\resizebox{!}{2.2cm}{\begin{tikzpicture}[scale=0.32, line cap=round]
  \draw (0.5,1) -- (0,0.5);
  \draw (0.5,1) -- (1,0.5);
  \draw (1,1.5) -- (0.5,1);
  \draw (1,1.5) -- (2,0.5);
  \draw (1.5,-1) -- (1,-0.5);
  \draw (1.5,-1) -- (2,-0.5);
  \draw (1,-1.5) -- (0,-0.5);
  \draw (1,-1.5) -- (1.5,-1);
  \fill (0,0.5) circle (2.5pt) (0,-0.5) circle (2.5pt);
  \fill (1,0.5) circle (2.5pt) (1,-0.5) circle (2.5pt);
  \fill (2,0.5) circle (2.5pt) (2,-0.5) circle (2.5pt);
\end{tikzpicture}} & \resizebox{!}{2.6cm}{\begin{tikzpicture}[scale=0.32, line cap=round]
  \draw (1.5,1) -- (1,0.5);
  \draw (1.5,1) -- (2,0.5);
  \draw (2,1.5) -- (1.5,1);
  \draw (2,1.5) -- (3,0.5);
  \draw (1.5,2) -- (0,0.5);
  \draw (1.5,2) -- (2,1.5);
  \draw (2.5,-1) -- (2,-0.5);
  \draw (2.5,-1) -- (3,-0.5);
  \draw (2,-1.5) -- (1,-0.5);
  \draw (2,-1.5) -- (2.5,-1);
  \draw (1.5,-2) -- (0,-0.5);
  \draw (1.5,-2) -- (2,-1.5);
  \fill (0,0.5) circle (2.5pt) (0,-0.5) circle (2.5pt);
  \fill (1,0.5) circle (2.5pt) (1,-0.5) circle (2.5pt);
  \fill (2,0.5) circle (2.5pt) (2,-0.5) circle (2.5pt);
  \fill (3,0.5) circle (2.5pt) (3,-0.5) circle (2.5pt);
\end{tikzpicture}} \\[1mm]
$x_0$ & $x_1$
\end{tabular}
\caption{The binary tree diagrams of the generators $x_0$ and $x_1$ of $F$.}
\label{fig:trees-F}
\end{figure}

\begin{figure}[ht]
\centering
\begin{tikzpicture}[scale=0.6, line cap=round]
  \draw (0,0.4) -- (0.5,1.1) -- (1,0.4);
  \draw (0,-0.4) -- (0.5,-1.1) -- (1,-0.4);
  \draw (0.5,1.1) -- (0.5,1.6); \draw (0.5,-1.1) -- (0.5,-1.6);
  \fill (0,0.4) circle (2pt) (1,0.4) circle (2pt) (0,-0.4) circle (2pt) (1,-0.4) circle (2pt);
  \node at (2,0) {$\longleftrightarrow$};
  \draw (3,1.6) -- (3,0.4); \draw (3,-0.4) -- (3,-1.6);
  \fill (3,0.4) circle (2pt) (3,-0.4) circle (2pt);
  \begin{scope}[xshift=6.5cm]
  \draw (0,0.4) -- (1,1.1) -- (2,0.4); \draw (1,1.1) -- (1,0.4);
  \draw (0,-0.4) -- (1,-1.1) -- (2,-0.4); \draw (1,-1.1) -- (1,-0.4);
  \draw (1,1.1) -- (1,1.6); \draw (1,-1.1) -- (1,-1.6);
  \fill (0,0.4) circle (2pt) (1,0.4) circle (2pt) (2,0.4) circle (2pt);
  \fill (0,-0.4) circle (2pt) (1,-0.4) circle (2pt) (2,-0.4) circle (2pt);
  \node at (3,0) {$\longleftrightarrow$};
  \draw (4,1.6) -- (4,0.4); \draw (4,-0.4) -- (4,-1.6);
  \fill (4,0.4) circle (2pt) (4,-0.4) circle (2pt);
  \end{scope}
\end{tikzpicture}
\caption{A pair of opposing carets in a binary tree diagram (left) and in a ternary tree diagram (right). Adding or
removing such a pair does not change the element of $F$ or of $F_3$.}
\label{fig:carets}
\end{figure}
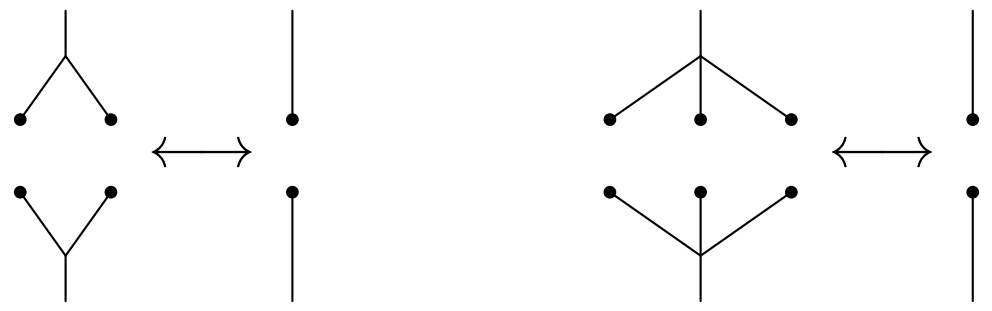

Elements of $F$ are represented by \emph{binary tree diagrams}: pairs $(T_+,T_-)$ of planar rooted binary trees
with the same number $n$ of leaves~\cite{CFP1996,Belk2007}. We draw $T_+$ above and $T_-$ below, upside down, with
the leaves on a horizontal line; $n$ is the \emph{size} of the diagram. Figure~\ref{fig:trees-F} shows the tree
diagrams of the generators $x_0$ and $x_1$ of $F$. Two tree diagrams represent the same element of $F$ if and only
if they are related by adding or removing pairs of \emph{opposing carets}, that is, a caret of $T_+$ and a caret of
$T_-$ joining the same two consecutive leaves (Figure~\ref{fig:carets}); every element of $F$ has a unique
\emph{reduced} tree diagram, without opposing carets~\cite{CFP1996}. Taking the inverse of an element amounts to exchanging the two trees: $(T_+,T_-)^{-1}=(T_-,T_+)$. To multiply $g=(T_+,T_-)$ and $h=(S_+,S_-)$, one first adds pairs of opposing carets to the two diagrams until $T_-=S_+$, and then $gh=(T_+,S_-)$; here $gh$ means that $g$ is applied first, $T_+$ being the domain tree~\cite{CFP1996,Belk2007}.

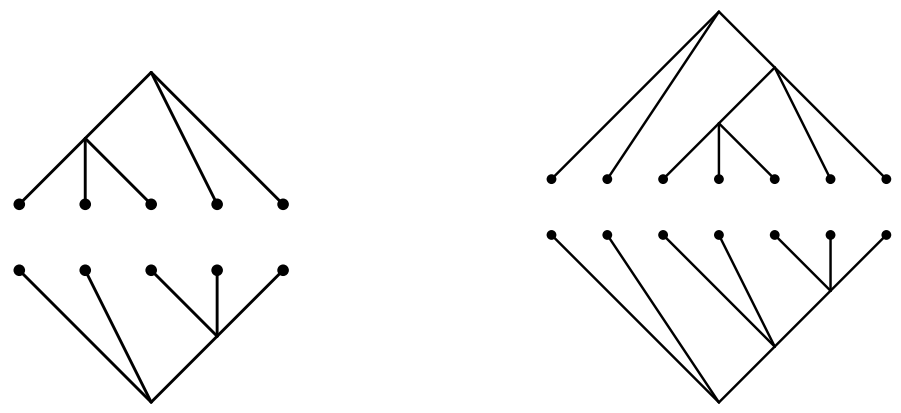
\begin{figure}[ht]
\centering
\begin{tabular}{c@{\qquad\qquad}c}
\resizebox{!}{2.2cm}{\begin{tikzpicture}[scale=0.32, line cap=round]
  \draw (1,1.5) -- (0,0.5);
  \draw (1,1.5) -- (1,0.5);
  \draw (1,1.5) -- (2,0.5);
  \draw (2,2.5) -- (1,1.5);
  \draw (2,2.5) -- (3,0.5);
  \draw (2,2.5) -- (4,0.5);
  \draw (2,-2.5) -- (0,-0.5);
  \draw (2,-2.5) -- (1,-0.5);
  \draw (3,-1.5) -- (2,-0.5);
  \draw (3,-1.5) -- (3,-0.5);
  \draw (3,-1.5) -- (4,-0.5);
  \draw (2,-2.5) -- (3,-1.5);
  \fill (0,0.5) circle (2.5pt) (0,-0.5) circle (2.5pt);
  \fill (1,0.5) circle (2.5pt) (1,-0.5) circle (2.5pt);
  \fill (2,0.5) circle (2.5pt) (2,-0.5) circle (2.5pt);
  \fill (3,0.5) circle (2.5pt) (3,-0.5) circle (2.5pt);
  \fill (4,0.5) circle (2.5pt) (4,-0.5) circle (2.5pt);
\end{tikzpicture}} & \resizebox{!}{2.6cm}{\begin{tikzpicture}[scale=0.32, line cap=round]
  \draw (3,3.5) -- (0,0.5);
  \draw (3,3.5) -- (1,0.5);
  \draw (3,1.5) -- (2,0.5);
  \draw (3,1.5) -- (3,0.5);
  \draw (3,1.5) -- (4,0.5);
  \draw (4,2.5) -- (3,1.5);
  \draw (4,2.5) -- (5,0.5);
  \draw (4,2.5) -- (6,0.5);
  \draw (3,3.5) -- (4,2.5);
  \draw (3,-3.5) -- (0,-0.5);
  \draw (3,-3.5) -- (1,-0.5);
  \draw (4,-2.5) -- (2,-0.5);
  \draw (4,-2.5) -- (3,-0.5);
  \draw (5,-1.5) -- (4,-0.5);
  \draw (5,-1.5) -- (5,-0.5);
  \draw (5,-1.5) -- (6,-0.5);
  \draw (4,-2.5) -- (5,-1.5);
  \draw (3,-3.5) -- (4,-2.5);
  \fill (0,0.5) circle (2.5pt) (0,-0.5) circle (2.5pt);
  \fill (1,0.5) circle (2.5pt) (1,-0.5) circle (2.5pt);
  \fill (2,0.5) circle (2.5pt) (2,-0.5) circle (2.5pt);
  \fill (3,0.5) circle (2.5pt) (3,-0.5) circle (2.5pt);
  \fill (4,0.5) circle (2.5pt) (4,-0.5) circle (2.5pt);
  \fill (5,0.5) circle (2.5pt) (5,-0.5) circle (2.5pt);
  \fill (6,0.5) circle (2.5pt) (6,-0.5) circle (2.5pt);
\end{tikzpicture}} \\[1mm]
$y_0=\iota(x_0)$ & $y_2=\iota(x_1)$
\end{tabular}
\caption{The ternary tree diagrams of the elements $y_0=\iota(x_0)$ and $y_2=\iota(x_1)$ of $F_3$.}
\label{fig:trees-F3}
\end{figure}

Similarly, elements of the Brown--Thompson group $F_3$ are represented by pairs of planar rooted ternary trees with
the same number of leaves, up to adding or removing pairs of opposing ternary carets (Figure~\ref{fig:carets}), and
every element of $F_3$ has a unique reduced tree diagram~\cite{Brown87}. A ternary tree with $n-1$ internal
vertices has $2n-1$ leaves, and we call $n$ the size of the diagram. Jones' embedding $\iota\colon F\to F_3$
replaces every binary caret by a ternary caret whose middle edge ends in a leaf; it preserves the size.
Figure~\ref{fig:trees-F3} shows the ternary tree diagrams of $y_0=\iota(x_0)$ and $y_2=\iota(x_1)$.

\subsection{The unoriented link}
We recall from~\cite[Section~2]{Aiello2025} Jones' construction of a link from a binary tree diagram
$(T_+,T_-)\in F$ of size $n$~\cite{Jones2017}. Place the leaves on the points $1,\dots,n$ of the $x$-axis.

\begin{definition}
The \emph{Tait diagram} $\Gamma(T_+,T_-)$ is the plane graph with vertices $(j-\frac12,0)$, $j=1,\dots,n$, that is,
one vertex immediately to the left of each leaf, and with the following edges:
\begin{itemize}
\item for every edge of $T_+$ sloping up from left to right (a \emph{West--North} edge, i.e.\ the left edge of a
  caret), an edge in the upper half-plane crossing it: if the left subtree of the caret has leaves $a,\dots,b$, it
  joins $(a-\frac12,0)$ and $(b+\frac12,0)$;
\item similarly, for every edge of $T_-$ sloping down from left to right (a \emph{West--South} edge), an edge in the
  lower half-plane crossing it.
\end{itemize}
We denote by $\Gamma_+(T_+)$ and $\Gamma_-(T_-)$ the subgraphs in the upper and in the lower half-plane. The Tait
diagram is a signed graph: the edges of $\Gamma_+(T_+)$ are positive and those of $\Gamma_-(T_-)$ are negative.
\end{definition}

Each of $\Gamma_+(T_+)$ and $\Gamma_-(T_-)$ has the same $n$ vertices, and every vertex other than the leftmost one
is joined to exactly one vertex on its left~\cite[Proposition~4.1.3]{Jones2017},~\cite[Remark~2.1]{Aiello2025}; in
particular both are spanning trees of $\Gamma(T_+,T_-)$, with $n-1$ edges each. Moreover $(T_+,T_-)\mapsto
\Gamma(T_+,T_-)$ is a bijection between binary tree diagrams and graphs of this form~\cite[Lemma~4.1.4]{Jones2017}.

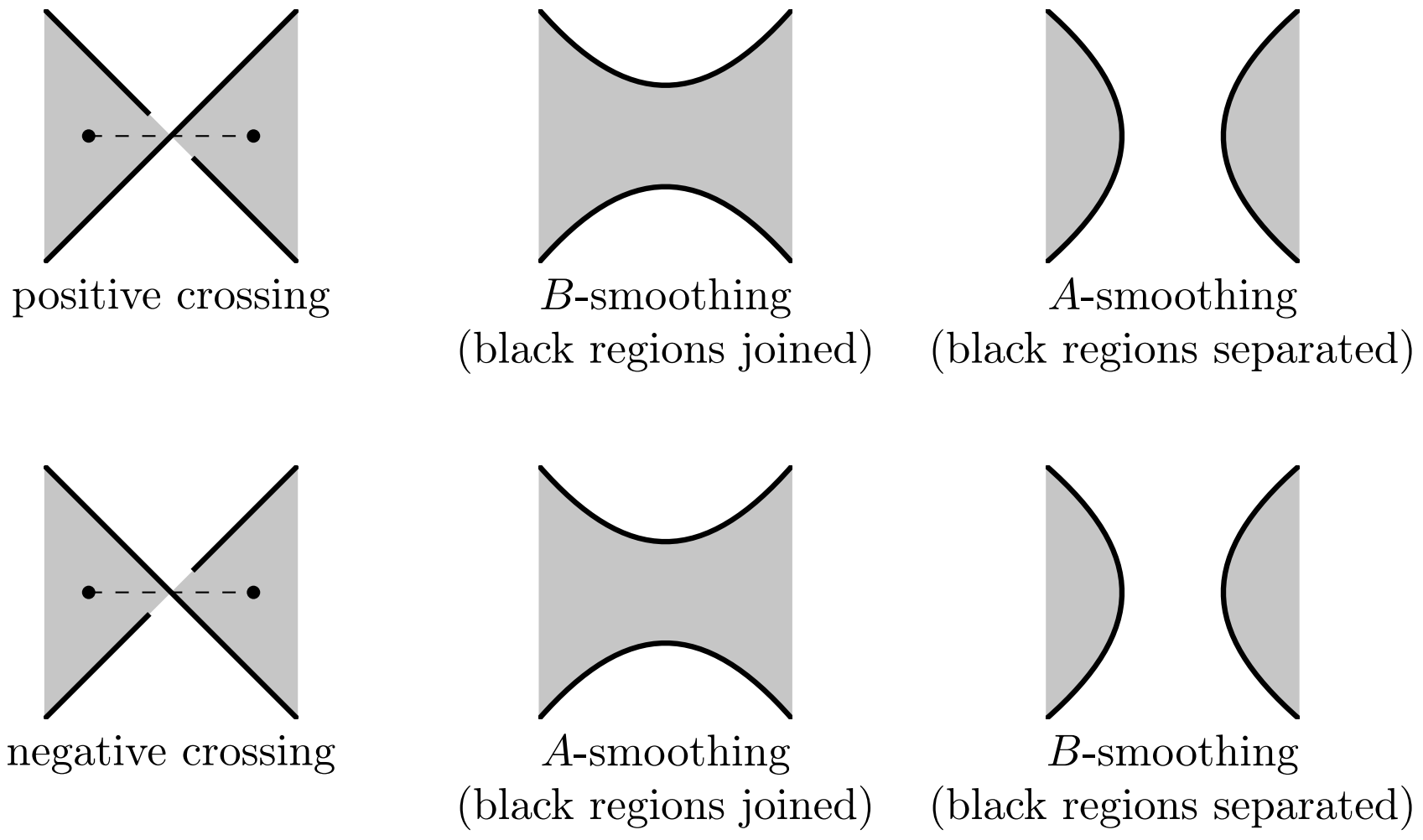
\begin{figure}[ht]
\centering
\begin{tikzpicture}[scale=1.0]
  \begin{scope}
    \clip (-1,-1) rectangle (1,1);
    \fill[gray!45] (-1,1) -- (0,0) -- (-1,-1) -- cycle;
    \fill[gray!45] (1,1) -- (0,0) -- (1,-1) -- cycle;
    \draw[very thick] (-1,-1) -- (1,1);
    \draw[very thick] (-1,1) -- (-0.17,0.17);
    \draw[very thick] (0.17,-0.17) -- (1,-1);
    \draw[dashed] (-0.65,0) -- (0.65,0);
    \fill (-0.65,0) circle (1.5pt) (0.65,0) circle (1.5pt);
  \end{scope}
  \node[font=\small,align=center] at (0,-1.5) {positive crossing\\\phantom{(}};
  \begin{scope}[xshift=3.9cm]
    \clip (-1,-1) rectangle (1,1);
    \fill[gray!45] (-1,1) .. controls (-0.3,0.2) and (0.3,0.2) .. (1,1)
                  -- (1,-1) .. controls (0.3,-0.2) and (-0.3,-0.2) .. (-1,-1) -- cycle;
    \draw[very thick] (-1,1) .. controls (-0.3,0.2) and (0.3,0.2) .. (1,1);
    \draw[very thick] (-1,-1) .. controls (-0.3,-0.2) and (0.3,-0.2) .. (1,-1);
  \end{scope}
  \node[font=\small,align=center] at (3.9,-1.5) {$B$-smoothing\\(black regions joined)};
  \begin{scope}[xshift=7.9cm]
    \clip (-1,-1) rectangle (1,1);
    \fill[gray!45] (-1,1) .. controls (-0.2,0.3) and (-0.2,-0.3) .. (-1,-1) -- cycle;
    \fill[gray!45] (1,1) .. controls (0.2,0.3) and (0.2,-0.3) .. (1,-1) -- cycle;
    \draw[very thick] (-1,1) .. controls (-0.2,0.3) and (-0.2,-0.3) .. (-1,-1);
    \draw[very thick] (1,1) .. controls (0.2,0.3) and (0.2,-0.3) .. (1,-1);
  \end{scope}
  \node[font=\small,align=center] at (7.9,-1.5) {$A$-smoothing\\(black regions separated)};
  \begin{scope}[yshift=-3.6cm]
    \clip (-1,-1) rectangle (1,1);
    \fill[gray!45] (-1,1) -- (0,0) -- (-1,-1) -- cycle;
    \fill[gray!45] (1,1) -- (0,0) -- (1,-1) -- cycle;
    \draw[very thick] (-1,1) -- (1,-1);
    \draw[very thick] (-1,-1) -- (-0.17,-0.17);
    \draw[very thick] (0.17,0.17) -- (1,1);
    \draw[dashed] (-0.65,0) -- (0.65,0);
    \fill (-0.65,0) circle (1.5pt) (0.65,0) circle (1.5pt);
  \end{scope}
  \node[font=\small,align=center] at (0,-5.1) {negative crossing\\\phantom{(}};
  \begin{scope}[xshift=3.9cm,yshift=-3.6cm]
    \clip (-1,-1) rectangle (1,1);
    \fill[gray!45] (-1,1) .. controls (-0.3,0.2) and (0.3,0.2) .. (1,1)
                  -- (1,-1) .. controls (0.3,-0.2) and (-0.3,-0.2) .. (-1,-1) -- cycle;
    \draw[very thick] (-1,1) .. controls (-0.3,0.2) and (0.3,0.2) .. (1,1);
    \draw[very thick] (-1,-1) .. controls (-0.3,-0.2) and (0.3,-0.2) .. (1,-1);
  \end{scope}
  \node[font=\small,align=center] at (3.9,-5.1) {$A$-smoothing\\(black regions joined)};
  \begin{scope}[xshift=7.9cm,yshift=-3.6cm]
    \clip (-1,-1) rectangle (1,1);
    \fill[gray!45] (-1,1) .. controls (-0.2,0.3) and (-0.2,-0.3) .. (-1,-1) -- cycle;
    \fill[gray!45] (1,1) .. controls (0.2,0.3) and (0.2,-0.3) .. (1,-1) -- cycle;
    \draw[very thick] (-1,1) .. controls (-0.2,0.3) and (-0.2,-0.3) .. (-1,-1);
    \draw[very thick] (1,1) .. controls (0.2,0.3) and (0.2,-0.3) .. (1,-1);
  \end{scope}
  \node[font=\small,align=center] at (7.9,-5.1) {$B$-smoothing\\(black regions separated)};
\end{tikzpicture}
\caption{A positive and a negative crossing of a shaded diagram, with the edge of the Tait graph joining the two
black regions (dashed), and their two smoothings. The $A$-smoothing joins the two regions swept out by the
overcrossing strand when it is rotated counterclockwise.}
\label{fig:smoothings}
\end{figure}

The link $L(T_+,T_-)$ is the link whose Tait graph is $\Gamma(T_+,T_-)$: its diagram is obtained from the medial
graph of $\Gamma(T_+,T_-)$ (one $4$-valent vertex on each edge of $\Gamma$, and an edge between two such vertices
whenever their edges are consecutive along a face of $\Gamma$) by turning every vertex into a crossing, so that in the
checkerboard shading in which the vertices of $\Gamma$ are the black regions, the crossings on the edges of
$\Gamma_+(T_+)$ are positive and those on the edges of $\Gamma_-(T_-)$ are negative, in the sense of
Figure~\ref{fig:smoothings}. Hence the diagram of $L(T_+,T_-)$ has $2(n-1)$ crossings.

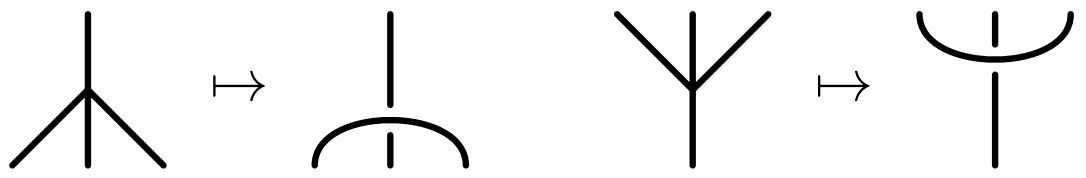
\begin{figure}[ht]
\centering
\begin{tikzpicture}[line width=1.2pt,line cap=round]
  \begin{scope}[xshift=4cm]
    \draw (0,0.5) -- (0,-0.5) (-0.5,-0.5) -- (0,0) -- (0.5,-0.5);
    \node at (1,0) {$\mapsto$};
    \draw (2,0.5) -- (2,-0.1) (2,-0.3) -- (2,-0.5);
    \draw (1.5,-0.5) .. controls (1.5,-0.1) and (2.5,-0.1) .. (2.5,-0.5);
  \end{scope}
  \begin{scope}[xshift=8cm]
    \draw (0,0.5) -- (0,-0.5) (-0.5,0.5) -- (0,0) -- (0.5,0.5);
    \node at (1,0) {$\mapsto$};
    \draw (2,-0.5) -- (2,0.1) (2,0.3) -- (2,0.5);
    \draw (1.5,0.5) .. controls (1.5,0.1) and (2.5,0.1) .. (2.5,0.5);
  \end{scope}
\end{tikzpicture}
\caption{The rules for obtaining $L(g)$:  
 the vertices of the top tree (left)
 and of the bottom tree (right) are replaced by ``forks''.}
\label{fig:forks}
\end{figure}

The construction extends to $F_3$~\cite{Jones2019}. Let $(T_+,T_-)$ be a pair of ternary trees with $2n-1$ leaves.
Join the two roots by an edge passing around the left of the picture: the result is a planar $4$-valent graph whose
vertices are the internal vertices of the two trees. The link $L(T_+,T_-)$ is obtained by turning every vertex into
a crossing, the strand joining the left and right children passing over the strand joining the parent and the middle
child (Figure~\ref{fig:forks}). For an element of $F$, $L(g)=L(\iota(g))$. Equivalently, the regions between consecutive leaves are
alternately black and white (the leftmost one being black); the $n$ black regions are the vertices of the Tait
graph $\Gamma(T_+,T_-)$, which has one edge for every internal vertex of $T_+$ (drawn above the line) and of $T_-$
(drawn below), and $L(g)$ is the link with Tait graph $\Gamma$. For $g\in F$ the Tait graph of $\iota(g)$ is the
Tait diagram defined above; for a general element of $F_3$, however, the edges in the upper (resp.\ lower)
half-plane need not be all positive (resp.\ negative)~\cite[Section~2]{Aiello2025}.

\begin{definition}
The \emph{$F$-index} $\ind_F(K)$ of a knot $K$ is the smallest number $n$ of leaves of a binary tree diagram
$(T_+,T_-)$ with $L(T_+,T_-)=K$. The \emph{$F_3$-index} $\ind_{F_3}(K)$ is the smallest $n$ such that $K$ is
$L(T_+,T_-)$ for ternary trees with $n-1$ internal vertices each~\cite{Aiello2025}. With this normalization
$\iota$ preserves $n$, so $\ind_{F_3}(K)\le\ind_F(K)$.
\end{definition}

A tree diagram of size $n$ has $2(n-1)$ crossings in both cases, and $(T_+,T_-)\mapsto(T_-,T_+)$ as well as
the vertical reflection of both trees replace $L(g)$ by its mirror image, so the two unoriented indices do
not distinguish a knot from its mirror image.

Concrete examples are in Section~\ref{sec:example817}: Figure~\ref{fig:817-unoriented} shows minimal tree
diagrams of the knot $8_{17}$ in $F$ and in $F_3$, together with the knot diagrams obtained from them.

\subsection{The oriented link}
We call an edge of a tree from a vertex to its left, middle or right child a \emph{left}, \emph{middle} or
\emph{right edge} (binary trees have no middle edges), and we number the leaves $1,2,\dots$ from left to right.
The element $g=(T_+,T_-)$ belongs to the oriented subgroup $\vF$ (resp.\ $\vFt$) if and only if its Tait graph
$\Gamma$ is \emph{$2$-colourable}, that is, its vertices can be coloured $+$ and $-$ so that adjacent vertices
have different colours~\cite{GolanSapir2017,Jones2019,AielloNagnibeda2022}. 
Since $\Gamma$ is connected, when it is $2$-colourable,  there is exactly one
such colouring in which the leftmost vertex has colour $+$. In terms of the trees:
\begin{itemize}
\item for $\vF$, let $r_\pm(j)$ be the number of right edges on the path from the root of $T_\pm$ to the leaf
  $j$. Then $g\in\vF$ if and only if $r_+(j)\equiv r_-(j)\pmod 2$ for every leaf $j$, and the Tait vertex to
  the left of the leaf $j$ has colour $+$ if and only if $r_\pm(j)$ is even;
\item for $\vFt$, cut the path from the root of $T_\pm$ to the leaf $j$ at its middle edges into pieces,
  numbered $1,2,\dots$ starting from the root, and let the weight $c_\pm(j)$ be the number of middle edges of
  the path, plus the number of right edges in the odd-numbered pieces and of left edges in the even-numbered
  ones. Then $g\in\vFt$ if and only if $c_+(j)\equiv c_-(j)\pmod 2$ for every leaf
  $j$~\cite[Prop.~1]{AielloNagnibeda2022}, and the Tait vertex to the left of an odd leaf $j$ has colour $+$
  if and only if $c_\pm(j)$ is even.
\end{itemize}
The black regions, joined by a twisted band at each crossing, form the \emph{checkerboard surface} of $L(g)$,
which is orientable when $g$ is in the oriented subgroup. The link $L(g)$ is then oriented as its boundary: the
boundary of a black region whose Tait vertex has colour $+$ runs counterclockwise, that of a region of colour
$-$ clockwise.

\begin{definition}
The \emph{$\vF$-index} $\ind_{\vF}(\vec K)$ and the \emph{$\vFt$-index} $\ind_{\vFt}(\vec K)$ of an oriented
knot $\vec K$ are defined as above, with tree diagrams in $\vF$ and $\vFt$ and $L(g)$ oriented by the
colouring.
\end{definition}

We write $\ind$ for any of the four indices when a statement holds for all of them, or for all the indices
under consideration (for instance the two oriented ones).

A tree diagram of size $\ind(K)$ realizing $K$ is a \emph{minimal diagram} of $K$; we also call it a
\emph{witness} of the index, and more generally a diagram realizing $K$ is a witness of the upper bound given by
its size.

By Theorem~\ref{thm:oriented}\ref{it:parity} below, the $\vF$- and $\vFt$-indices of a knot are odd: knots
only arise from tree diagrams of odd size, so only odd $n$ need to be considered for the oriented groups.

Examples of oriented knots obtained in this way are in Section~\ref{sec:example817}: Figure~\ref{fig:817-F3vec}
shows diagrams of $\vFt$ realizing the two orientations of $8_{17}$, and Figure~\ref{fig:1088-Fvec} diagrams of
$\vF$ realizing the two orientations of $10_{88}$.

\subsection{Example: the knot $8_{17}$}\label{sec:example817}
\begin{figure}[tbp]
\centering
\begin{tabular}{cc}
\resizebox{0.32\linewidth}{!}{\begin{tikzpicture}[scale=0.32, line cap=round]
  \draw (5,5.5) -- (0,0.5);
  \draw (2,1.5) -- (1,0.5);
  \draw (2,1.5) -- (2,0.5);
  \draw (2,1.5) -- (3,0.5);
  \draw (3,2.5) -- (2,1.5);
  \draw (3,2.5) -- (4,0.5);
  \draw (3,2.5) -- (5,0.5);
  \draw (5,5.5) -- (3,2.5);
  \draw (7,1.5) -- (6,0.5);
  \draw (7,1.5) -- (7,0.5);
  \draw (7,1.5) -- (8,0.5);
  \draw (8,2.5) -- (7,1.5);
  \draw (8,2.5) -- (9,0.5);
  \draw (8,2.5) -- (10,0.5);
  \draw (5,5.5) -- (8,2.5);
  \draw (6,6.5) -- (5,5.5);
  \draw (6,6.5) -- (11,0.5);
  \draw (6,6.5) -- (12,0.5);
  \draw (6,-6.5) -- (0,-0.5);
  \draw (6,-6.5) -- (1,-0.5);
  \draw (3,-1.5) -- (2,-0.5);
  \draw (3,-1.5) -- (3,-0.5);
  \draw (3,-1.5) -- (4,-0.5);
  \draw (7,-5.5) -- (3,-1.5);
  \draw (6,-1.5) -- (5,-0.5);
  \draw (6,-1.5) -- (6,-0.5);
  \draw (6,-1.5) -- (7,-0.5);
  \draw (7,-5.5) -- (6,-1.5);
  \draw (10,-2.5) -- (8,-0.5);
  \draw (10,-1.5) -- (9,-0.5);
  \draw (10,-1.5) -- (10,-0.5);
  \draw (10,-1.5) -- (11,-0.5);
  \draw (10,-2.5) -- (10,-1.5);
  \draw (10,-2.5) -- (12,-0.5);
  \draw (7,-5.5) -- (10,-2.5);
  \draw (6,-6.5) -- (7,-5.5);
  \fill (0,0.5) circle (2.5pt) (0,-0.5) circle (2.5pt);
  \fill (1,0.5) circle (2.5pt) (1,-0.5) circle (2.5pt);
  \fill (2,0.5) circle (2.5pt) (2,-0.5) circle (2.5pt);
  \fill (3,0.5) circle (2.5pt) (3,-0.5) circle (2.5pt);
  \fill (4,0.5) circle (2.5pt) (4,-0.5) circle (2.5pt);
  \fill (5,0.5) circle (2.5pt) (5,-0.5) circle (2.5pt);
  \fill (6,0.5) circle (2.5pt) (6,-0.5) circle (2.5pt);
  \fill (7,0.5) circle (2.5pt) (7,-0.5) circle (2.5pt);
  \fill (8,0.5) circle (2.5pt) (8,-0.5) circle (2.5pt);
  \fill (9,0.5) circle (2.5pt) (9,-0.5) circle (2.5pt);
  \fill (10,0.5) circle (2.5pt) (10,-0.5) circle (2.5pt);
  \fill (11,0.5) circle (2.5pt) (11,-0.5) circle (2.5pt);
  \fill (12,0.5) circle (2.5pt) (12,-0.5) circle (2.5pt);
\end{tikzpicture}} &
\resizebox{0.40\linewidth}{!}{\begin{tikzpicture}[scale=0.32, line cap=round, line join=round, >=stealth, thick]
  \draw[rounded corners=4pt] (5.682,5.682) -- (5.318,5.318);
  \draw[rounded corners=4pt] (6.288,5.654) -- (10.712,0.346) .. controls (11,0) .. (10.682,-0.318) -- (10.318,-0.682);
  \draw[rounded corners=4pt] (6.318,5.682) -- (11.682,0.318) .. controls (12,0) .. (11.682,-0.318) -- (10.318,-1.682);
  \draw[rounded corners=4pt] (4.682,4.682) -- (0.318,0.318) .. controls (0,0) .. (0.318,-0.318) -- (5.682,-5.682);
  \draw[rounded corners=4pt] (4.75,4.626) -- (3.25,2.374);
  \draw[rounded corners=4pt] (5.318,4.682) -- (7.682,2.318);
  \draw[rounded corners=4pt] (2.682,1.682) -- (2.318,1.318);
  \draw[rounded corners=4pt] (3.201,1.598) -- (3.799,0.402) .. controls (4,0) .. (3.682,-0.318) -- (3.318,-0.682);
  \draw[rounded corners=4pt] (3.318,1.682) -- (4.682,0.318) .. controls (5,0) .. (5.318,-0.318) -- (5.682,-0.682);
  \draw[rounded corners=4pt] (1.682,0.682) -- (1.318,0.318) .. controls (1,0) .. (1.288,-0.346) -- (5.712,-5.654);
  \draw[rounded corners=4pt] (2,0.55) -- (2,0.45) .. controls (2,0) .. (2.318,-0.318) -- (2.682,-0.682);
  \draw[rounded corners=4pt] (2.318,0.682) -- (2.682,0.318) .. controls (3,0) .. (3,-0.45) -- (3,-0.55);
  \draw[rounded corners=4pt] (7.682,1.682) -- (7.318,1.318);
  \draw[rounded corners=4pt] (8.201,1.598) -- (8.799,0.402) .. controls (9,0) .. (9.318,-0.318) -- (9.682,-0.682);
  \draw[rounded corners=4pt] (8.318,1.682) -- (9.682,0.318) .. controls (10,0) .. (10,-0.45) -- (10,-0.55);
  \draw[rounded corners=4pt] (6.682,0.682) -- (6.318,0.318) .. controls (6,0) .. (6,-0.45) -- (6,-0.55);
  \draw[rounded corners=4pt] (7,0.55) -- (7,0.45) .. controls (7,0) .. (6.682,-0.318) -- (6.318,-0.682);
  \draw[rounded corners=4pt] (7.318,0.682) -- (7.682,0.318) .. controls (8,0) .. (8.318,-0.318) -- (9.682,-1.682);
  \draw[rounded corners=4pt] (6.318,-5.682) -- (6.682,-5.318);
  \draw[rounded corners=4pt] (6.682,-4.682) -- (3.318,-1.318);
  \draw[rounded corners=4pt] (6.891,-4.563) -- (6.109,-1.437);
  \draw[rounded corners=4pt] (7.318,-4.682) -- (9.682,-2.318);
  \draw[rounded corners=4pt] (10,-1.55) -- (10,-1.45);
  \draw[rounded corners=4pt] (6,6.45) -- (6,7) -- (-1,7) -- (-1,-7) -- (6,-7) -- (6,-6.45);
  \draw (6,6.45) .. controls (6,6) .. (6.288,5.654);
  \draw (5.318,5.318) .. controls (5,5) .. (4.75,4.626);
  \draw (3.25,2.374) .. controls (3,2) .. (3.201,1.598);
  \draw (2.318,1.318) .. controls (2,1) .. (2,0.55);
  \draw (7.682,2.318) .. controls (8,2) .. (8.201,1.598);
  \draw (7.318,1.318) .. controls (7,1) .. (7,0.55);
  \draw (6,-6.45) .. controls (6,-6) .. (5.712,-5.654);
  \draw (6.682,-5.318) .. controls (7,-5) .. (6.891,-4.563);
  \draw (3.318,-1.318) .. controls (3,-1) .. (3,-0.55);
  \draw (6.109,-1.437) .. controls (6,-1) .. (6,-0.55);
  \draw (9.682,-2.318) .. controls (10,-2) .. (10,-1.55);
  \draw (10,-1.45) .. controls (10,-1) .. (10,-0.55);
  \draw[white, line width=6pt, line cap=butt] (5.682,5.682) .. controls (6,6) .. (6.318,5.682);
  \draw (5.682,5.682) .. controls (6,6) .. (6.318,5.682);
  \draw[white, line width=6pt, line cap=butt] (4.682,4.682) .. controls (5,5) .. (5.318,4.682);
  \draw (4.682,4.682) .. controls (5,5) .. (5.318,4.682);
  \draw[white, line width=6pt, line cap=butt] (2.682,1.682) .. controls (3,2) .. (3.318,1.682);
  \draw (2.682,1.682) .. controls (3,2) .. (3.318,1.682);
  \draw[white, line width=6pt, line cap=butt] (1.682,0.682) .. controls (2,1) .. (2.318,0.682);
  \draw (1.682,0.682) .. controls (2,1) .. (2.318,0.682);
  \draw[white, line width=6pt, line cap=butt] (7.682,1.682) .. controls (8,2) .. (8.318,1.682);
  \draw (7.682,1.682) .. controls (8,2) .. (8.318,1.682);
  \draw[white, line width=6pt, line cap=butt] (6.682,0.682) .. controls (7,1) .. (7.318,0.682);
  \draw (6.682,0.682) .. controls (7,1) .. (7.318,0.682);
  \draw[white, line width=6pt, line cap=butt] (5.682,-5.682) .. controls (6,-6) .. (6.318,-5.682);
  \draw (5.682,-5.682) .. controls (6,-6) .. (6.318,-5.682);
  \draw[white, line width=6pt, line cap=butt] (6.682,-4.682) .. controls (7,-5) .. (7.318,-4.682);
  \draw (6.682,-4.682) .. controls (7,-5) .. (7.318,-4.682);
  \draw[white, line width=6pt, line cap=butt] (2.682,-0.682) .. controls (3,-1) .. (3.318,-0.682);
  \draw (2.682,-0.682) .. controls (3,-1) .. (3.318,-0.682);
  \draw[white, line width=6pt, line cap=butt] (5.682,-0.682) .. controls (6,-1) .. (6.318,-0.682);
  \draw (5.682,-0.682) .. controls (6,-1) .. (6.318,-0.682);
  \draw[white, line width=6pt, line cap=butt] (9.682,-1.682) .. controls (10,-2) .. (10.318,-1.682);
  \draw (9.682,-1.682) .. controls (10,-2) .. (10.318,-1.682);
  \draw[white, line width=6pt, line cap=butt] (9.682,-0.682) .. controls (10,-1) .. (10.318,-0.682);
  \draw (9.682,-0.682) .. controls (10,-1) .. (10.318,-0.682);
  \draw[->] (5.5,5.5) -- (5.535,5.535);
  \draw[->] (8.5,3) -- (8.468,3.038);
  \draw[->] (9,3) -- (9.035,2.965);
  \draw[->] (3,-3) -- (2.965,-2.965);
  \draw[->] (4,3.5) -- (4.028,3.542);
  \draw[->] (6.5,3.5) -- (6.535,3.465);
  \draw[->] (2.5,1.5) -- (2.465,1.465);
  \draw[->] (3.5,1) -- (3.478,1.045);
  \draw[->] (4,1) -- (3.965,1.035);
  \draw[->] (3.5,-3) -- (3.468,-2.962);
  \draw[->] (2.5,-0.5) -- (2.535,-0.535);
  \draw[->] (2.5,0.5) -- (2.535,0.465);
  \draw[->] (7.5,1.5) -- (7.465,1.465);
  \draw[->] (8.5,1) -- (8.522,0.955);
  \draw[->] (9,1) -- (8.965,1.035);
  \draw[->] (6.5,0.5) -- (6.465,0.465);
  \draw[->] (6.5,-0.5) -- (6.465,-0.535);
  \draw[->] (9,-1) -- (8.965,-0.965);
  \draw[->] (6.5,-5.5) -- (6.465,-5.535);
  \draw[->] (5,-3) -- (5.035,-3.035);
  \draw[->] (6.5,-3) -- (6.512,-3.049);
  \draw[->] (8.5,-3.5) -- (8.535,-3.465);
  \draw[->] (10,-1.5) -- (10,-1.45);
  \draw[->] (-1,0) -- (-1,-0.05);
\end{tikzpicture}} \\[2mm]
\resizebox{0.32\linewidth}{!}{\begin{tikzpicture}[scale=0.32, line cap=round]
  \draw (6,6.5) -- (0,0.5);
  \draw (2,1.5) -- (1,0.5);
  \draw (2,1.5) -- (2,0.5);
  \draw (2,1.5) -- (3,0.5);
  \draw (5,4.5) -- (2,1.5);
  \draw (5,4.5) -- (4,0.5);
  \draw (6,1.5) -- (5,0.5);
  \draw (6,1.5) -- (6,0.5);
  \draw (6,1.5) -- (7,0.5);
  \draw (7,2.5) -- (6,1.5);
  \draw (7,2.5) -- (8,0.5);
  \draw (7,2.5) -- (9,0.5);
  \draw (5,4.5) -- (7,2.5);
  \draw (6,6.5) -- (5,4.5);
  \draw (11,1.5) -- (10,0.5);
  \draw (11,1.5) -- (11,0.5);
  \draw (11,1.5) -- (12,0.5);
  \draw (6,6.5) -- (11,1.5);
  \draw (8,8.5) -- (6,6.5);
  \draw (8,8.5) -- (13,0.5);
  \draw (15,1.5) -- (14,0.5);
  \draw (15,1.5) -- (15,0.5);
  \draw (15,1.5) -- (16,0.5);
  \draw (8,8.5) -- (15,1.5);
  \draw (8,-8.5) -- (0,-0.5);
  \draw (8,-8.5) -- (1,-0.5);
  \draw (9,-7.5) -- (2,-0.5);
  \draw (9,-7.5) -- (3,-0.5);
  \draw (10,-6.5) -- (4,-0.5);
  \draw (10,-6.5) -- (5,-0.5);
  \draw (9,-3.5) -- (6,-0.5);
  \draw (8,-1.5) -- (7,-0.5);
  \draw (8,-1.5) -- (8,-0.5);
  \draw (8,-1.5) -- (9,-0.5);
  \draw (9,-2.5) -- (8,-1.5);
  \draw (9,-2.5) -- (10,-0.5);
  \draw (9,-2.5) -- (11,-0.5);
  \draw (9,-3.5) -- (9,-2.5);
  \draw (9,-3.5) -- (12,-0.5);
  \draw (11,-5.5) -- (9,-3.5);
  \draw (14,-1.5) -- (13,-0.5);
  \draw (14,-1.5) -- (14,-0.5);
  \draw (14,-1.5) -- (15,-0.5);
  \draw (11,-5.5) -- (14,-1.5);
  \draw (11,-5.5) -- (16,-0.5);
  \draw (10,-6.5) -- (11,-5.5);
  \draw (9,-7.5) -- (10,-6.5);
  \draw (8,-8.5) -- (9,-7.5);
  \fill (0,0.5) circle (2.5pt) (0,-0.5) circle (2.5pt);
  \fill (1,0.5) circle (2.5pt) (1,-0.5) circle (2.5pt);
  \fill (2,0.5) circle (2.5pt) (2,-0.5) circle (2.5pt);
  \fill (3,0.5) circle (2.5pt) (3,-0.5) circle (2.5pt);
  \fill (4,0.5) circle (2.5pt) (4,-0.5) circle (2.5pt);
  \fill (5,0.5) circle (2.5pt) (5,-0.5) circle (2.5pt);
  \fill (6,0.5) circle (2.5pt) (6,-0.5) circle (2.5pt);
  \fill (7,0.5) circle (2.5pt) (7,-0.5) circle (2.5pt);
  \fill (8,0.5) circle (2.5pt) (8,-0.5) circle (2.5pt);
  \fill (9,0.5) circle (2.5pt) (9,-0.5) circle (2.5pt);
  \fill (10,0.5) circle (2.5pt) (10,-0.5) circle (2.5pt);
  \fill (11,0.5) circle (2.5pt) (11,-0.5) circle (2.5pt);
  \fill (12,0.5) circle (2.5pt) (12,-0.5) circle (2.5pt);
  \fill (13,0.5) circle (2.5pt) (13,-0.5) circle (2.5pt);
  \fill (14,0.5) circle (2.5pt) (14,-0.5) circle (2.5pt);
  \fill (15,0.5) circle (2.5pt) (15,-0.5) circle (2.5pt);
  \fill (16,0.5) circle (2.5pt) (16,-0.5) circle (2.5pt);
\end{tikzpicture}} &
\resizebox{0.40\linewidth}{!}{\begin{tikzpicture}[scale=0.32, line cap=round, line join=round, >=stealth, thick]
  \draw[rounded corners=4pt] (7.682,7.682) -- (6.318,6.318);
  \draw[rounded corners=4pt] (8.238,7.618) -- (12.762,0.382) .. controls (13,0) .. (13.318,-0.318) -- (13.682,-0.682);
  \draw[rounded corners=4pt] (8.318,7.682) -- (14.682,1.318);
  \draw[rounded corners=4pt] (5.682,5.682) -- (0.318,0.318) .. controls (0,0) .. (0.318,-0.318) -- (7.682,-7.682);
  \draw[rounded corners=4pt] (5.799,5.598) -- (5.201,4.402);
  \draw[rounded corners=4pt] (6.318,5.682) -- (10.682,1.318);
  \draw[rounded corners=4pt] (4.682,3.682) -- (2.318,1.318);
  \draw[rounded corners=4pt] (4.891,3.563) -- (4.109,0.437) .. controls (4,0) .. (4.318,-0.318) -- (9.682,-5.682);
  \draw[rounded corners=4pt] (5.318,3.682) -- (6.682,2.318);
  \draw[rounded corners=4pt] (1.682,0.682) -- (1.318,0.318) .. controls (1,0) .. (1.296,-0.339) -- (7.704,-7.661);
  \draw[rounded corners=4pt] (2,0.55) -- (2,0.45) .. controls (2,0) .. (2.318,-0.318) -- (8.682,-6.682);
  \draw[rounded corners=4pt] (2.318,0.682) -- (2.682,0.318) .. controls (3,0) .. (3.293,-0.342) -- (8.707,-6.658);
  \draw[rounded corners=4pt] (6.682,1.682) -- (6.318,1.318);
  \draw[rounded corners=4pt] (7.201,1.598) -- (7.799,0.402) .. controls (8,0) .. (8,-0.45) -- (8,-0.55);
  \draw[rounded corners=4pt] (7.318,1.682) -- (8.682,0.318) .. controls (9,0) .. (8.682,-0.318) -- (8.318,-0.682);
  \draw[rounded corners=4pt] (5.682,0.682) -- (5.318,0.318) .. controls (5,0) .. (5.288,-0.346) -- (9.712,-5.654);
  \draw[rounded corners=4pt] (6,0.55) -- (6,0.45) .. controls (6,0) .. (6.318,-0.318) -- (8.682,-2.682);
  \draw[rounded corners=4pt] (6.318,0.682) -- (6.682,0.318) .. controls (7,0) .. (7.318,-0.318) -- (7.682,-0.682);
  \draw[rounded corners=4pt] (10.682,0.682) -- (10.318,0.318) .. controls (10,0) .. (9.799,-0.402) -- (9.201,-1.598);
  \draw[rounded corners=4pt] (11,0.55) -- (11,0.45) .. controls (11,0) .. (10.682,-0.318) -- (9.318,-1.682);
  \draw[rounded corners=4pt] (11.318,0.682) -- (11.682,0.318) .. controls (12,0) .. (11.682,-0.318) -- (9.318,-2.682);
  \draw[rounded corners=4pt] (14.682,0.682) -- (14.318,0.318) .. controls (14,0) .. (14,-0.45) -- (14,-0.55);
  \draw[rounded corners=4pt] (15,0.55) -- (15,0.45) .. controls (15,0) .. (14.682,-0.318) -- (14.318,-0.682);
  \draw[rounded corners=4pt] (15.318,0.682) -- (15.682,0.318) .. controls (16,0) .. (15.682,-0.318) -- (11.318,-4.682);
  \draw[rounded corners=4pt] (8.318,-7.682) -- (8.682,-7.318);
  \draw[rounded corners=4pt] (9.318,-6.682) -- (9.682,-6.318);
  \draw[rounded corners=4pt] (10.318,-5.682) -- (10.682,-5.318);
  \draw[rounded corners=4pt] (10.682,-4.682) -- (9.318,-3.318);
  \draw[rounded corners=4pt] (11.27,-4.64) -- (13.73,-1.36);
  \draw[rounded corners=4pt] (9,-2.55) -- (9,-2.45);
  \draw[rounded corners=4pt] (8.682,-1.682) -- (8.318,-1.318);
  \draw[rounded corners=4pt] (8,8.45) -- (8,9) -- (-1,9) -- (-1,-9) -- (8,-9) -- (8,-8.45);
  \draw (8,8.45) .. controls (8,8) .. (8.238,7.618);
  \draw (6.318,6.318) .. controls (6,6) .. (5.799,5.598);
  \draw (5.201,4.402) .. controls (5,4) .. (4.891,3.563);
  \draw (2.318,1.318) .. controls (2,1) .. (2,0.55);
  \draw (6.682,2.318) .. controls (7,2) .. (7.201,1.598);
  \draw (6.318,1.318) .. controls (6,1) .. (6,0.55);
  \draw (10.682,1.318) .. controls (11,1) .. (11,0.55);
  \draw (14.682,1.318) .. controls (15,1) .. (15,0.55);
  \draw (8,-8.45) .. controls (8,-8) .. (7.704,-7.661);
  \draw (8.682,-7.318) .. controls (9,-7) .. (8.707,-6.658);
  \draw (9.682,-6.318) .. controls (10,-6) .. (9.712,-5.654);
  \draw (10.682,-5.318) .. controls (11,-5) .. (11.27,-4.64);
  \draw (9.318,-3.318) .. controls (9,-3) .. (9,-2.55);
  \draw (9,-2.45) .. controls (9,-2) .. (9.201,-1.598);
  \draw (8.318,-1.318) .. controls (8,-1) .. (8,-0.55);
  \draw (13.73,-1.36) .. controls (14,-1) .. (14,-0.55);
  \draw[white, line width=6pt, line cap=butt] (7.682,7.682) .. controls (8,8) .. (8.318,7.682);
  \draw (7.682,7.682) .. controls (8,8) .. (8.318,7.682);
  \draw[white, line width=6pt, line cap=butt] (5.682,5.682) .. controls (6,6) .. (6.318,5.682);
  \draw (5.682,5.682) .. controls (6,6) .. (6.318,5.682);
  \draw[white, line width=6pt, line cap=butt] (4.682,3.682) .. controls (5,4) .. (5.318,3.682);
  \draw (4.682,3.682) .. controls (5,4) .. (5.318,3.682);
  \draw[white, line width=6pt, line cap=butt] (1.682,0.682) .. controls (2,1) .. (2.318,0.682);
  \draw (1.682,0.682) .. controls (2,1) .. (2.318,0.682);
  \draw[white, line width=6pt, line cap=butt] (6.682,1.682) .. controls (7,2) .. (7.318,1.682);
  \draw (6.682,1.682) .. controls (7,2) .. (7.318,1.682);
  \draw[white, line width=6pt, line cap=butt] (5.682,0.682) .. controls (6,1) .. (6.318,0.682);
  \draw (5.682,0.682) .. controls (6,1) .. (6.318,0.682);
  \draw[white, line width=6pt, line cap=butt] (10.682,0.682) .. controls (11,1) .. (11.318,0.682);
  \draw (10.682,0.682) .. controls (11,1) .. (11.318,0.682);
  \draw[white, line width=6pt, line cap=butt] (14.682,0.682) .. controls (15,1) .. (15.318,0.682);
  \draw (14.682,0.682) .. controls (15,1) .. (15.318,0.682);
  \draw[white, line width=6pt, line cap=butt] (7.682,-7.682) .. controls (8,-8) .. (8.318,-7.682);
  \draw (7.682,-7.682) .. controls (8,-8) .. (8.318,-7.682);
  \draw[white, line width=6pt, line cap=butt] (8.682,-6.682) .. controls (9,-7) .. (9.318,-6.682);
  \draw (8.682,-6.682) .. controls (9,-7) .. (9.318,-6.682);
  \draw[white, line width=6pt, line cap=butt] (9.682,-5.682) .. controls (10,-6) .. (10.318,-5.682);
  \draw (9.682,-5.682) .. controls (10,-6) .. (10.318,-5.682);
  \draw[white, line width=6pt, line cap=butt] (10.682,-4.682) .. controls (11,-5) .. (11.318,-4.682);
  \draw (10.682,-4.682) .. controls (11,-5) .. (11.318,-4.682);
  \draw[white, line width=6pt, line cap=butt] (8.682,-2.682) .. controls (9,-3) .. (9.318,-2.682);
  \draw (8.682,-2.682) .. controls (9,-3) .. (9.318,-2.682);
  \draw[white, line width=6pt, line cap=butt] (8.682,-1.682) .. controls (9,-2) .. (9.318,-1.682);
  \draw (8.682,-1.682) .. controls (9,-2) .. (9.318,-1.682);
  \draw[white, line width=6pt, line cap=butt] (7.682,-0.682) .. controls (8,-1) .. (8.318,-0.682);
  \draw (7.682,-0.682) .. controls (8,-1) .. (8.318,-0.682);
  \draw[white, line width=6pt, line cap=butt] (13.682,-0.682) .. controls (14,-1) .. (14.318,-0.682);
  \draw (13.682,-0.682) .. controls (14,-1) .. (14.318,-0.682);
  \draw[->] (7,7) -- (7.035,7.035);
  \draw[->] (10.5,4) -- (10.474,4.042);
  \draw[->] (11.5,4.5) -- (11.535,4.465);
  \draw[->] (4,-4) -- (3.965,-3.965);
  \draw[->] (5.5,5) -- (5.522,5.045);
  \draw[->] (8.5,3.5) -- (8.535,3.465);
  \draw[->] (3.5,2.5) -- (3.465,2.465);
  \draw[->] (7,-3) -- (6.965,-2.965);
  \draw[->] (6,3) -- (5.965,3.035);
  \draw[->] (4.5,-4) -- (4.467,-3.962);
  \draw[->] (5.5,-3.5) -- (5.535,-3.535);
  \draw[->] (6,-3.5) -- (6.033,-3.538);
  \draw[->] (6.5,1.5) -- (6.465,1.465);
  \draw[->] (7.5,1) -- (7.478,1.045);
  \draw[->] (8,1) -- (7.965,1.035);
  \draw[->] (7.5,-3) -- (7.468,-2.962);
  \draw[->] (7.5,-1.5) -- (7.535,-1.535);
  \draw[->] (6.5,0.5) -- (6.535,0.465);
  \draw[->] (9.5,-1) -- (9.478,-1.045);
  \draw[->] (10,-1) -- (9.965,-1.035);
  \draw[->] (10.5,-1.5) -- (10.535,-1.465);
  \draw[->] (14.5,0.5) -- (14.465,0.465);
  \draw[->] (14.5,-0.5) -- (14.465,-0.535);
  \draw[->] (13.5,-2.5) -- (13.535,-2.465);
  \draw[->] (8.5,-7.5) -- (8.465,-7.535);
  \draw[->] (9.5,-6.5) -- (9.535,-6.465);
  \draw[->] (10.5,-5.5) -- (10.465,-5.535);
  \draw[->] (10,-4) -- (10.035,-4.035);
  \draw[->] (12.5,-3) -- (12.47,-3.04);
  \draw[->] (9,-2.5) -- (9,-2.55);
  \draw[->] (8.5,-1.5) -- (8.465,-1.465);
  \draw[->] (-1,0) -- (-1,-0.05);
\end{tikzpicture}}
\end{tabular}
\caption{The two orientations of $8_{17}$ in $\vFt$. Top: $K$ with $n=7$. Bottom: $K^r$ with $n=9$; no
diagram of $\vFt$ with $n\le 7$ realizes $K^r$.}
\label{fig:817-F3vec}
\end{figure}
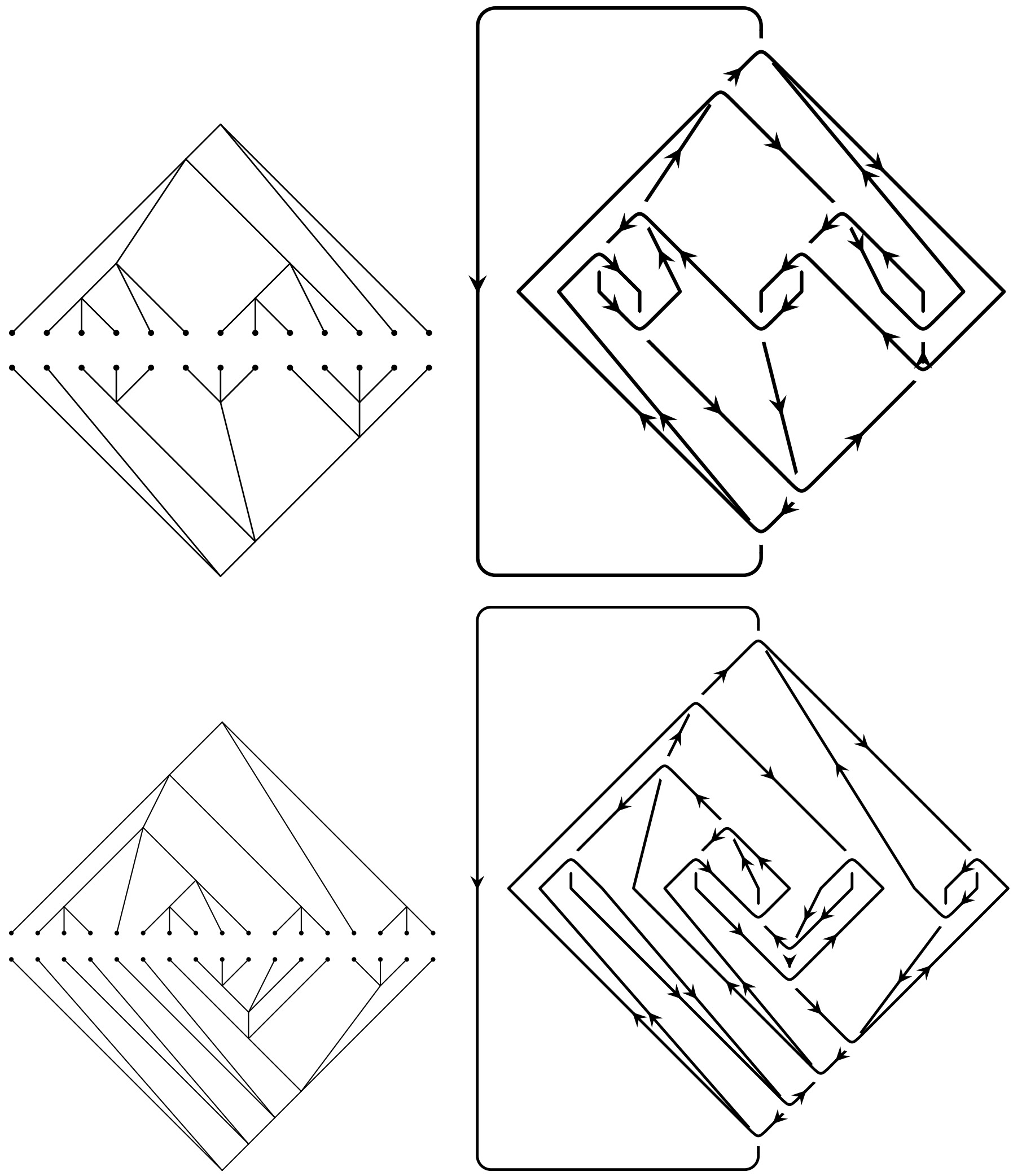

\begin{figure}[tbp]
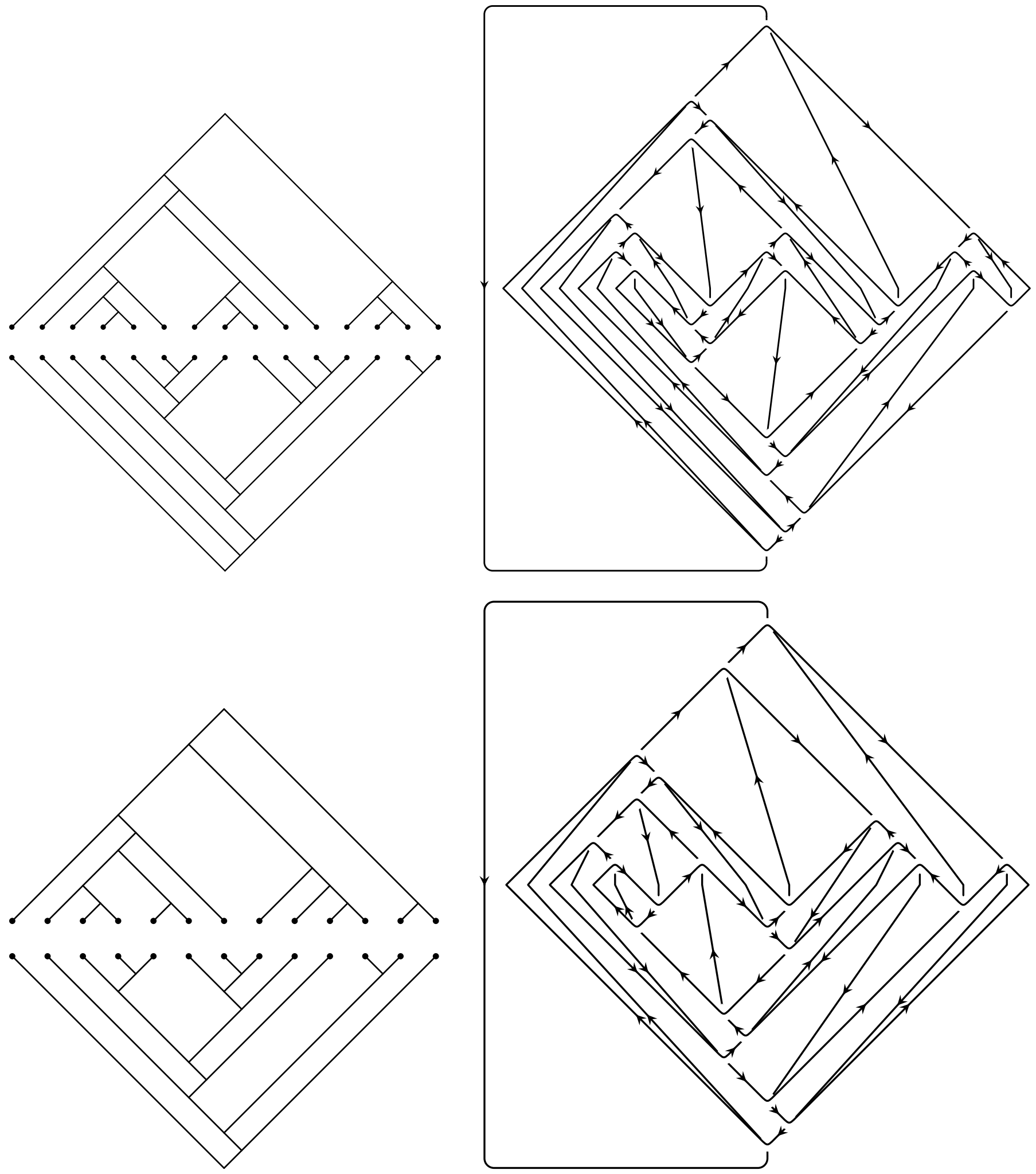

\centering
\begin{tabular}{cc}
\resizebox{0.32\linewidth}{!}{\begin{tikzpicture}[scale=0.32, line cap=round]
  \draw (3.5,1) -- (3,0.5);
  \draw (3.5,1) -- (4,0.5);
  \draw (3,1.5) -- (2,0.5);
  \draw (3,1.5) -- (3.5,1);
  \draw (3.5,2) -- (3,1.5);
  \draw (3.5,2) -- (5,0.5);
  \draw (3,2.5) -- (1,0.5);
  \draw (3,2.5) -- (3.5,2);
  \draw (7.5,1) -- (7,0.5);
  \draw (7.5,1) -- (8,0.5);
  \draw (7,1.5) -- (6,0.5);
  \draw (7,1.5) -- (7.5,1);
  \draw (7.5,2) -- (7,1.5);
  \draw (7.5,2) -- (9,0.5);
  \draw (5,4.5) -- (3,2.5);
  \draw (5,4.5) -- (7.5,2);
  \draw (5.5,5) -- (5,4.5);
  \draw (5.5,5) -- (10,0.5);
  \draw (5,5.5) -- (0,0.5);
  \draw (5,5.5) -- (5.5,5);
  \draw (12.5,1) -- (12,0.5);
  \draw (12.5,1) -- (13,0.5);
  \draw (12,1.5) -- (11,0.5);
  \draw (12,1.5) -- (12.5,1);
  \draw (12.5,2) -- (12,1.5);
  \draw (12.5,2) -- (14,0.5);
  \draw (7,7.5) -- (5,5.5);
  \draw (7,7.5) -- (12.5,2);
  \draw (5.5,-1) -- (5,-0.5);
  \draw (5.5,-1) -- (6,-0.5);
  \draw (5,-1.5) -- (4,-0.5);
  \draw (5,-1.5) -- (5.5,-1);
  \draw (5.5,-2) -- (5,-1.5);
  \draw (5.5,-2) -- (7,-0.5);
  \draw (5,-2.5) -- (3,-0.5);
  \draw (5,-2.5) -- (5.5,-2);
  \draw (10.5,-1) -- (10,-0.5);
  \draw (10.5,-1) -- (11,-0.5);
  \draw (10,-1.5) -- (9,-0.5);
  \draw (10,-1.5) -- (10.5,-1);
  \draw (9.5,-2) -- (8,-0.5);
  \draw (9.5,-2) -- (10,-1.5);
  \draw (7,-4.5) -- (5,-2.5);
  \draw (7,-4.5) -- (9.5,-2);
  \draw (7.5,-5) -- (7,-4.5);
  \draw (7.5,-5) -- (12,-0.5);
  \draw (7,-5.5) -- (2,-0.5);
  \draw (7,-5.5) -- (7.5,-5);
  \draw (13.5,-1) -- (13,-0.5);
  \draw (13.5,-1) -- (14,-0.5);
  \draw (8,-6.5) -- (7,-5.5);
  \draw (8,-6.5) -- (13.5,-1);
  \draw (7.5,-7) -- (1,-0.5);
  \draw (7.5,-7) -- (8,-6.5);
  \draw (7,-7.5) -- (0,-0.5);
  \draw (7,-7.5) -- (7.5,-7);
  \fill (0,0.5) circle (2.5pt) (0,-0.5) circle (2.5pt);
  \fill (1,0.5) circle (2.5pt) (1,-0.5) circle (2.5pt);
  \fill (2,0.5) circle (2.5pt) (2,-0.5) circle (2.5pt);
  \fill (3,0.5) circle (2.5pt) (3,-0.5) circle (2.5pt);
  \fill (4,0.5) circle (2.5pt) (4,-0.5) circle (2.5pt);
  \fill (5,0.5) circle (2.5pt) (5,-0.5) circle (2.5pt);
  \fill (6,0.5) circle (2.5pt) (6,-0.5) circle (2.5pt);
  \fill (7,0.5) circle (2.5pt) (7,-0.5) circle (2.5pt);
  \fill (8,0.5) circle (2.5pt) (8,-0.5) circle (2.5pt);
  \fill (9,0.5) circle (2.5pt) (9,-0.5) circle (2.5pt);
  \fill (10,0.5) circle (2.5pt) (10,-0.5) circle (2.5pt);
  \fill (11,0.5) circle (2.5pt) (11,-0.5) circle (2.5pt);
  \fill (12,0.5) circle (2.5pt) (12,-0.5) circle (2.5pt);
  \fill (13,0.5) circle (2.5pt) (13,-0.5) circle (2.5pt);
  \fill (14,0.5) circle (2.5pt) (14,-0.5) circle (2.5pt);
\end{tikzpicture}} &
\resizebox{0.40\linewidth}{!}{\input{figures/ex1088_Fvec_K_knot}} \\[2mm]
\resizebox{0.32\linewidth}{!}{\begin{tikzpicture}[scale=0.32, line cap=round]
  \draw (2.5,1) -- (2,0.5);
  \draw (2.5,1) -- (3,0.5);
  \draw (2,1.5) -- (1,0.5);
  \draw (2,1.5) -- (2.5,1);
  \draw (4.5,1) -- (4,0.5);
  \draw (4.5,1) -- (5,0.5);
  \draw (3,2.5) -- (2,1.5);
  \draw (3,2.5) -- (4.5,1);
  \draw (3.5,3) -- (3,2.5);
  \draw (3.5,3) -- (6,0.5);
  \draw (3,3.5) -- (0,0.5);
  \draw (3,3.5) -- (3.5,3);
  \draw (9.5,1) -- (9,0.5);
  \draw (9.5,1) -- (10,0.5);
  \draw (9,1.5) -- (8,0.5);
  \draw (9,1.5) -- (9.5,1);
  \draw (8.5,2) -- (7,0.5);
  \draw (8.5,2) -- (9,1.5);
  \draw (5,5.5) -- (3,3.5);
  \draw (5,5.5) -- (8.5,2);
  \draw (11.5,1) -- (11,0.5);
  \draw (11.5,1) -- (12,0.5);
  \draw (6,6.5) -- (5,5.5);
  \draw (6,6.5) -- (11.5,1);
  \draw (3.5,-1) -- (3,-0.5);
  \draw (3.5,-1) -- (4,-0.5);
  \draw (3,-1.5) -- (2,-0.5);
  \draw (3,-1.5) -- (3.5,-1);
  \draw (6.5,-1) -- (6,-0.5);
  \draw (6.5,-1) -- (7,-0.5);
  \draw (6,-1.5) -- (5,-0.5);
  \draw (6,-1.5) -- (6.5,-1);
  \draw (6.5,-2) -- (6,-1.5);
  \draw (6.5,-2) -- (8,-0.5);
  \draw (5,-3.5) -- (3,-1.5);
  \draw (5,-3.5) -- (6.5,-2);
  \draw (5.5,-4) -- (5,-3.5);
  \draw (5.5,-4) -- (9,-0.5);
  \draw (5,-4.5) -- (1,-0.5);
  \draw (5,-4.5) -- (5.5,-4);
  \draw (10.5,-1) -- (10,-0.5);
  \draw (10.5,-1) -- (11,-0.5);
  \draw (6,-5.5) -- (5,-4.5);
  \draw (6,-5.5) -- (10.5,-1);
  \draw (6.5,-6) -- (6,-5.5);
  \draw (6.5,-6) -- (12,-0.5);
  \draw (6,-6.5) -- (0,-0.5);
  \draw (6,-6.5) -- (6.5,-6);
  \fill (0,0.5) circle (2.5pt) (0,-0.5) circle (2.5pt);
  \fill (1,0.5) circle (2.5pt) (1,-0.5) circle (2.5pt);
  \fill (2,0.5) circle (2.5pt) (2,-0.5) circle (2.5pt);
  \fill (3,0.5) circle (2.5pt) (3,-0.5) circle (2.5pt);
  \fill (4,0.5) circle (2.5pt) (4,-0.5) circle (2.5pt);
  \fill (5,0.5) circle (2.5pt) (5,-0.5) circle (2.5pt);
  \fill (6,0.5) circle (2.5pt) (6,-0.5) circle (2.5pt);
  \fill (7,0.5) circle (2.5pt) (7,-0.5) circle (2.5pt);
  \fill (8,0.5) circle (2.5pt) (8,-0.5) circle (2.5pt);
  \fill (9,0.5) circle (2.5pt) (9,-0.5) circle (2.5pt);
  \fill (10,0.5) circle (2.5pt) (10,-0.5) circle (2.5pt);
  \fill (11,0.5) circle (2.5pt) (11,-0.5) circle (2.5pt);
  \fill (12,0.5) circle (2.5pt) (12,-0.5) circle (2.5pt);
\end{tikzpicture}} &
\resizebox{0.40\linewidth}{!}{\input{figures/ex1088_Fvec_Kr_knot}}
\end{tabular}
\caption{The two orientations of $10_{88}$ in $\vF$. Top: $K$ (KnotInfo orientation) with 15 leaves. Bottom:
$K^r$ with 13 leaves; no diagram of $\vF$ with at most 13 leaves realizes $K$.}
\label{fig:1088-Fvec}
\end{figure}

The knot $8_{17}$ is the non-invertible knot with the fewest crossings; it is negative amphicheiral, that is
$K=K^{mr}$ and $K^r=K^m$. Figure~\ref{fig:817-unoriented} shows minimal tree diagrams of $8_{17}$ in $F$
(9 leaves) and $F_3$ ($n=7$) and the corresponding knot diagrams. Figure~\ref{fig:817-F3vec} shows minimal
diagrams of the two orientations of $8_{17}$ in $\vFt$ ($n=7$ for $K$, $n=9$ for $K^r$): \emph{the
$\vFt$-index distinguishes $8_{17}$ from its reverse}. Its $\vF$-index is instead $13$ in both orientations.
The knots with fewest crossings whose two orientations have different $\vF$-index have 10 crossings;
Figure~\ref{fig:1088-Fvec} shows minimal diagrams of the two orientations of $10_{88}$ in $\vF$, with 15
leaves for $K$ and 13 for $K^r$.

\begin{figure}[tbp]
\centering
\begin{tabular}{cc}
\resizebox{0.32\linewidth}{!}{\begin{tikzpicture}[scale=0.32, line cap=round]
  \draw (1.5,1) -- (1,0.5);
  \draw (1.5,1) -- (2,0.5);
  \draw (3.5,1) -- (3,0.5);
  \draw (3.5,1) -- (4,0.5);
  \draw (6.5,1) -- (6,0.5);
  \draw (6.5,1) -- (7,0.5);
  \draw (6,1.5) -- (5,0.5);
  \draw (6,1.5) -- (6.5,1);
  \draw (6.5,2) -- (6,1.5);
  \draw (6.5,2) -- (8,0.5);
  \draw (5.5,3) -- (3.5,1);
  \draw (5.5,3) -- (6.5,2);
  \draw (4.5,4) -- (1.5,1);
  \draw (4.5,4) -- (5.5,3);
  \draw (4,4.5) -- (0,0.5);
  \draw (4,4.5) -- (4.5,4);
  \draw (0.5,-1) -- (0,-0.5);
  \draw (0.5,-1) -- (1,-0.5);
  \draw (2.5,-1) -- (2,-0.5);
  \draw (2.5,-1) -- (3,-0.5);
  \draw (4.5,-1) -- (4,-0.5);
  \draw (4.5,-1) -- (5,-0.5);
  \draw (5,-1.5) -- (4.5,-1);
  \draw (5,-1.5) -- (6,-0.5);
  \draw (4,-2.5) -- (2.5,-1);
  \draw (4,-2.5) -- (5,-1.5);
  \draw (3,-3.5) -- (0.5,-1);
  \draw (3,-3.5) -- (4,-2.5);
  \draw (7.5,-1) -- (7,-0.5);
  \draw (7.5,-1) -- (8,-0.5);
  \draw (4,-4.5) -- (3,-3.5);
  \draw (4,-4.5) -- (7.5,-1);
  \fill (0,0.5) circle (2.5pt) (0,-0.5) circle (2.5pt);
  \fill (1,0.5) circle (2.5pt) (1,-0.5) circle (2.5pt);
  \fill (2,0.5) circle (2.5pt) (2,-0.5) circle (2.5pt);
  \fill (3,0.5) circle (2.5pt) (3,-0.5) circle (2.5pt);
  \fill (4,0.5) circle (2.5pt) (4,-0.5) circle (2.5pt);
  \fill (5,0.5) circle (2.5pt) (5,-0.5) circle (2.5pt);
  \fill (6,0.5) circle (2.5pt) (6,-0.5) circle (2.5pt);
  \fill (7,0.5) circle (2.5pt) (7,-0.5) circle (2.5pt);
  \fill (8,0.5) circle (2.5pt) (8,-0.5) circle (2.5pt);
\end{tikzpicture}} &
\resizebox{0.40\linewidth}{!}{\begin{tikzpicture}[scale=0.32, line cap=round, line join=round, >=stealth, thick]
  \draw[rounded corners=4pt] (7.682,7.682) -- (0.318,0.318) .. controls (0,0) .. (0.318,-0.318) -- (0.682,-0.682);
  \draw[rounded corners=4pt] (7.704,7.661) -- (1.296,0.339) .. controls (1,0) .. (1,-0.45) -- (1,-0.55);
  \draw[rounded corners=4pt] (8.318,7.682) -- (8.682,7.318);
  \draw[rounded corners=4pt] (8.682,6.682) -- (3.318,1.318);
  \draw[rounded corners=4pt] (8.777,6.609) -- (5.223,0.391) .. controls (5,0) .. (5,-0.45) -- (5,-0.55);
  \draw[rounded corners=4pt] (9.318,6.682) -- (10.682,5.318);
  \draw[rounded corners=4pt] (2.682,0.682) -- (2.318,0.318) .. controls (2,0) .. (1.682,-0.318) -- (1.318,-0.682);
  \draw[rounded corners=4pt] (3,0.55) -- (3,0.45) .. controls (3,0) .. (3.201,-0.402) -- (5.799,-5.598);
  \draw[rounded corners=4pt] (3.318,0.682) -- (3.682,0.318) .. controls (4,0) .. (4.318,-0.318) -- (4.682,-0.682);
  \draw[rounded corners=4pt] (10.682,4.682) -- (7.318,1.318);
  \draw[rounded corners=4pt] (10.833,4.582) -- (9.167,0.418) .. controls (9,0) .. (9,-0.45) -- (9,-0.55);
  \draw[rounded corners=4pt] (11.318,4.682) -- (12.682,3.318);
  \draw[rounded corners=4pt] (6.682,0.682) -- (6.318,0.318) .. controls (6,0) .. (5.682,-0.318) -- (5.318,-0.682);
  \draw[rounded corners=4pt] (7,0.55) -- (7,0.45) .. controls (7,0) .. (7.109,-0.437) -- (7.891,-3.563);
  \draw[rounded corners=4pt] (7.318,0.682) -- (7.682,0.318) .. controls (8,0) .. (8.318,-0.318) -- (8.682,-0.682);
  \draw[rounded corners=4pt] (12.682,2.682) -- (12.318,2.318);
  \draw[rounded corners=4pt] (13.25,2.626) -- (14.75,0.374) .. controls (15,0) .. (15,-0.45) -- (15,-0.55);
  \draw[rounded corners=4pt] (13.318,2.682) -- (15.682,0.318) .. controls (16,0) .. (15.682,-0.318) -- (15.318,-0.682);
  \draw[rounded corners=4pt] (11.682,1.682) -- (10.318,0.318) .. controls (10,0) .. (9.682,-0.318) -- (9.318,-0.682);
  \draw[rounded corners=4pt] (11.799,1.598) -- (11.201,0.402) .. controls (11,0) .. (10.799,-0.402) -- (10.201,-1.598);
  \draw[rounded corners=4pt] (12.318,1.682) -- (12.682,1.318);
  \draw[rounded corners=4pt] (12.682,0.682) -- (12.318,0.318) .. controls (12,0) .. (11.682,-0.318) -- (10.318,-1.682);
  \draw[rounded corners=4pt] (13,0.55) -- (13,0.45) .. controls (13,0) .. (12.762,-0.382) -- (8.238,-7.618);
  \draw[rounded corners=4pt] (13.318,0.682) -- (13.682,0.318) .. controls (14,0) .. (14.318,-0.318) -- (14.682,-0.682);
  \draw[rounded corners=4pt] (7.682,-7.682) -- (6.318,-6.318);
  \draw[rounded corners=4pt] (8.318,-7.682) -- (14.682,-1.318);
  \draw[rounded corners=4pt] (5.682,-5.682) -- (1.318,-1.318);
  \draw[rounded corners=4pt] (6.318,-5.682) -- (7.682,-4.318);
  \draw[rounded corners=4pt] (7.682,-3.682) -- (5.318,-1.318);
  \draw[rounded corners=4pt] (8.318,-3.682) -- (9.682,-2.318);
  \draw[rounded corners=4pt] (9.682,-1.682) -- (9.318,-1.318);
  \draw[rounded corners=4pt] (8,8.45) -- (8,9) -- (-1,9) -- (-1,-9) -- (8,-9) -- (8,-8.45);
  \draw (8,8.45) .. controls (8,8) .. (7.704,7.661);
  \draw (8.682,7.318) .. controls (9,7) .. (8.777,6.609);
  \draw (3.318,1.318) .. controls (3,1) .. (3,0.55);
  \draw (10.682,5.318) .. controls (11,5) .. (10.833,4.582);
  \draw (7.318,1.318) .. controls (7,1) .. (7,0.55);
  \draw (12.682,3.318) .. controls (13,3) .. (13.25,2.626);
  \draw (12.318,2.318) .. controls (12,2) .. (11.799,1.598);
  \draw (12.682,1.318) .. controls (13,1) .. (13,0.55);
  \draw (8,-8.45) .. controls (8,-8) .. (8.238,-7.618);
  \draw (6.318,-6.318) .. controls (6,-6) .. (5.799,-5.598);
  \draw (1.318,-1.318) .. controls (1,-1) .. (1,-0.55);
  \draw (7.682,-4.318) .. controls (8,-4) .. (7.891,-3.563);
  \draw (5.318,-1.318) .. controls (5,-1) .. (5,-0.55);
  \draw (9.682,-2.318) .. controls (10,-2) .. (10.201,-1.598);
  \draw (9.318,-1.318) .. controls (9,-1) .. (9,-0.55);
  \draw (14.682,-1.318) .. controls (15,-1) .. (15,-0.55);
  \draw[white, line width=6pt, line cap=butt] (7.682,7.682) .. controls (8,8) .. (8.318,7.682);
  \draw (7.682,7.682) .. controls (8,8) .. (8.318,7.682);
  \draw[white, line width=6pt, line cap=butt] (8.682,6.682) .. controls (9,7) .. (9.318,6.682);
  \draw (8.682,6.682) .. controls (9,7) .. (9.318,6.682);
  \draw[white, line width=6pt, line cap=butt] (2.682,0.682) .. controls (3,1) .. (3.318,0.682);
  \draw (2.682,0.682) .. controls (3,1) .. (3.318,0.682);
  \draw[white, line width=6pt, line cap=butt] (10.682,4.682) .. controls (11,5) .. (11.318,4.682);
  \draw (10.682,4.682) .. controls (11,5) .. (11.318,4.682);
  \draw[white, line width=6pt, line cap=butt] (6.682,0.682) .. controls (7,1) .. (7.318,0.682);
  \draw (6.682,0.682) .. controls (7,1) .. (7.318,0.682);
  \draw[white, line width=6pt, line cap=butt] (12.682,2.682) .. controls (13,3) .. (13.318,2.682);
  \draw (12.682,2.682) .. controls (13,3) .. (13.318,2.682);
  \draw[white, line width=6pt, line cap=butt] (11.682,1.682) .. controls (12,2) .. (12.318,1.682);
  \draw (11.682,1.682) .. controls (12,2) .. (12.318,1.682);
  \draw[white, line width=6pt, line cap=butt] (12.682,0.682) .. controls (13,1) .. (13.318,0.682);
  \draw (12.682,0.682) .. controls (13,1) .. (13.318,0.682);
  \draw[white, line width=6pt, line cap=butt] (7.682,-7.682) .. controls (8,-8) .. (8.318,-7.682);
  \draw (7.682,-7.682) .. controls (8,-8) .. (8.318,-7.682);
  \draw[white, line width=6pt, line cap=butt] (5.682,-5.682) .. controls (6,-6) .. (6.318,-5.682);
  \draw (5.682,-5.682) .. controls (6,-6) .. (6.318,-5.682);
  \draw[white, line width=6pt, line cap=butt] (0.682,-0.682) .. controls (1,-1) .. (1.318,-0.682);
  \draw (0.682,-0.682) .. controls (1,-1) .. (1.318,-0.682);
  \draw[white, line width=6pt, line cap=butt] (7.682,-3.682) .. controls (8,-4) .. (8.318,-3.682);
  \draw (7.682,-3.682) .. controls (8,-4) .. (8.318,-3.682);
  \draw[white, line width=6pt, line cap=butt] (4.682,-0.682) .. controls (5,-1) .. (5.318,-0.682);
  \draw (4.682,-0.682) .. controls (5,-1) .. (5.318,-0.682);
  \draw[white, line width=6pt, line cap=butt] (9.682,-1.682) .. controls (10,-2) .. (10.318,-1.682);
  \draw (9.682,-1.682) .. controls (10,-2) .. (10.318,-1.682);
  \draw[white, line width=6pt, line cap=butt] (8.682,-0.682) .. controls (9,-1) .. (9.318,-0.682);
  \draw (8.682,-0.682) .. controls (9,-1) .. (9.318,-0.682);
  \draw[white, line width=6pt, line cap=butt] (14.682,-0.682) .. controls (15,-1) .. (15.318,-0.682);
  \draw (14.682,-0.682) .. controls (15,-1) .. (15.318,-0.682);
\end{tikzpicture}} \\[2mm]
\resizebox{0.32\linewidth}{!}{\begin{tikzpicture}[scale=0.32, line cap=round]
  \draw (6,6.5) -- (0,0.5);
  \draw (6,6.5) -- (1,0.5);
  \draw (7,5.5) -- (2,0.5);
  \draw (7,5.5) -- (3,0.5);
  \draw (8,4.5) -- (4,0.5);
  \draw (8,4.5) -- (5,0.5);
  \draw (7,1.5) -- (6,0.5);
  \draw (7,1.5) -- (7,0.5);
  \draw (7,1.5) -- (8,0.5);
  \draw (9,3.5) -- (7,1.5);
  \draw (10,1.5) -- (9,0.5);
  \draw (10,1.5) -- (10,0.5);
  \draw (10,1.5) -- (11,0.5);
  \draw (9,3.5) -- (10,1.5);
  \draw (9,3.5) -- (12,0.5);
  \draw (8,4.5) -- (9,3.5);
  \draw (7,5.5) -- (8,4.5);
  \draw (6,6.5) -- (7,5.5);
  \draw (6,-6.5) -- (0,-0.5);
  \draw (5,-4.5) -- (1,-0.5);
  \draw (5,-4.5) -- (2,-0.5);
  \draw (5,-2.5) -- (3,-0.5);
  \draw (5,-2.5) -- (4,-0.5);
  \draw (6,-1.5) -- (5,-0.5);
  \draw (6,-1.5) -- (6,-0.5);
  \draw (6,-1.5) -- (7,-0.5);
  \draw (5,-2.5) -- (6,-1.5);
  \draw (6,-3.5) -- (5,-2.5);
  \draw (6,-3.5) -- (8,-0.5);
  \draw (6,-3.5) -- (9,-0.5);
  \draw (5,-4.5) -- (6,-3.5);
  \draw (6,-5.5) -- (5,-4.5);
  \draw (6,-5.5) -- (10,-0.5);
  \draw (6,-5.5) -- (11,-0.5);
  \draw (6,-6.5) -- (6,-5.5);
  \draw (6,-6.5) -- (12,-0.5);
  \fill (0,0.5) circle (2.5pt) (0,-0.5) circle (2.5pt);
  \fill (1,0.5) circle (2.5pt) (1,-0.5) circle (2.5pt);
  \fill (2,0.5) circle (2.5pt) (2,-0.5) circle (2.5pt);
  \fill (3,0.5) circle (2.5pt) (3,-0.5) circle (2.5pt);
  \fill (4,0.5) circle (2.5pt) (4,-0.5) circle (2.5pt);
  \fill (5,0.5) circle (2.5pt) (5,-0.5) circle (2.5pt);
  \fill (6,0.5) circle (2.5pt) (6,-0.5) circle (2.5pt);
  \fill (7,0.5) circle (2.5pt) (7,-0.5) circle (2.5pt);
  \fill (8,0.5) circle (2.5pt) (8,-0.5) circle (2.5pt);
  \fill (9,0.5) circle (2.5pt) (9,-0.5) circle (2.5pt);
  \fill (10,0.5) circle (2.5pt) (10,-0.5) circle (2.5pt);
  \fill (11,0.5) circle (2.5pt) (11,-0.5) circle (2.5pt);
  \fill (12,0.5) circle (2.5pt) (12,-0.5) circle (2.5pt);
\end{tikzpicture}} &
\resizebox{0.40\linewidth}{!}{\begin{tikzpicture}[scale=0.32, line cap=round, line join=round, >=stealth, thick]
  \draw[rounded corners=4pt] (5.682,5.682) -- (0.318,0.318) .. controls (0,0) .. (0.318,-0.318) -- (5.682,-5.682);
  \draw[rounded corners=4pt] (5.712,5.654) -- (1.288,0.346) .. controls (1,0) .. (1.318,-0.318) -- (4.682,-3.682);
  \draw[rounded corners=4pt] (6.318,5.682) -- (6.682,5.318);
  \draw[rounded corners=4pt] (6.682,4.682) -- (2.318,0.318) .. controls (2,0) .. (2.27,-0.36) -- (4.73,-3.64);
  \draw[rounded corners=4pt] (6.719,4.649) -- (3.281,0.351) .. controls (3,0) .. (3.318,-0.318) -- (4.682,-1.682);
  \draw[rounded corners=4pt] (7.318,4.682) -- (7.682,4.318);
  \draw[rounded corners=4pt] (7.682,3.682) -- (4.318,0.318) .. controls (4,0) .. (4.201,-0.402) -- (4.799,-1.598);
  \draw[rounded corners=4pt] (7.73,3.64) -- (5.27,0.36) .. controls (5,0) .. (5.318,-0.318) -- (5.682,-0.682);
  \draw[rounded corners=4pt] (8.318,3.682) -- (8.682,3.318);
  \draw[rounded corners=4pt] (8.682,2.682) -- (7.318,1.318);
  \draw[rounded corners=4pt] (9.201,2.598) -- (9.799,1.402);
  \draw[rounded corners=4pt] (9.318,2.682) -- (11.682,0.318) .. controls (12,0) .. (11.682,-0.318) -- (6.318,-5.682);
  \draw[rounded corners=4pt] (6.682,0.682) -- (6.318,0.318) .. controls (6,0) .. (6,-0.45) -- (6,-0.55);
  \draw[rounded corners=4pt] (7,0.55) -- (7,0.45) .. controls (7,0) .. (6.682,-0.318) -- (6.318,-0.682);
  \draw[rounded corners=4pt] (7.318,0.682) -- (7.682,0.318) .. controls (8,0) .. (7.75,-0.374) -- (6.25,-2.626);
  \draw[rounded corners=4pt] (9.682,0.682) -- (9.318,0.318) .. controls (9,0) .. (8.682,-0.318) -- (6.318,-2.682);
  \draw[rounded corners=4pt] (10,0.55) -- (10,0.45) .. controls (10,0) .. (9.719,-0.351) -- (6.281,-4.649);
  \draw[rounded corners=4pt] (10.318,0.682) -- (10.682,0.318) .. controls (11,0) .. (10.682,-0.318) -- (6.318,-4.682);
  \draw[rounded corners=4pt] (6,-5.55) -- (6,-5.45);
  \draw[rounded corners=4pt] (5.682,-4.682) -- (5.318,-4.318);
  \draw[rounded corners=4pt] (5.318,-3.682) -- (5.682,-3.318);
  \draw[rounded corners=4pt] (5.682,-2.682) -- (5.318,-2.318);
  \draw[rounded corners=4pt] (5.318,-1.682) -- (5.682,-1.318);
  \draw[rounded corners=4pt] (6,6.45) -- (6,7) -- (-1,7) -- (-1,-7) -- (6,-7) -- (6,-6.45);
  \draw (6,6.45) .. controls (6,6) .. (5.712,5.654);
  \draw (6.682,5.318) .. controls (7,5) .. (6.719,4.649);
  \draw (7.682,4.318) .. controls (8,4) .. (7.73,3.64);
  \draw (8.682,3.318) .. controls (9,3) .. (9.201,2.598);
  \draw (7.318,1.318) .. controls (7,1) .. (7,0.55);
  \draw (9.799,1.402) .. controls (10,1) .. (10,0.55);
  \draw (6,-6.45) .. controls (6,-6) .. (6,-5.55);
  \draw (6,-5.45) .. controls (6,-5) .. (6.281,-4.649);
  \draw (5.318,-4.318) .. controls (5,-4) .. (4.73,-3.64);
  \draw (5.682,-3.318) .. controls (6,-3) .. (6.25,-2.626);
  \draw (5.318,-2.318) .. controls (5,-2) .. (4.799,-1.598);
  \draw (5.682,-1.318) .. controls (6,-1) .. (6,-0.55);
  \draw[white, line width=6pt, line cap=butt] (5.682,5.682) .. controls (6,6) .. (6.318,5.682);
  \draw (5.682,5.682) .. controls (6,6) .. (6.318,5.682);
  \draw[white, line width=6pt, line cap=butt] (6.682,4.682) .. controls (7,5) .. (7.318,4.682);
  \draw (6.682,4.682) .. controls (7,5) .. (7.318,4.682);
  \draw[white, line width=6pt, line cap=butt] (7.682,3.682) .. controls (8,4) .. (8.318,3.682);
  \draw (7.682,3.682) .. controls (8,4) .. (8.318,3.682);
  \draw[white, line width=6pt, line cap=butt] (8.682,2.682) .. controls (9,3) .. (9.318,2.682);
  \draw (8.682,2.682) .. controls (9,3) .. (9.318,2.682);
  \draw[white, line width=6pt, line cap=butt] (6.682,0.682) .. controls (7,1) .. (7.318,0.682);
  \draw (6.682,0.682) .. controls (7,1) .. (7.318,0.682);
  \draw[white, line width=6pt, line cap=butt] (9.682,0.682) .. controls (10,1) .. (10.318,0.682);
  \draw (9.682,0.682) .. controls (10,1) .. (10.318,0.682);
  \draw[white, line width=6pt, line cap=butt] (5.682,-5.682) .. controls (6,-6) .. (6.318,-5.682);
  \draw (5.682,-5.682) .. controls (6,-6) .. (6.318,-5.682);
  \draw[white, line width=6pt, line cap=butt] (5.682,-4.682) .. controls (6,-5) .. (6.318,-4.682);
  \draw (5.682,-4.682) .. controls (6,-5) .. (6.318,-4.682);
  \draw[white, line width=6pt, line cap=butt] (4.682,-3.682) .. controls (5,-4) .. (5.318,-3.682);
  \draw (4.682,-3.682) .. controls (5,-4) .. (5.318,-3.682);
  \draw[white, line width=6pt, line cap=butt] (5.682,-2.682) .. controls (6,-3) .. (6.318,-2.682);
  \draw (5.682,-2.682) .. controls (6,-3) .. (6.318,-2.682);
  \draw[white, line width=6pt, line cap=butt] (4.682,-1.682) .. controls (5,-2) .. (5.318,-1.682);
  \draw (4.682,-1.682) .. controls (5,-2) .. (5.318,-1.682);
  \draw[white, line width=6pt, line cap=butt] (5.682,-0.682) .. controls (6,-1) .. (6.318,-0.682);
  \draw (5.682,-0.682) .. controls (6,-1) .. (6.318,-0.682);
\end{tikzpicture}}
\end{tabular}
\caption{$8_{17}$ as an unoriented knot. Top: a minimal tree diagram in $F$ (9 leaves) and its knot.
Bottom: a minimal ternary tree diagram ($n=7$) and its knot.}
\label{fig:817-unoriented}
\end{figure}
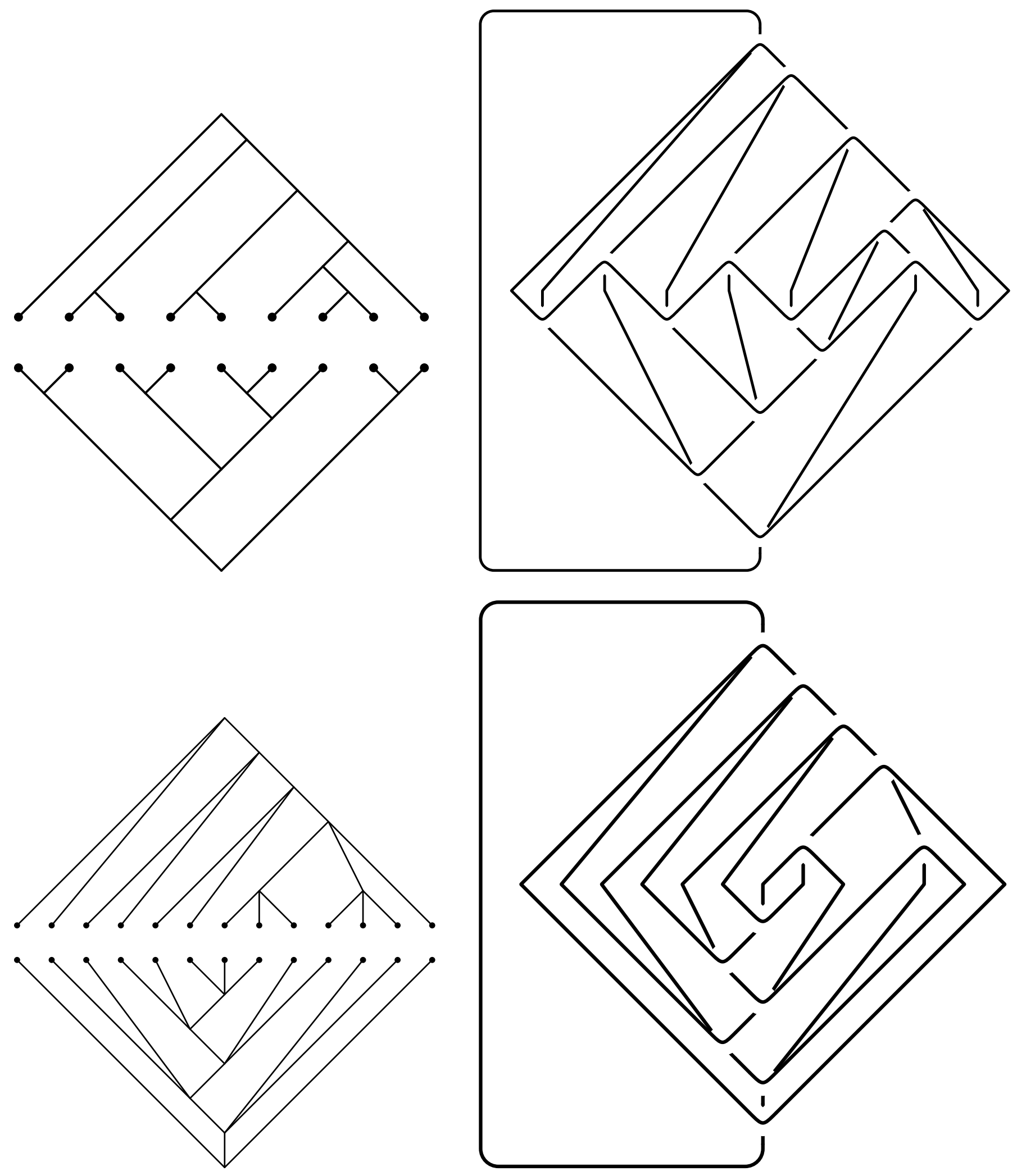

\section{Some old and some new results on the indices}\label{sec:known}

We begin this section with a  collection of results on the four indices: Theorem~\ref{thm:unoriented} for unoriented knots and
links, Theorem~\ref{thm:oriented} for oriented ones. Most of them are known; in some cases there are no explicit statements in the literature for the groups $F_3$ and $\vFt$, but only for $F$ and $\vF$, however the proofs adapt easily, as explained below.


\begin{theorem}\label{thm:unoriented}
Let $L$ be an unoriented link and $K$ a knot with crossing number $c$.
\begin{enumerate}[(i),ref=(\roman*)]
\item\label{it:alexander} $L=L(g)$ for some $g\in F$, and for some $g\in F_3$
  \cite{Jones2017,Jones2019}.
\item\label{it:mirror} For $g$ in $F$ or in $F_3$, the element $g^{-1}$ is represented by a tree diagram of the same
  size as $g$, and $L(g^{-1})=L(g)^m$. Hence $\ind_F(L^m)=\ind_F(L)$ and $\ind_{F_3}(L^m)=\ind_{F_3}(L)$.
\item\label{it:golan-sapir} If $L$ contains $u$ unlinked unknots and has a diagram with $n$ crossings, then
  $\ind_F(L)\le 12n+u+3$~\cite[Theorem~6.13]{GolanSapir2017}. If $L$ contains no unlinked unknots, then $L$ is
  represented by a reduced tree diagram with at most $12n$ leaves~\cite[Lemma~6.12]{GolanSapir2017}; in particular
  $\ind_F(K)\le 12c$ for every non-trivial knot $K$.
\item\label{it:grid} $\operatorname{bridge}(L)\le\operatorname{braid}(L)\le\operatorname{grid}(L)/2\le
  \ind_F(L)$~\cite[Theorem~1.15]{LuoWan2024}.
\item\label{it:crossings} $c\le 2(\ind_F(K)-1)$ and $c\le 2(\ind_{F_3}(K)-1)$.
\item\label{it:iota} $\ind_{F_3}(K)\le\ind_F(K)$.
\item\label{it:sum} For knots $K_1$ and $K_2$, $\ind_F(K_1\#K_2)\le\ind_F(K_1)+\ind_F(K_2)-1$ and
  $\ind_{F_3}(K_1\#K_2)\le\ind_{F_3}(K_1)+\ind_{F_3}(K_2)-1$, for both connected sums of $K_1$ and $K_2$.
\end{enumerate}
\end{theorem}
\begin{proof}
Statements~\ref{it:alexander}, \ref{it:golan-sapir} and~\ref{it:grid} are proved in the references (the
diagrams of~\cite{GolanSapir2017} over $\langle x\mid x=x^2\rangle$ have one vertex more than leaves). For~\ref{it:mirror},
the diagram of $g^{-1}=(T_-,T_+)$ is the reflection of that of $g$ in the line of the leaves, with the same over-
and under-strands, that is, the image of $L(g)$ under a reflection of $\mathbb R^3$; since $g\mapsto g^{-1}$ is a
bijection between the tree diagrams of size $n$ of $L$ and those of $L^m$, the two links have the same index.
Statement~\ref{it:crossings} holds because a tree diagram of size $n$ gives a link diagram with $2(n-1)$
crossings, and~\ref{it:iota} because Jones' map $\iota\colon F\to F_3$ preserves the size and the link.

For~\ref{it:sum}, let $(T_+,T_-)$ and $(S_+,S_-)$ be minimal tree diagrams of $K_1$ and $K_2$, of sizes $n_1$ and
$n_2$. Inserting $S_+$ and $S_-$ in place of the first leaf of $T_+$ and of $T_-$ (a \emph{graft}) gives a tree
diagram of size $n_1+n_2-1$, both in $F$ and in $F_3$. Its link is a connected sum of $K_1$ and $K_2$, since
the part of the diagram coming from $(S_+,S_-)$ is a $1$-tangle whose closure is $K_2$. Grafting $(S_-,S_+)$
instead, which by~\ref{it:mirror} has the same size and gives $K_2^m$, yields the other connected sum.
\end{proof}

\begin{theorem}\label{thm:oriented}
Let $K$ be an oriented knot and $L$ an oriented link.
\begin{enumerate}[(i),ref=(\roman*)]
\item\label{it:alexander-or} $L=L(g)$ for some $g\in\vF$~\cite{Aiello2020}, and for some $g\in\vFt$.
\item\label{it:parity} For $g$ in $\vF$ or in $\vFt$ of size $n$, the number of components of $L(g)$ is
  congruent to $n$ modulo $2$~\cite[Proposition~3.11]{AielloConti2019b}; in particular knots only arise for
  odd $n$.
\item\label{it:inverse} For $g$ in $\vF$ or in $\vFt$, the element $g^{-1}$ is represented by a tree diagram of the
  same size as $g$, and $L(g^{-1})=L(g)^{mr}$~\cite[Remark~3.8]{AielloConti2019b}. Hence
  $\ind_{\vF}(L^{mr})=\ind_{\vF}(L)$ and $\ind_{\vF}(L^{r})=\ind_{\vF}(L^{m})$, and the same holds for $\ind_{\vFt}$.
\item\label{it:inequalities} $\ind_F(K)\le\ind_{\vF}(K)$ and $\ind_{F_3}(K)\le\ind_{\vFt}(K)\le\ind_{\vF}(K)$.
\item\label{it:genus} $\chi(L)\ge 2-\ind_{\vF}(L)$, where $\chi(L)$ is the maximal Euler characteristic of a
  Seifert surface of $L$; hence $g(K)\le(\ind_{\vF}(K)-1)/2$~\cite[Theorem~1.13]{LuoWan2024};
  moreover $g(K)\le(\ind_{\vFt}(K)-1)/2$.
\item\label{it:reverse} $|\ind_{\vF}(K)-\ind_{\vF}(K^r)|\le 4$~\cite[Theorem~4.2]{Aiello2026}. More precisely,
  if $g$ is a minimal diagram of $K$ in $\vF$, then $\ind_{\vF}(K^r)\le\ind_{\vF}(K)+2$ when the first and last
  vertices of $\Gamma(g)$ have different colours, and $\ind_{\vF}(K^r)\le\ind_{\vF}(K)+4$ when they have the
  same colour~\cite[proof of Theorem~4.2]{Aiello2026}.
\item\label{it:reverse3} $|\ind_{\vFt}(K)-\ind_{\vFt}(K^r)|\le 4$.\footnote{The proof of the corresponding
  statement for $\vF$ in the published version of~\cite{Aiello2026} contains a mistake in one of the two cases.
  We therefore give a complete proof below, which applies to both $\vF$ and $\vFt$.} More precisely, if $g$ is a minimal diagram of
  $K$ in $\vFt$, then $\ind_{\vFt}(K^r)\le\ind_{\vFt}(K)+2$ when the first and last vertices of $\Gamma(g)$ have
  different colours, and $\ind_{\vFt}(K^r)\le\ind_{\vFt}(K)+4$ when they have the same colour.
\item\label{it:sum-or} For oriented knots $K_1$ and $K_2$, it holds $\ind_{\vF}(K_1\#K_2)\le\ind_{\vF}(K_1)+\ind_{\vF}(K_2)-1$
  and $\ind_{\vFt}(K_1\#K_2)\le\ind_{\vFt}(K_1)+\ind_{\vFt}(K_2)-1$.
\end{enumerate}
\end{theorem}
\begin{proof}
Statements~\ref{it:alexander-or}--\ref{it:reverse} are proved in the references for $\vF$, and the proofs adapt
easily to $\vFt$. For~\ref{it:alexander-or}, $\iota$ maps $\vF$ into $\vFt$, preserving $n$ and the oriented link,
since the Tait graphs of $g$ and $\iota(g)$ coincide. For~\ref{it:parity}, the Tait graph of $g\in\vFt$ has $n$
vertices and $2(n-1)$ edges and is bipartite, as for $\vF$ (a second proof, via the Laplacian of the dual graph
modulo $2$, is in the repository). For~\ref{it:inverse}, the diagram of $g^{-1}=(T_-,T_+)$, which has the same size
as $g$, is the reflection of that of $g$ in the line of the leaves, with the same over- and under-strands, i.e.\
the image of $L(g)$ under an orientation-reversing reflection of $\mathbb R^3$; the reflection fixes the regions and
their colours but reverses the sense of rotation of the plane, so it reverses the orientation. Hence
$L(g^{-1})=L(g)^{mr}$, and the equalities between the indices follow as in
Theorem~\ref{thm:unoriented}\ref{it:mirror}, using $(L^{mr})^{mr}=L$ and $(L^{r})^{mr}=L^{m}$.
Statement~\ref{it:inequalities} follows from $\vF\subset F$, $\vFt\subset F_3$ and from $\iota$ as
in~\ref{it:alexander-or}. For~\ref{it:genus}, the checkerboard surface of $g\in\vFt$ is a Seifert surface made of
$n$ disks and $2(n-1)$ bands. For~\ref{it:reverse} and~\ref{it:reverse3}, the Tait graph of $g$ is modified near its first and last vertices,
adding two or four vertices according to whether these have different or equal colours, by moves that can be undone
by Reidemeister moves of type I and II but exchange the colours of the old vertices, hence reverse the orientation
(see~\cite[Theorem~4.2]{Aiello2026} for $\vF$). We give the details for $\vFt$; they apply verbatim to $\vF$, as
explained at the end of the proof.

For~\ref{it:reverse3}, let $(T_+,T_-)$ be a minimal diagram of $K$ in $\vFt$; its first Tait vertex has colour
$+$. If the last one has colour $-$, replace $T_+$ and $T_-$ by the trees $T'_+$ and $T'_-$ below (the triangles
are $T_+$ and $T_-$, both drawn with the root at the top):
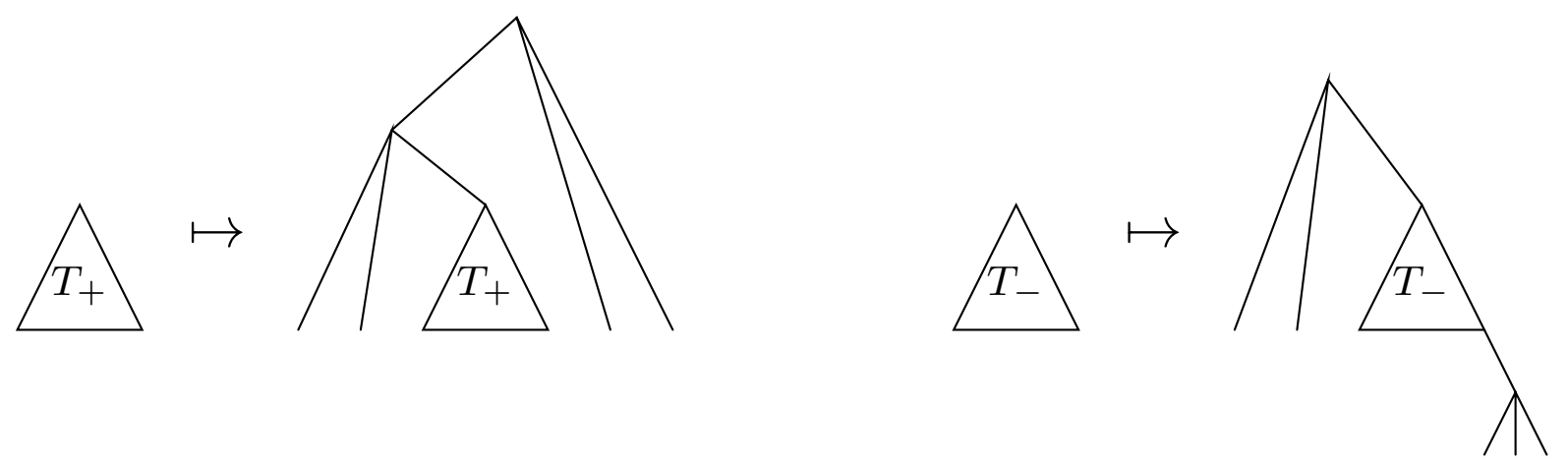
\begin{center}
\begin{tikzpicture}[scale=0.42,line cap=round]
  \draw (1,2)--(0,0)--(2,0)--cycle; \node at (1,0.7) {\scriptsize $T_+$};
  \node at (3.2,1.5) {$\mapsto$};
  \begin{scope}[shift={(4.5,0)}]
    \draw (3,2)--(2,0)--(4,0)--cycle; \node at (3,0.7) {\scriptsize $T_+$};
    \draw (0,0)--(1.5,3.2)--(1,0) (1.5,3.2)--(3,2);
    \draw (1.5,3.2)--(3.5,5)--(5,0) (3.5,5)--(6,0);
  \end{scope}
  \begin{scope}[shift={(15,0)}]
    \draw (1,2)--(0,0)--(2,0)--cycle; \node at (1,0.7) {\scriptsize $T_-$};
    \node at (3.2,1.5) {$\mapsto$};
  \end{scope}
  \begin{scope}[shift={(19.5,0)}]
    \draw (3,2)--(2,0)--(4,0)--cycle; \node at (3,0.7) {\scriptsize $T_-$};
    \draw (0,0)--(1.5,4)--(1,0) (1.5,4)--(3,2);
    \draw (4,0)--(4.5,-1) (4,-2)--(4.5,-1)--(5,-2) (4.5,-1)--(4.5,-2);
  \end{scope}
\end{tikzpicture}
\end{center}
that is, $T'_+$ is obtained by hanging $T_+$ from the right child of a new caret, itself the left child of a new
root, and $T'_-$ by hanging $T_-$ from the right child of a new root and adding a caret on its last leaf. The
diagram $(T'_+,T'_-)$ has size $n+2$. Drawing it, one sees that its Tait graph has two new vertices, the first one
of colour $+$ and the last one of colour $-$, and that the colours of all the old vertices are exchanged; hence it
belongs to $\vFt$, its checkerboard surface is orientable, and its link is isotopic to $L(T_+,T_-)$ with the opposite
orientation, that is, $L(T'_+,T'_-)=K^r$. If the last vertex has colour $+$, replace instead $T_+$ and $T_-$ by
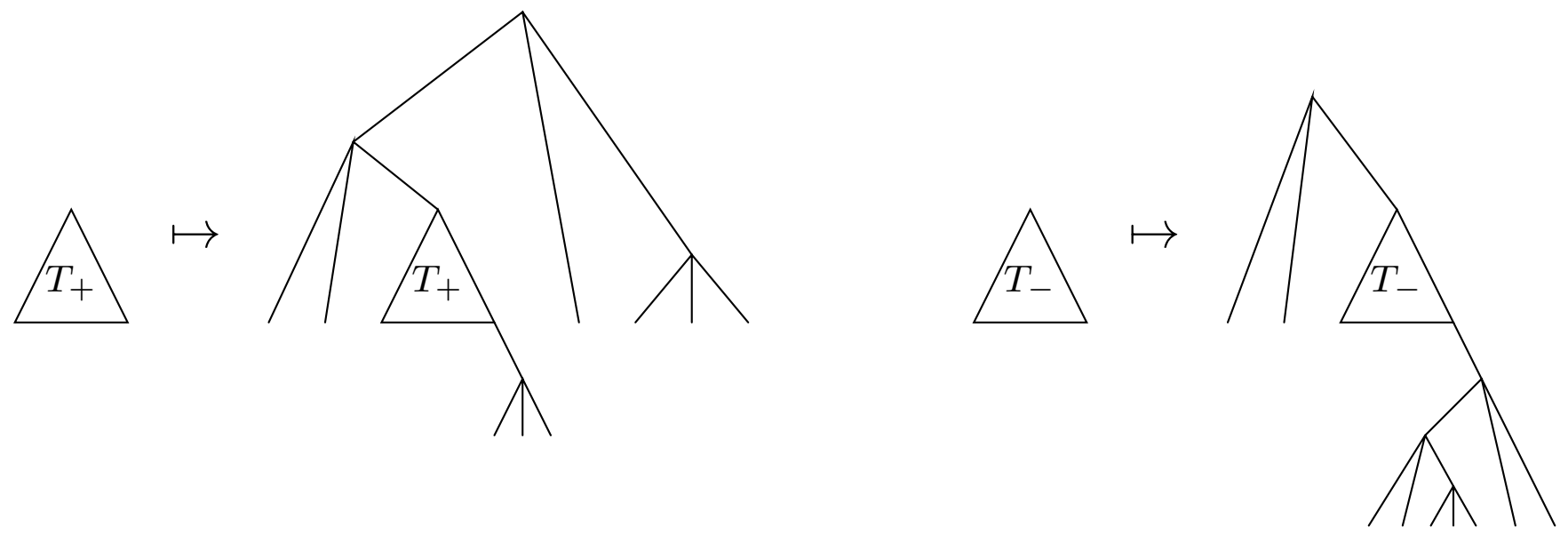
\begin{center}
\begin{tikzpicture}[scale=0.42,line cap=round]
  \draw (1,2)--(0,0)--(2,0)--cycle; \node at (1,0.7) {\scriptsize $T_+$};
  \node at (3.2,1.5) {$\mapsto$};
  \begin{scope}[shift={(4.5,0)}]
    \draw (3,2)--(2,0)--(4,0)--cycle; \node at (3,0.7) {\scriptsize $T_+$};
    \draw (4,0)--(4.5,-1) (4,-2)--(4.5,-1)--(5,-2) (4.5,-1)--(4.5,-2);
    \draw (0,0)--(1.5,3.2)--(1,0) (1.5,3.2)--(3,2);
    \draw (1.5,3.2)--(4.5,5.5)--(5.5,0) (4.5,5.5)--(7.5,1.2);
    \draw (6.5,0)--(7.5,1.2)--(8.5,0) (7.5,1.2)--(7.5,0);
  \end{scope}
  \begin{scope}[shift={(17,0)}]
    \draw (1,2)--(0,0)--(2,0)--cycle; \node at (1,0.7) {\scriptsize $T_-$};
    \node at (3.2,1.5) {$\mapsto$};
  \end{scope}
  \begin{scope}[shift={(21.5,0)}]
    \draw (3,2)--(2,0)--(4,0)--cycle; \node at (3,0.7) {\scriptsize $T_-$};
    \draw (0,0)--(1.5,4)--(1,0) (1.5,4)--(3,2);
    \draw (4,0)--(4.5,-1);
    \draw (4.5,-1)--(3.5,-2) (4.5,-1)--(5.1,-3.6) (4.5,-1)--(5.8,-3.6);
    \draw (3.5,-2)--(2.5,-3.6) (3.5,-2)--(3.1,-3.6) (3.5,-2)--(4,-2.9);
    \draw (4,-2.9)--(3.6,-3.6) (4,-2.9)--(4,-3.6) (4,-2.9)--(4.4,-3.6);
  \end{scope}
\end{tikzpicture}
\end{center}
that is, $T'_+$ is obtained by adding a caret on the last leaf of $T_+$, hanging the result from the right child
of a new caret, itself the left child of a new root whose right child is a further new caret, and $T'_-$ by
hanging $T_-$ from the right child of a new root and attaching to its last leaf the tree with three internal
vertices shown above. This gives a diagram of size $n+4$; its Tait graph has one new vertex on the left and three
on the right, and drawing it one sees in the same way that it belongs to $\vFt$ and that its link is $K^r$. For~\ref{it:reverse}, if $T_\pm$ are binary trees then $T'_\pm$ are the images under $\iota$ of binary trees
(replace each new ternary caret by a binary one), and the argument above only uses the Tait graph, which is the
same for a binary tree diagram and its image under $\iota$; hence the same two moves, on binary trees, prove the
statement for $\vF$. Similar moves were employed in~\cite{Terron2025}, but they do not preserve $\vFt$.

For~\ref{it:sum-or}, consider the graft at the first leaf of minimal diagrams of $K_1$ and $K_2$, as in the proof
of Theorem~\ref{thm:unoriented}\ref{it:sum}. The path from the root to the first leaf of $T_\pm$ consists of left
edges only; hence for every new leaf the number $r_\pm$ of right edges (in $\vF$), and the weight $c_\pm$ (in
$\vFt$: these left edges lie in the first piece of the path), are those of the corresponding leaf of $S_\pm$, while
the leaves of $T_\pm$ after the first one keep their values, and in $\vFt$ also their parity, since they move by
$2(n_2-1)$ positions. Hence the graft belongs to $\vF$ (resp.\ $\vFt$), and the Tait vertices of the two diagrams keep
their colours, so that both pieces keep their orientations: the link of the graft is the oriented connected sum
$K_1\#K_2$.
\end{proof}

By~\ref{it:parity}, \ref{it:reverse} of Theorem \ref{thm:oriented} we get $\ind_{\vF}(K^r)-\ind_{\vF}(K)\in\{-4,-2,0,2,4\}$; all the differences we found are equal to $\pm2$. By~\ref{it:inverse}
an oriented knot has at most two oriented index values, $\ind(K)=\ind(K^{mr})$ and $\ind(K^r)=\ind(K^m)$; for
an invertible knot they coincide.

In~\cite{Aiello2026} two questions are asked: whether some knot $K$ has $\ind_{\vF}(K)\ne\ind_{\vF}(K^r)$
(Question~1), and whether every knot has minimal diagrams of both types, with first and last vertices of the
same or of different colours (Question~2). Both are answered in Section~\ref{sec:oriented}.

\medskip\noindent\textbf{Non-prime diagrams.}\label{sec:nonprime} The enumeration of
Section~\ref{sec:computation} skips the non-prime diagrams; the following lemma and corollary show that this loses
nothing, also in the oriented groups. A tree diagram $(T_+,T_-)$ is \emph{non-prime} if a subtree of $T_+$ and a subtree of $T_-$, both rooted at
internal vertices other than the roots, have the same leaves, that is, the same set of consecutive positions among
the leaves $1,2,\dots$. It is then a graft: it is obtained by inserting the smaller diagram formed by the two
subtrees (the \emph{inner piece}) in place of a leaf of another diagram (the \emph{outer piece}), and its link is
the connected sum of the links of the two pieces (proof of Theorem~\ref{thm:unoriented}\ref{it:sum}). Hence a prime
knot realized by a non-prime diagram of $F$ or $F_3$ is already realized by the smaller piece. In the oriented groups
the pieces may lie outside the group or carry the opposite orientation; nevertheless the following holds.

\begin{lemma}\label{lem:nonprime}
Let $g$ be a non-prime diagram of size $k$ of $\vF$ (resp.\ $\vFt$) whose oriented link is a prime knot $J$. Then
one of $J$, $J^r$, $J^m$, $J^{mr}$ is realized in $\vF$ (resp.\ $\vFt$) by a diagram of size at most $k-4$; in
particular $\min\{\ind_{\vF}(J),\ind_{\vF}(J^r)\}\le k-4$ (resp.\ $\min\{\ind_{\vFt}(J),\ind_{\vFt}(J^r)\}\le k-4$).
\end{lemma}

\begin{proof}
We give the proof for $\vFt$; for $\vF$ it is the same, with the number $r_\pm$ of right edges in place of the
weight $c_\pm$ and without middle edges. Since $g$ is non-prime, there are internal vertices $u_+$ of $T_+$ and $u_-$ of $T_-$,
different from the roots, such that the leaves below $u_+$ in $T_+$ are exactly the leaves below $u_-$ in $T_-$; these
leaves are consecutive, and we denote by $I$ the set of their positions, an interval of integers in the numbering
$1,2,\dots$ of the leaves (not a subinterval of $[0,1]$). Let $\pi_\pm$ be the path from the root of $T_\pm$ to $u_\pm$, and let $A$ and $B$ be the outer and
the inner piece. For a leaf $\ell\in I$ the path from the root of $T_\pm$ to $\ell$ is $\pi_\pm$ followed by the path
$\sigma_\pm(\ell)$ from $u_\pm$ to $\ell$ in $B$, and
\[
c_\pm(\ell)=c(\pi_\pm)+c^{(\mu)}(\sigma_\pm(\ell)),
\]
where $\mu$ is the number of middle edges of $\pi_\pm$ and $c^{(\mu)}$ is the weight computed with the pieces of the
path numbered from $\mu+1$. The parity of $\mu$ is the same for $\pi_+$ and $\pi_-$: in a ternary tree the first leaf
below a vertex is odd exactly when the path to the vertex contains an even number of middle edges, and $u_+$ and
$u_-$ have the same first leaf. If $\mu$ is even, $c^{(\mu)}=c$; if $\mu$ is odd, $c^{(\mu)}(\sigma)=c(\bar\sigma)$,
where $\bar\sigma$ is the path obtained by exchanging left and right edges, since this exchanges the roles of the
odd and the even pieces.

Apply $c_+\equiv c_-$ to the first leaf $\ell$ of $I$ if $\mu$ is even, and to its last leaf $\ell$ if $\mu$ is
odd: the path $\sigma_\pm(\ell)$ from $u_\pm$ to this leaf consists of left (resp.\ right) edges only. Since it
contains no middle edge, it lies in the last piece of the path from the root to $\ell$, which has number $\mu+1$,
odd (resp.\ even); hence $c^{(\mu)}(\sigma_\pm(\ell))=0$, and we get $c(\pi_+)\equiv c(\pi_-)$. Hence $c^{(\mu)}(\sigma_+(\ell))\equiv
c^{(\mu)}(\sigma_-(\ell))$ for every $\ell\in I$. Since the leaves of $A$ other than the one replacing $I$ keep their
paths, and the one replacing $I$ has weights $c(\pi_\pm)$, we get $A\in\vFt$. If $\mu$ is even, that is, if $I$
starts at an odd leaf, the weights simply add up, and $g\in\vFt$ if and only if $A\in\vFt$ and $B\in\vFt$. If $\mu$ is
odd, that is, if $I$ starts at an even leaf, $B$ need not belong to $\vFt$, but its vertical reflection $\bar
B=(\sigma B_+,\sigma B_-)$ does, since $c^{(\mu)}(\sigma)=c(\bar\sigma)$.

As unoriented knots, $L(g)=L(A)\#L(B)$, as in the proof of Theorem~\ref{thm:unoriented}\ref{it:sum} (the graft
may be at any leaf), and $L(\bar B)=L(B)^m$. Since $J$ is prime, one of $L(A)$ and $L(B)$ is the unknot. The smallest
diagrams of $\vF$ and of $\vFt$ realizing a knot have size $5$ (with $n=3$ only links with $3$ components occur), and
the sizes of $A$ and $B$ add up to $k+1$; hence, applying this to the piece giving the unknot (or to $\bar B$), the
piece giving $J$ has size at most $k-4$. If it is $A$, or $B$ with
$\mu$ even, it belongs to the group and its oriented link is $J$ or $J^r$; if it is $B$ with $\mu$ odd, then $\bar B$
belongs to the group, has the same size, and its oriented link is $J^m$ or $J^{mr}$. The last claim follows from
Theorem~\ref{thm:oriented}\ref{it:inverse}.
\end{proof}

\begin{corollary}\label{cor:nonprime}
Let $J$ be a non-invertible prime knot, let $a=\min\{\ind_{\vF}(J),\ind_{\vF}(J^r)\}$ (resp.\ $a=\min\{\ind_{\vFt}(J),\ind_{\vFt}(J^r)\}$), and let $b'$ be the
smallest size of a prime diagram of $\vF$ (resp.\ $\vFt$) realizing the other orientation. Then $a$ is attained by
a prime diagram, and the $\vF$-index (resp.\ $\vFt$-index) of the other orientation is $\min\{b',a+4\}$.
\end{corollary}

\begin{proof}
A non-prime diagram of size $k$ realizing $J$ or $J^r$ has $k\ge a+4$ by Lemma~\ref{lem:nonprime}, so the
minimum $a$ is attained by a prime diagram, and the other orientation is realized by a non-prime diagram only at
sizes $\ge a+4$. By Theorem~\ref{thm:oriented}\ref{it:reverse} and~\ref{it:reverse3} its index is at most $a+4$.
\end{proof}

Hence enumerating only the prime diagrams loses nothing: it gives the index up to orientation exactly, and, with
Theorem~\ref{thm:oriented}\ref{it:reverse}--\ref{it:reverse3}, the index of both orientations.

We end this section with two new lower bounds for the indices in terms of polynomial invariants: the span of
the Jones polynomial bounds the $F$- and the $\vF$-index (Theorem~\ref{thm:span}), and the degrees of the
HOMFLY-PT polynomial bound the oriented indices (Proposition~\ref{prop:homfly}).

Recall from Section~\ref{sec:construction} that the Tait diagram $\Gamma(T_+,T_-)$ of a binary tree diagram of
size $n$ is the union of two spanning trees $\Gamma_+(T_+)$ and $\Gamma_-(T_-)$, each with $n-1$ edges, the first
one with only positive edges and the second one with only negative edges, in the sense of
Figure~\ref{fig:smoothings}.

\begin{theorem}\label{thm:span}
For every link $L$,
\[
\ind_{\vF}(L)\geq \ind_F(L)\geq\operatorname{span}V_L(t)+1 ,
\]
where $V_L$ is the Jones polynomial and $\operatorname{span}$ denotes the difference between the maximal and the
minimal degree.
\end{theorem}

\begin{proof}
Let $(T_+,T_-)$ be a tree diagram with $n$ leaves such that $L(T_+,T_-)=L$, let $D$ be the corresponding
link diagram, and let $\Gamma=\Gamma(T_+,T_-)$ be its Tait graph. As recalled above, $\Gamma$ is the union of the
two spanning trees $\Gamma_+(T_+)$ and $\Gamma_-(T_-)$, each with $n-1$ edges, whose edges correspond to positive and
to negative crossings, respectively. Hence $D$ has $k=2(n-1)$ crossings.

Shade the regions of $D$ in a checkerboard fashion, so that the vertices of $\Gamma$ are the black regions. Since
$\Gamma$ is connected, the projection of $D$ is connected, so every bounded region is a disc; in particular the black
regions, which are bounded, are discs. At each crossing two black corners and two white corners meet, and the
corresponding edge of $\Gamma$ joins the two black regions. One of the two smoothings of the crossing opens a channel
between the two black corners, and the other one opens a channel between the two white corners: at a positive
crossing the $B$-smoothing joins the black corners, while at a negative crossing the $A$-smoothing does
(Figure~\ref{fig:smoothings}).

Fix a state $S$, and let $\Gamma_S\subseteq\Gamma$ be the subgraph with the same vertices as $\Gamma$ and with only
those edges at whose crossing $S$ chooses the smoothing that joins the two black corners. After smoothing all the
crossings according to $S$, keep the shading: the black part of the plane is then the union of the black regions,
which are discs, together with a thin black band at each crossing corresponding to an edge of $\Gamma_S$. We denote
this black part by $\Sigma_S$; it is a thickening of the graph $\Gamma_S$, each vertex becoming a disc and each edge a
band. The smoothed diagram has no crossings, and each of its curves has black on one side and white on the other.
Hence the circles of the state $S$ are exactly the boundary curves of $\Sigma_S$.

In the state in which all the smoothings are of type $A$, the graph $\Gamma_S$ consists of all the vertices and of
the negative edges, that is, $\Gamma_S=\Gamma_-(T_-)$; in the state in which all the smoothings are of type $B$, it
consists of all the vertices and of the positive edges, that is, $\Gamma_S=\Gamma_+(T_+)$. If $\Gamma_S$ is a tree,
then $\Sigma_S$ is a disc: indeed, adding a leaf to a tree corresponds to attaching a new disc to the previous
surface by a band glued along an arc of each boundary, which gives again a disc. Alternatively, $\Sigma_S$ is a
connected planar surface with Euler characteristic $n-(n-1)=1$. In both cases the boundary of $\Sigma_S$ is a single
circle, so the numbers of circles of the two extreme states are $s_A=s_B=1$.

With the normalisation $\langle\bigcirc\rangle=1$, the Kauffman bracket is
$\langle D\rangle=\sum_S A^{a(S)-b(S)}\delta^{|S|-1}$, with $\delta=-A^2-A^{-2}$, where $a(S)$ and $b(S)$ are the
numbers of $A$- and $B$-smoothings of $S$ and $|S|$ is its number of circles. The maximal degree of the term of a
state $S$ is $a(S)-b(S)+2(|S|-1)$, and changing one smoothing from $A$ to $B$ decreases $a-b$ by $2$ and changes
$|S|$ by $\pm1$, so this quantity does not increase. This gives the standard estimates
\[
\max\deg\langle D\rangle\leq k+2(s_A-1),\qquad \min\deg\langle D\rangle\geq -k-2(s_B-1),
\]
and therefore
\[
\operatorname{span}\langle D\rangle\leq 2k+2(s_A+s_B)-4=2k=4(n-1).
\]

Note that $\langle D\rangle$ does not depend on any orientation of $L$; changing the orientation of some components
of $L$ only changes the writhe, hence multiplies $V_L$ by a monomial. In particular $\operatorname{span}V_L$ is
well defined for unoriented links. For any orientation, $V_L(t)=(-A^3)^{-\operatorname{wr}(D)}\langle D\rangle$, where $\operatorname{wr}(D)$ is the writhe of $D$, with $A=t^{-1/4}$,
so that $\operatorname{span}_A\langle D\rangle=4\operatorname{span}_tV_L$. Hence $\operatorname{span}V_L\leq n-1$,
and taking the minimum over all tree diagrams of $L$ gives $\ind_F(L)\geq\operatorname{span}V_L+1$.

Finally, if $L$ is oriented and $g\in\vF$ is represented by a tree diagram with $n$ leaves giving $L$, the same tree
diagram, viewed as a diagram of an element of $F$, gives the underlying unoriented link; hence $\ind_{\vF}(L)\geq
\ind_F(L)$, where on the right $L$ denotes the underlying unoriented link.
\end{proof}

\begin{corollary}\label{cor:alternating}\label{cor:parity}
Let $L$ be a non-split alternating link, for instance an alternating knot, with crossing number $c(L)$. Then
\[
\ind_{\vF}(L)\geq \ind_F(L)\geq c(L)+1 .
\]
If, moreover, $L$ is a knot and $c(L)$ is odd, then $\ind_{\vF}(L)\geq c(L)+2$.
\end{corollary}

\begin{proof}
For a non-split alternating link $\operatorname{span}V_L=c(L)$, by the theorem of Kauffman, Murasugi and
Thistlethwaite~\cite{Kauffman1987,Murasugi1987,Thistlethwaite1987}, so the first claim follows from Theorem~\ref{thm:span}. If $g\in\vF$ is represented by a tree
diagram with $n$ leaves, then the number of components of the associated link is at most $n$ and congruent to $n$
modulo $2$~\cite[Proposition~3.11]{AielloConti2019b}. Hence, if
$L$ is a knot, $\ind_{\vF}(L)$ is odd, and when $c(L)$ is odd the even number $c(L)+1$ cannot be attained.
\end{proof}

\begin{remark}
Without using the Jones polynomial one only gets $k=2(n-1)\geq c(L)$, that is $\ind_F(L)\geq c(L)/2+1$; in the
alternating case, Theorem~\ref{thm:span} doubles this estimate. A computer enumeration of all tree diagrams with at
most $7$ leaves confirms the inequality of Theorem~\ref{thm:span} for all the $39$ link types that occur, but it
also shows that the bound $\ind_F(K)\geq c(K)+1$ is not attained by the first knots:
\begin{center}
\begin{tabular}{@{}lccc@{}}
\toprule
knot & $c(K)=\operatorname{span}V_K$ & $c(K)+1$ & $\ind_F(K)$\\
\midrule
trefoil & $3$ & $4$ & $5$\\
figure-eight knot & $4$ & $5$ & $7$\\
knots with $5$ crossings ($5_1$ or $5_2$) & $5$ & $6$ & $7$\\
\bottomrule
\end{tabular}
\end{center}
Since all tree diagrams with at most $7$ leaves have been examined, the values in the last column are the actual
indices. The bound is nevertheless attained by larger knots. The first one is $7_4$, with $c(7_4)=7$ and
$\ind_F(7_4)=8$; another one is $8_{17}$, with $\ind_F(8_{17})=9$. Over the 6\,051 knots with at most 13 crossings
whose $F$-index we know, the equality $\ind_F(K)=\operatorname{span}V_K+1$ holds for 1\,440 knots (1 with 7
crossings, 7 with 8, 12 with 9, 43 with 10, 135 with 11, 583 with 12 and 659 with 13), about 650 of them
alternating. For $\vF$ the bound of Theorem~\ref{thm:span}, raised to the next odd integer since the $\vF$-index of a
knot is odd (Theorem~\ref{thm:oriented}\ref{it:parity}), is attained by 8 knots: the alternating knots $11a_5$,
$11a_{128}$, $13a_{31}$, $13a_{212}$, $13a_{757}$, $13a_{1713}$, $13a_{2227}$, for which it is the bound
$\ind_{\vF}(K)\ge c(K)+2$ of Corollary~\ref{cor:parity}, and the non-alternating knot $13n_{2273}$, with
$\operatorname{span}V=11$ and $\ind_{\vF}=13=c$.
For all these knots the index is proved by a single tree diagram, without any
enumeration. The estimate in the proof of Theorem~\ref{thm:span} is sharp exactly when the extreme terms of
$\langle D\rangle$ do not cancel. It would be interesting to characterize the knots for which this happens.

The bound does not hold for $F_3$ and $\vFt$: 105 knots have a minimal crossing diagram coming from $F_3$, that
is $c(K)=2(\ind_{F_3}(K)-1)$ (Section~\ref{sec:bounds}).
\end{remark}

\begin{proposition}\label{prop:homfly}
Let $L$ be an oriented link, and write its HOMFLY-PT polynomial, normalised as in~\cite{Morton1986} by
$v^{-1}P_{L_+}-vP_{L_-}=zP_{L_0}$, as $P_L(v,z)=\sum_{k=e}^{E}a_k(z)v^k$ with $a_e\neq0\neq a_E$; let $M$ be its
maximal degree in $z$. Then
\[
\ind_{\vF}(L)\geq \max\{\operatorname{braid}(L),\,1+E,\,1-e,\,1+M\},
\qquad
\ind_{\vFt}(L)\geq \max\{\operatorname{braid}(L),\,1+M\},
\]
and in particular $\ind_{\vFt}(L)\geq 1+\max\{(E-e)/2,\,M\}$.
\end{proposition}

\begin{proof}
Let $g$ be an element of $\vF$ or of $\vFt$ represented by a tree diagram of size $n$ giving $L$, and let $D$ be the
corresponding oriented link diagram, which has $c(D)=2(n-1)$ crossings. Moreover, $D$ has $s(D)=n$ Seifert circles.
Indeed, the boundary of each black region is oriented coherently, counterclockwise or clockwise according to the
colour of the corresponding vertex of the Tait graph, and adjacent black regions have opposite colours; hence at
every crossing the smoothing respecting the orientation is the one separating the two black regions, and the
Seifert circles of $D$ are the boundaries of the $n$ black regions. By~\cite[Theorem~2]{Morton1986},
$M\leq c(D)-(s(D)-1)=n-1$. By Yamada's theorem~\cite{Yamada1987}, the braid index of $L$ is the minimal number of
Seifert circles of a diagram of $L$, so that $\operatorname{braid}(L)\leq n$; the inequality
$\operatorname{braid}(L)\geq 1+(E-e)/2$ is the Morton--Franks--Williams inequality~\cite{Morton1986,FranksWilliams1987}.
If $g\in\vF$, the writhe $\operatorname{wr}(D)$ vanishes \cite[Proposition~3.6]{AielloConti2019b}, and
\cite[Theorem~1]{Morton1986} gives $-(n-1)=\operatorname{wr}(D)-(s(D)-1)\leq e\leq E\leq\operatorname{wr}(D)+(s(D)-1)=n-1$. Taking
the minimum over all such tree diagrams gives the claim.
\end{proof}

\begin{remark}
Proposition~\ref{prop:homfly} and Theorem~\ref{thm:span} are complementary. For the trefoil, $e=2$, $E=4$ and
$M=2$, so that $\ind_{\vF}\geq 5$, while Theorem~\ref{thm:span} only gives $4$; for the figure-eight knot, $e=-2$,
$E=2$ and $M=2$, so that Proposition~\ref{prop:homfly} only gives $3$, while Theorem~\ref{thm:span} gives $5$.
\end{remark}

\FloatBarrier
\section{Methods}\label{sec:computation}

All the code used for this paper, together with the resulting tables, is available in the GitHub repository
\url{https://github.com/valerianoaiello/Thompson-knot-theory}~\cite{Repository}. Its \texttt{README} explains
how to reproduce every computation; the main files are named below. The data archived on Zenodo~\cite{Data}
contain, for each group and each $n$ enumerated, every HOMFLY-PT class of knot diagrams with its number of
diagrams and one representative (a witness when the class contains a knot whose index is $n$), and the final
tables; the full lists of diagrams are not stored, since the code regenerates them. The exhaustive enumerations and the certification
were run on a laptop with 8 cores; the enumeration of $\vF$ with $n=19$ and part of the targeted searches of
Section~\ref{sec:targeted} ran on the server Lagrange of the Dipartimento di Matematica Guido Castelnuovo, Sapienza Universit\`a di Roma, with 12
threads. The present code was developed with the assistance of Claude Code (Anthropic), starting from a first
implementation 
written in 2022--2023; the results were checked as described below.

\subsection{Overview}
For each of the four groups and each $n$ we proceed in four steps.
\begin{enumerate}
\item \emph{Enumeration} of all tree diagrams of size $n$, up to the symmetries below, keeping the prime ones
  (Section~\ref{sec:reductions}) that give a knot, and grouping them by their HOMFLY-PT polynomial.
\item \emph{Identification} of the knots of each HOMFLY-PT class with KnotInfo and SnapPy.
\item \emph{Complete check} of the classes in which a knot may be hidden among the diagrams that were not
  examined in step~2.
\item For $\vF$ and $\vFt$, \emph{oriented scan} of the non-invertible knots, which separates the two
  orientations.
\end{enumerate}
The index of a knot is the first $n$ at which it is realized. The first $n$ at which its HOMFLY-PT class
appears is an exact lower bound, since the HOMFLY-PT polynomial is computed for every diagram. Single cases
beyond the enumeration are settled by targeted searches (Section~\ref{sec:targeted}), and all the results are
then certified (Sections~\ref{sec:certified-tables} and~\ref{sec:certified}).

\subsection{Enumeration}
Two C++ programs enumerate the diagrams: \texttt{src/fknots.cpp} for $F$ and $\vF$ (binary trees), and
\texttt{src/iota\_knots.cpp} for $F_3$ and $\vFt$ (ternary trees). The script \texttt{run.py} runs them in
parallel. For a pair of trees $(T_+,T_-)$ the programs:
\begin{enumerate}
\item build the Tait graph $\Gamma(T_+,T_-)$ described in Section~\ref{sec:construction}, and from it the
  link diagram, as a signed Gauss code;
\item count the components of the link and keep only the knots;
\item simplify the diagram by Reidemeister moves of type I and II;
\item compute the HOMFLY-PT polynomial with the library libhomfly, keeping a cache indexed by the simplified
  Gauss codes.
\end{enumerate}
For the oriented groups only the pairs of trees that lie in the group are generated: by
Section~\ref{sec:construction} these are the pairs with the same parities of right edges at every leaf for $\vF$,
and of the weights $c_\pm$ for $\vFt$, so the trees are grouped by these parities and only pairs within a group
are formed. For $\vF$ with $n=17$ the program builds all the trees once and shares them among the worker
processes: it examines about $10^8$ knot diagrams in 23 minutes. For $\vF$ with $n=19$ (about $3\cdot 10^9$ knot
diagrams; 9 hours and 61\,GB of memory on the server, with 12 threads) it records only the HOMFLY-PT classes
with their number of diagrams, and stores the diagrams of a prescribed list of classes
(Section~\ref{sec:targeted}). Two independent runs gave the same 160\,473 classes. The output is, for every $n$, one line per
HOMFLY-PT class (up to mirror image), with the number of diagrams and up to 50 representatives. A reference
implementation in Python (\texttt{ftrees.py}, \texttt{oriented.py}) is used to test the C++ programs.

\subsection{Reductions}\label{sec:reductions}
The following reductions are used during the enumeration.

\begin{center}
\begin{tabular}{lcccc}
\toprule
reduction & $F$ & $F_3$ & $\vF$ & $\vFt$ \\
\midrule
$g\mapsto g^{-1}$ & yes & yes & yes & yes \\
vertical reflection $(\sigma T_+,\sigma T_-)$ & yes & yes & no & no \\
skip non-prime diagrams & yes & yes & yes & yes \\
knots only for odd $n$ & -- & -- & yes & yes \\
\bottomrule
\end{tabular}

\end{center}

\begin{itemize}
\item \emph{Symmetries.} $(T_+,T_-)\mapsto(T_-,T_+)$ replaces $L(g)$ by its mirror image in $F$ and $F_3$, and
  by the mirror image of its reverse in $\vF$ and $\vFt$ (Theorem~\ref{thm:oriented}\ref{it:inverse}); one diagram of each
  pair is kept. The vertical reflection $(\sigma T_+,\sigma T_-)$ also gives the mirror image, but it does not
  preserve $\vF$ nor $\vFt$. In $\vFt$ it exchanges black and white regions; for $n=4$, 16 of the 28 elements
  of $\vFt$ are mapped outside $\vFt$.
\item \emph{Non-prime diagrams} (Section~\ref{sec:nonprime}). A prime knot realized by a non-prime diagram is
  already realized by the smaller piece in $F$ and $F_3$, and by Lemma~\ref{lem:nonprime} and
  Corollary~\ref{cor:nonprime} the non-prime diagrams can also be skipped in $\vF$ and $\vFt$, so they are skipped
  in all four groups. (The first computation of the $\vFt$ tables skipped only those whose common interval starts
  at an odd leaf; as predicted by Corollary~\ref{cor:nonprime}, none of the others is a witness: all the
  $12\,708$ witnesses of the $\vFt$-index and the $16\,006$ witnesses of the oriented values are prime diagrams,
  and the second computation, which skips all of them, gives the same values (Section~\ref{sec:certified-tables}).
  These non-prime diagrams almost never give knots: for $n=9$ they are $9\,800$ of the $685\,727$ knot diagrams.)
\item \emph{Parity.} In $\vF$ and $\vFt$ knots only occur for odd $n$ (Theorem~\ref{thm:oriented}\ref{it:parity}), so even $n$
  are skipped.
\end{itemize}

\subsection{Identification}
The scripts \texttt{identify.sage} and \texttt{identify\_lib.py} identify the representatives of each
HOMFLY-PT class.
\begin{itemize}
\item A diagram with at most 13 crossings after a further simplification with spherogram, the link diagram module of SnapPy~\cite{SnapPy}, is
  recognized by
  Sage's function \texttt{get\_knotinfo} from KnotInfo~\cite{KnotInfo}. Its knot then has crossing number at
  most 13 and belongs to the table, so the identification is a proof.
\item A larger diagram is identified by SnapPy~\cite{SnapPy} from its complement (census up to 15
  crossings), which determines the knot up to mirror image by the theorem of Gordon and
  Luecke~\cite{GordonLuecke1989}. A HOMFLY-PT match alone is not enough on larger diagrams: at $n=10$ a class
  has the polynomial of $13n_{2638}$, but its knot is $14n_{2554}$.
\item All the stored representatives of a class are tried, since one class can contain several knots (for
  instance $10_{129}$ and $8_8$).
\end{itemize}
SnapPy crashes on a few diagrams; the computation runs in a pool of processes that survives the crashes and
reports every diagram it could not examine.

\subsection{Complete check}
Since only 50 representatives of each HOMFLY-PT class are stored, a knot can hide among the unstored
diagrams of a large class. Without the check described below some indices would be overestimated: for instance $10_{132}$ is realized at $n=9$
only by diagrams that are not among the 50 stored representatives of its class. The script \texttt{complete\_check.py} therefore considers every knot whose class appears
before the knot is found, or which is not found at all. For such a knot it enumerates again \emph{all} the
diagrams of its class, from the first appearance of the class up to the index, and compares each of them
with the knot: first the hyperbolic volume, then a SnapPy isometry. Diagrams without a hyperbolic structure are
treated by \texttt{resolve\_undecided.py}:
\begin{itemize}
\item a different knot Floer homology (computed with Szab\'o's program~\cite{SzaboHFK}, through spherogram)
  shows that the diagram is not the knot;
\item otherwise an essential torus (SnapPy's \texttt{splitting\_surfaces}) shows that it is a satellite knot,
  whereas the knots concerned are hyperbolic.
\end{itemize}
No diagram is left undecided. All the conclusions of the complete check were then checked again with exact or
certified methods only (Section~\ref{sec:certified-tables}).

\subsection{Certified verification of the tables}\label{sec:certified-tables}
The volumes and isometries computed by SnapPy in floating point are very reliable, but they are not proofs. We
therefore checked every value of the four tables again with exact or certified methods only
(scripts \texttt{certify\_indices.py} and \texttt{certify\_retry.py}). A computation is \emph{exact} when it only involves integers and polynomials (HOMFLY-PT polynomials, knot
Floer homology, triangulations and their isomorphism signatures, normal surfaces), and \emph{certified} when it
uses interval arithmetic: every real quantity is represented by an interval that provably contains its true value,
rounding errors included. In this way SnapPy proves that a hyperbolic structure exists, encloses its volume in an
interval and certifies the canonical retriangulation of the complement~\cite{HIKMOT2016,SnapPy}. A value $\ind(K)=n$ rests on two kinds of statements.

\emph{Lower bound.} No tree diagram of size smaller than $n$ realizes $K$. If the HOMFLY-PT class of $K$ does not occur
before $n$, this follows from the enumeration, since the HOMFLY-PT polynomial is computed exactly for every
diagram; for $F$ and $\vF$ it also follows from Theorem~\ref{thm:span} (for $\vF$ raised to the next odd
integer, by Theorem~\ref{thm:oriented}\ref{it:parity}) whenever $\ind(K)$ equals the bound, which happens for 1\,440 knots in $F$ and for 8 knots in $\vF$, without any enumeration.
Otherwise every tree diagram of the class of size between its first appearance and $n-1$ (or the largest size
enumerated, if $K$ is not found) is shown not to be $K$, by the first of the following arguments that applies:
\begin{enumerate}
\item after simplification by Reidemeister moves it has fewer than $c(K)$ crossings, so its knot has crossing
  number smaller than $c(K)$;
\item it is a connected sum of two diagrams with non-trivial Alexander polynomial, so its knot is composite,
  whereas $K$ is prime;
\item its knot Floer homology differs from that of $K$ up to mirror image;
\item its hyperbolic volume, computed in interval arithmetic with a verified hyperbolic structure, is disjoint
  from that of $K$;
\item the verified canonical retriangulations of the two complements are not combinatorially isomorphic;
\item it has no hyperbolic structure and Regina~\cite{Regina} finds in its complement a closed normal torus that is
  incompressible and not boundary parallel, so its knot is a satellite, whereas $K$ is hyperbolic.
\end{enumerate}

\emph{Witness.} A tree diagram of size $n$ realizes $K$. This is shown by a combinatorial isomorphism between the verified
canonical retriangulations of the two complements or, when these are too slow to compute or $K$ is not
hyperbolic, by random retriangulations of the two complements reaching the same triangulation (equal isomorphism
signatures). In both cases the complements are homeomorphic, hence the knots coincide up to mirror image by the
theorem of Gordon and Luecke~\cite{GordonLuecke1989}. For a few torus knots the witness has at most 13 crossings
and is recognized in KnotInfo.

\begin{table}[htbp]
\centering
\begin{tabular}{lrrrr}
\toprule
 & $F$ & $F_3$ & $\vF$ & $\vFt$ \\
\midrule
witnesses from the enumeration certified & 6\,049 & 11\,320 & 2\,155 & 12\,708 \\
\quad by canonical triangulations & 6\,036 & 11\,303 & 2\,152 & 12\,690 \\
\quad by isomorphism signatures & 6 & 10 & 1 & 11 \\
\quad by KnotInfo (torus knots) & 7 & 7 & 2 & 7 \\
diagrams excluded & 27\,857 & 26\,847 & 390\,776 & 14\,265 \\
\quad by the crossing number & 24\,257 & 14\,402 & 367\,460 & 4\,323 \\
\quad as composite knots & 5 & 1 & 218 & 73 \\
\quad by knot Floer homology & 438 & 1\,078 & 8\,860 & 708 \\
\quad by the certified volume & 2\,467 & 10\,679 & 12\,650 & 8\,517 \\
\quad by canonical triangulations & 690 & 687 & 1\,588 & 644 \\
contradictions with the tables & 0 & 0 & 0 & 0 \\
\bottomrule
\end{tabular}
\end{table}

The witnesses found by the targeted searches of Section~\ref{sec:targeted} are certified in the same way.
A final check (\texttt{certify\_completeness.py}) confirms that every knot with a known index has a certified
witness, that every knot whose index exceeds the first appearance of its HOMFLY-PT class is in the complete check,
and that all the diagrams of these classes in the relevant range have been excluded. Hence all the values of the
four tables are certified, the only non-certified ingredient being the exact enumeration of the diagrams and of
their HOMFLY-PT polynomials.

\paragraph{Second computation.} All four tables were computed a second time, from scratch, with the reductions
of Section~\ref{sec:reductions} (all non-prime diagrams skipped in the four groups) and, for $F_3$, a local
$RIII$--$RII$ reduction (a ternary diagram whose signed Tait graph admits a Reidemeister III move followed by the
cancellation of a pair of edges gives the same link as a diagram with one caret less, and is skipped); the second
pass on the undecided diagrams used exact tests only. The indices, the
lower bounds and the oriented values of every knot agree with the first computation; some witnesses differ (51
for $F_3$, 37 for $\vFt$, 2 for $F$, 1 for $\vF$), and all of them are certified. The figures above refer to the
second computation. For $F_3$ the $RIII$--$RII$ reduction is only a minor speed-up (at $n=8$ it skips 27\% of
the knot diagrams and a quarter of the enumeration time, the later phases being per HOMFLY-PT class); it leaves
the 59\,076 HOMFLY-PT classes at $n=9$ unchanged, and the first computation, which does not use it, gives the
same tables. On an 8-core laptop the whole computation, certification included, takes 1~h~18~min for $F$
($n\le13$), 1~h~42~min for $F_3$ ($n\le9$), 4~h~50~min for $\vF$ ($n\le17$) and 2~h~13~min for $\vFt$ ($n\le11$).

\subsection{Oriented indices}
The script \texttt{oriented\_scan.py} computes the oriented indices of the non-invertible knots.
\begin{enumerate}
\item Every diagram of the HOMFLY-PT class of $K$ is oriented by the colouring of its Tait graph
  (\texttt{oriented.py}).
\item The diagram is compared with the oriented KnotInfo diagram of $K$ through the isometries of the
  complements and their action on meridian and longitude. This tells whether $g$ realizes $K$, $K^r$, $K^m$
  or $K^{mr}$; then $g^{-1}$ realizes $K^{mr}$, $K^m$, $K^r$ or $K$ respectively
  (Theorem~\ref{thm:oriented}\ref{it:inverse}).
\item A knot is done when both $K$ and $K^r$ have been realized.
\end{enumerate}
The rule for $g^{-1}$ was also checked directly, on 17\,740 diagrams.

\subsection{Certified verification of the oriented results}\label{sec:certified}
The results on non-invertible knots rest on statements of the form ``no diagram with fewer than $b$ leaves
realizes $K^r$''. The volumes and isometries computed by SnapPy in floating point are very reliable, but they
are not proofs. We therefore checked all these statements again with exact or certified methods only
(\texttt{verify\_oriented.py}). For each knot $K$ with $\ind(K)\ne\ind(K^r)$, every diagram of the HOMFLY-PT
class of $K$ with $n$ smaller than the index of the missing orientation is compared with $K$ as follows:
\begin{enumerate}
\item if its knot Floer homology (exact) differs from that of $K$ up to mirror image, it is another knot;
\item otherwise, if its hyperbolic volume, computed in interval arithmetic with a verified hyperbolic
  structure, is disjoint from that of $K$, it is another knot;
\item otherwise the verified canonical retriangulations of the two complements and all their combinatorial
  isomorphisms determine exactly whether the diagram realizes $K$, $K^r$, $K^m$ or $K^{mr}$.
\end{enumerate}
The witnesses of the realized orientations are certified in the same way (step 3).

\begin{center}
\begin{tabular}{lcc}
\toprule
 & $\vF$ & $\vFt$ \\
\midrule
knots with $\ind(K)\ne\ind(K^r)$ & 311 & 946 \\
witnesses certified & 396 & 1\,656 \\
diagrams excluded by knot Floer homology & 619 & 17\,774 \\
diagrams excluded by the certified volume & 195\,981 & 11\,390 \\
diagrams excluded by the canonical triangulations & 1\,288 & 516 \\
diagrams of $K$ certified in the allowed orientation & 1\,093 & 1\,855 \\
diagrams not certified & 0 & 0 \\
\bottomrule
\end{tabular}
\end{center}

All the 311 results for $\vF$ and all the 946 results for $\vFt$ are therefore certified. The exact value $19$
of 170 of the knots marked $19^*$ (and $13$ for 212 knots in $\vFt$) also uses
Theorem~\ref{thm:oriented}\ref{it:reverse} (resp.~\ref{it:reverse3}), whose moves consist of Reidemeister moves;
for the other 56 (resp.\ 24) knots it uses a diagram found by a targeted search, certified in the same way
(Section~\ref{sec:targeted}).

\subsection{Targeted searches}\label{sec:targeted}
Beyond the exhaustive enumeration, single open cases can be attacked by a local search
(\texttt{descent\_search.py}). Starting from known diagrams of a knot $K$ of size $n$, one deletes $s+1$ internal
vertices and inserts one new vertex in each tree, with $s=1$ in $F$ and $F_3$ and $s=2$ in $\vF$ and $\vFt$, so
that the size decreases by $s$ (a vertex can be deleted when all its children but at most one are leaves, and it
is replaced by that child; a new vertex is inserted above any subtree, which becomes one of its children). All the pairs of trees obtained in this way that lie in
the group are filtered by the determinant (computed from the signed Tait graph), by the number of components
and by the hyperbolic volume, and the remaining ones are compared with $K$ by SnapPy. Every diagram found is then
certified as in Sections~\ref{sec:certified-tables} and~\ref{sec:certified} (\texttt{certify\_descent.py}):
4\,569 diagrams, all certified. A diagram found proves an upper bound; not finding one proves nothing. Where
this search finds nothing, a deeper one deletes $s+2$ vertices and inserts two in each tree. The searches gave:
\begin{itemize}
\item $F_3$: the 14 knots with 12 crossings not realized with $n\le 9$ have a diagram of size $11$ in
  $\vFt\subset F_3$, hence $F_3$-index $10$ or $11$; starting from these diagrams (from all the stored ones of
  size $11$ in a second round), the search found diagrams of size $10$ for 12 of them, whose $F_3$-index is
  therefore $10$: $12a_{868}$, $12a_{875}$, $12a_{928}$, $12a_{1019}$, $12a_{1020}$, $12a_{1050}$, $12a_{1113}$,
  $12a_{1152}$, $12a_{1196}$, $12a_{1209}$, $12a_{1285}$, $12a_{1286}$, and the deeper search found diagrams of
  size $10$ for the last two, $12a_{883}$ and $12a_{975}$. Hence the $F_3$-index of every knot with at most 12
  crossings is known. In the same way, of the 1\,467 knots with 13 crossings not realized with $n\le 9$ that
  have a $\vFt$-witness of size $11$ and lower bound $10$, 897 have a diagram of size $10$, hence $F_3$-index
  $10$: 307 in a first round starting from the witness, 337 in a second round starting from at most five other
  stored diagrams of size $11$ of the knot, and 253 in a third round starting from all the remaining ones (for the
  other 570 the $F_3$-index is still $10$ or $11$);
\item $\vFt$: starting from $\iota$ of their $\vF$-witnesses (size $17$), $\ind_{\vFt}(12a_{1211})=13$ and, with
  the deeper search, $\ind_{\vFt}(12a_{1249})=13$. The knot $12a_{1202}$ has no known diagram in $\vF$; its
  $F_3$-index is $9$, and inserting four vertices in each tree of its diagrams of size $9$ in $F_3$ (all the
  pairs of trees obtained that lie in $\vFt$) gave a diagram of size $13$, so $\ind_{\vFt}(12a_{1202})=13$.
  The non-invertible knots $12a_{1211}$ and $12a_{1249}$ have the same index in both orientations: the
  diagrams of size $13$ of $12a_{1249}$ realize both, and for $12a_{1211}$, whose diagrams of size $13$ realize
  only $K^r=K^m$, the move with two new vertices of Theorem~\ref{thm:oriented}\ref{it:reverse3} gives $K$ at
  size $15$ and one step of the search gives it at size $13$.
  Hence the $\vFt$-index of every knot with at most 12 crossings is known, in both orientations;
\item $F$: the deeper search found diagrams of size $14$ for $11a_{249}$ (whose upper bound was $17$) and for
  $11a_{352}$ (upper bound $15$, by its $\vF$-witness), so $\ind_F(11a_{249})=\ind_F(11a_{352})=14$ and the
  $F$-index of every knot with at most 11 crossings is known;
\item the other orientation of the 56 knots of $\vF$ and of the 24 knots of $\vFt$ above (starting from the
  diagram of size $n+4$ given by Theorem~\ref{thm:oriented}\ref{it:reverse}, resp.~\ref{it:reverse3}): a
  diagram of size $n+2$ for all of them;
\item $\vF$, knots with at most 10 crossings: before these searches 24 knots with at most 9 crossings and 62
  with 10 crossings had no diagram with at most 17 leaves. The enumeration of $\vF$ with $n=19$
  (Section~\ref{sec:computation}) stored the diagrams of their HOMFLY-PT classes. The class is present at
  $n=19$ for 30 of them, and in all cases a diagram with 19 leaves realizes the knot: the $\vF$-index is $19$
  for $7_3$, $7_5$, $9_{11}$, $9_{20}$, $9_{36}$, $9_{49}$ and for 24 knots with 10 crossings. Three of them, $10_{67}$,
  $10_{85}$ and $10_{109}$, are non-invertible, and the oriented index of the other orientation is not computed
  (it is $19$, $21$ or $23$ by Theorem~\ref{thm:oriented}\ref{it:reverse}). For the other 56
  the class is absent, so their index is at least $21$ (Theorem~\ref{thm:oriented}\ref{it:parity}). In all,
  the class is absent at $n=19$ for 6\,908 of the 10\,810 knots with at most 13 crossings not realized with
  $n\le 17$, whose $\vF$-index is therefore at least $21$. Since these knots have
  no known diagram in $\vF$, we produced one with the proof of the Alexander-type theorem of~\cite{Aiello2020}
  turned into an algorithm, and then reduced its size by deleting vertices (greedy steps checked by SnapPy,
  alternated with one step of the search above). This gave diagrams with 21 leaves, hence $\vF$-index $21$,
  for $7_1$, $8_{15}$, $8_{19}$, $9_2$, $9_5$ and $10_{25}$. All the 36 diagrams were certified as above. For
  the remaining 13 knots with 9 crossings and 37 with 10 crossings the $\vF$-index is at least $21$; the
  descent stopped at 23 leaves for $9_7$, $9_{18}$, $9_{23}$, $10_6$ and $10_{14}$, even when started from
  several different constructions.
\end{itemize}

\subsection{Verification}
\begin{itemize}
\item The C++ programs were compared with the Python implementation, with Sage, and with an earlier
  implementation by the author. The ternary program, applied to $\iota(F)$, was compared with the binary one.
\item The orientation of the links of $\vFt$ was compared with the C++ program (same knot), with $\iota$ (same
  oriented knot), and with the rule for $g^{-1}$.
\item The final tables satisfy all the inequalities of Theorems~\ref{thm:unoriented} and~\ref{thm:oriented}. For every $\vF$
  witness $g$, $\iota(g)$ was checked to be a diagram of $\vFt$ of the same knot.
\item Every diagram drawn in this paper was compared with the data by SnapPy.
\end{itemize}

\section{Results}\label{sec:results}

\begin{table}[tbp]
\centering
\begin{tabular}{c|cccc}
\toprule
crossings & $F$ & $F_3$ & $\vec F$ & $\vec F_3$ \\
 & $n\le 13$ & $n\le 9$ & $n\le 17$ & $n\le 11$ \\
\midrule
3 & 1 / 1 & 1 / 1 & 1 / 1 & 1 / 1 \\
4 & 1 / 1 & 1 / 1 & 1 / 1 & 1 / 1 \\
5 & 2 / 2 & 2 / 2 & 2 / 2 & 2 / 2 \\
6 & 3 / 3 & 3 / 3 & 3 / 3 & 3 / 3 \\
7 & 7 / 7 & 7 / 7 & 7 / 7 & 7 / 7 \\
8 & 21 / 21 & 21 / 21 & 21 / 21 & 21 / 21 \\
9 & 49 / 49 & 49 / 49 & 36 / 49 & 49 / 49 \\
10 & 165 / 165 & 165 / 165 & 128 / 165 & 165 / 165 \\
11 & 552 / 552 & 552 / 552 & 278 / 552 & 552 / 552 \\
12 & 1348 / 2176 & 2176 / 2176 & 570 / 2176 & 2176 / 2176 \\
13 & 3902 / 9988 & 9254 / 9988 & 1144 / 9988 & 9734 / 9988 \\
\bottomrule
\end{tabular}

\caption{Number of knots with known index / number of prime knots, by crossing number, with the largest $n$
enumerated in each group. For the oriented groups the index of a knot is known when it is known for at least
one orientation. The counts include the knots found beyond the enumeration by the targeted searches of
Section~\ref{sec:targeted} (911 knots with $F_3$-index $10$, three with $\vFt$-index $13$, two with $F$-index $14$, 30 with
$\vF$-index $19$ and six with $\vF$-index $21$).}
\label{tab:summary}
\end{table}

The four indices of all the prime knots with at most 10 crossings are listed in Table~\ref{tab:small}
(Appendix~\ref{app:tables}).

\subsection{Non-invertible knots}\label{sec:oriented}
Among the 1\,308 non-invertible knots with at most 13 crossings whose HOMFLY-PT class occurs in $\vF$ with
at most 17 leaves, 997 have the same $\vF$-index in both orientations, 85 have both orientations realized
with different indices ($13$ and $15$ for $10_{88}$, $12n_{18}$, $12n_{49}$ and two more; $15$ and $17$ for
80 knots, the first being $10_{93}$ and $10_{102}$), and 226 have one orientation realized with at most 17
leaves and the other not (e.g.\ $10_{95}$). In all cases the realized orientation has index $17$. For 170 of
these 226 knots a minimal diagram of that orientation has first and last Tait vertices of different colours, so
by Theorem~\ref{thm:oriented}\ref{it:reverse} the other orientation has index exactly $19$. For the remaining 56
knots (the first ones being $12a_8$, $12a_{131}$, $12a_{328}$, $12a_{818}$, $12a_{951}$), every minimal diagram of
the realized orientation has first and last vertices of the same colour, so the theorem only gives $\le 21$;
for all of them a targeted search (Section~\ref{sec:targeted}) found a diagram with 19 leaves realizing the other
orientation, which therefore also has index exactly $19$.

This answers the two questions of~\cite{Aiello2026}: the $\vF$-indices of $K$ and $K^r$ can be different, the
smallest crossing number of such a knot being $10$ ($10_{88}$, $10_{93}$, $10_{95}$, $10_{102}$); and not every
knot has minimal diagrams of both types (the 56 knots above). Over all the knot diagrams of $\vF$ with
$n\le 13$ stored by the enumeration, 1\,104 have first and last vertices of different colours and 756 of the
same colour. All these results are certified by exact or interval-arithmetic methods (Section~\ref{sec:certified}). All the differences we found are equal to $2$, in $\vF$ as well as in $\vFt$. In $\vFt$, with $n\le 11$, 7\,175 of 8\,121 knots have equal
indices, 710 have different indices ($7$ and $9$ for $8_{17}$, $10_{147}$, $10_{149}$, $10_{150}$, $10_{153}$
and five more; $9$ and $11$ for 700 knots) and 236 have one orientation $>11$. For these 236 knots the realized
orientation has index $11$; for 212 of them some minimal diagram of the realized orientation has first and last
Tait vertices of different colours (for 182 the first one found, for 30 another one found by searching all the
minimal diagrams), so by Theorem~\ref{thm:oriented}\ref{it:reverse3} the other orientation has index exactly
$13$. For the remaining 24 every minimal diagram has ends of the same colour, so the theorem only gives
$\le 15$; for all of them a targeted search found a diagram of size $13$ realizing the other orientation, which
therefore also has index exactly $13$. The complete lists are in Tables~\ref{tab:oriented-Fvec}
and~\ref{tab:oriented-F3vec} (Appendix~\ref{app:tables}).

\subsection{Mutants}
A mutant of a knot $K$ is obtained by cutting $K$ along a Conway sphere, a sphere meeting $K$ in four points,
and regluing the ball it bounds after one of the three rotations by $\pi$ that preserve the four points. Mutation
preserves many invariants, among them the HOMFLY-PT and Kauffman polynomials, hence the Jones and Alexander
polynomials. In particular two mutants lie in the same HOMFLY-PT class, so they cannot be told apart by
the first $n$ at which this class occurs, and it is natural to ask whether the Thompson indices distinguish them.
The best-known pair is formed by the Conway knot $11n_{34}$ and the Kinoshita--Terasaka knot $11n_{42}$.

\begin{figure}[ht]
\centering
\resizebox{0.45\linewidth}{!}{\begin{tikzpicture}[scale=0.32, line cap=round]
  \draw (1.5,1) -- (1,0.5);
  \draw (1.5,1) -- (2,0.5);
  \draw (1,1.5) -- (0,0.5);
  \draw (1,1.5) -- (1.5,1);
  \draw (5.5,1) -- (5,0.5);
  \draw (5.5,1) -- (6,0.5);
  \draw (5,1.5) -- (4,0.5);
  \draw (5,1.5) -- (5.5,1);
  \draw (7.5,1) -- (7,0.5);
  \draw (7.5,1) -- (8,0.5);
  \draw (6,2.5) -- (5,1.5);
  \draw (6,2.5) -- (7.5,1);
  \draw (6.5,3) -- (6,2.5);
  \draw (6.5,3) -- (9,0.5);
  \draw (6,3.5) -- (3,0.5);
  \draw (6,3.5) -- (6.5,3);
  \draw (6.5,4) -- (6,3.5);
  \draw (6.5,4) -- (10,0.5);
  \draw (5,5.5) -- (1,1.5);
  \draw (5,5.5) -- (6.5,4);
  \draw (11.5,1) -- (11,0.5);
  \draw (11.5,1) -- (12,0.5);
  \draw (13.5,1) -- (13,0.5);
  \draw (13.5,1) -- (14,0.5);
  \draw (12.5,2) -- (11.5,1);
  \draw (12.5,2) -- (13.5,1);
  \draw (7,7.5) -- (5,5.5);
  \draw (7,7.5) -- (12.5,2);
  \draw (3.5,-1) -- (3,-0.5);
  \draw (3.5,-1) -- (4,-0.5);
  \draw (3,-1.5) -- (2,-0.5);
  \draw (3,-1.5) -- (3.5,-1);
  \draw (6.5,-1) -- (6,-0.5);
  \draw (6.5,-1) -- (7,-0.5);
  \draw (6,-1.5) -- (5,-0.5);
  \draw (6,-1.5) -- (6.5,-1);
  \draw (9.5,-1) -- (9,-0.5);
  \draw (9.5,-1) -- (10,-0.5);
  \draw (9,-1.5) -- (8,-0.5);
  \draw (9,-1.5) -- (9.5,-1);
  \draw (9.5,-2) -- (9,-1.5);
  \draw (9.5,-2) -- (11,-0.5);
  \draw (8,-3.5) -- (6,-1.5);
  \draw (8,-3.5) -- (9.5,-2);
  \draw (8.5,-4) -- (8,-3.5);
  \draw (8.5,-4) -- (12,-0.5);
  \draw (7,-5.5) -- (3,-1.5);
  \draw (7,-5.5) -- (8.5,-4);
  \draw (6.5,-6) -- (1,-0.5);
  \draw (6.5,-6) -- (7,-5.5);
  \draw (7,-6.5) -- (6.5,-6);
  \draw (7,-6.5) -- (13,-0.5);
  \draw (6.5,-7) -- (0,-0.5);
  \draw (6.5,-7) -- (7,-6.5);
  \draw (7,-7.5) -- (6.5,-7);
  \draw (7,-7.5) -- (14,-0.5);
  \fill (0,0.5) circle (2.5pt) (0,-0.5) circle (2.5pt);
  \fill (1,0.5) circle (2.5pt) (1,-0.5) circle (2.5pt);
  \fill (2,0.5) circle (2.5pt) (2,-0.5) circle (2.5pt);
  \fill (3,0.5) circle (2.5pt) (3,-0.5) circle (2.5pt);
  \fill (4,0.5) circle (2.5pt) (4,-0.5) circle (2.5pt);
  \fill (5,0.5) circle (2.5pt) (5,-0.5) circle (2.5pt);
  \fill (6,0.5) circle (2.5pt) (6,-0.5) circle (2.5pt);
  \fill (7,0.5) circle (2.5pt) (7,-0.5) circle (2.5pt);
  \fill (8,0.5) circle (2.5pt) (8,-0.5) circle (2.5pt);
  \fill (9,0.5) circle (2.5pt) (9,-0.5) circle (2.5pt);
  \fill (10,0.5) circle (2.5pt) (10,-0.5) circle (2.5pt);
  \fill (11,0.5) circle (2.5pt) (11,-0.5) circle (2.5pt);
  \fill (12,0.5) circle (2.5pt) (12,-0.5) circle (2.5pt);
  \fill (13,0.5) circle (2.5pt) (13,-0.5) circle (2.5pt);
  \fill (14,0.5) circle (2.5pt) (14,-0.5) circle (2.5pt);
\end{tikzpicture}}\quad
\resizebox{0.5\linewidth}{!}{\input{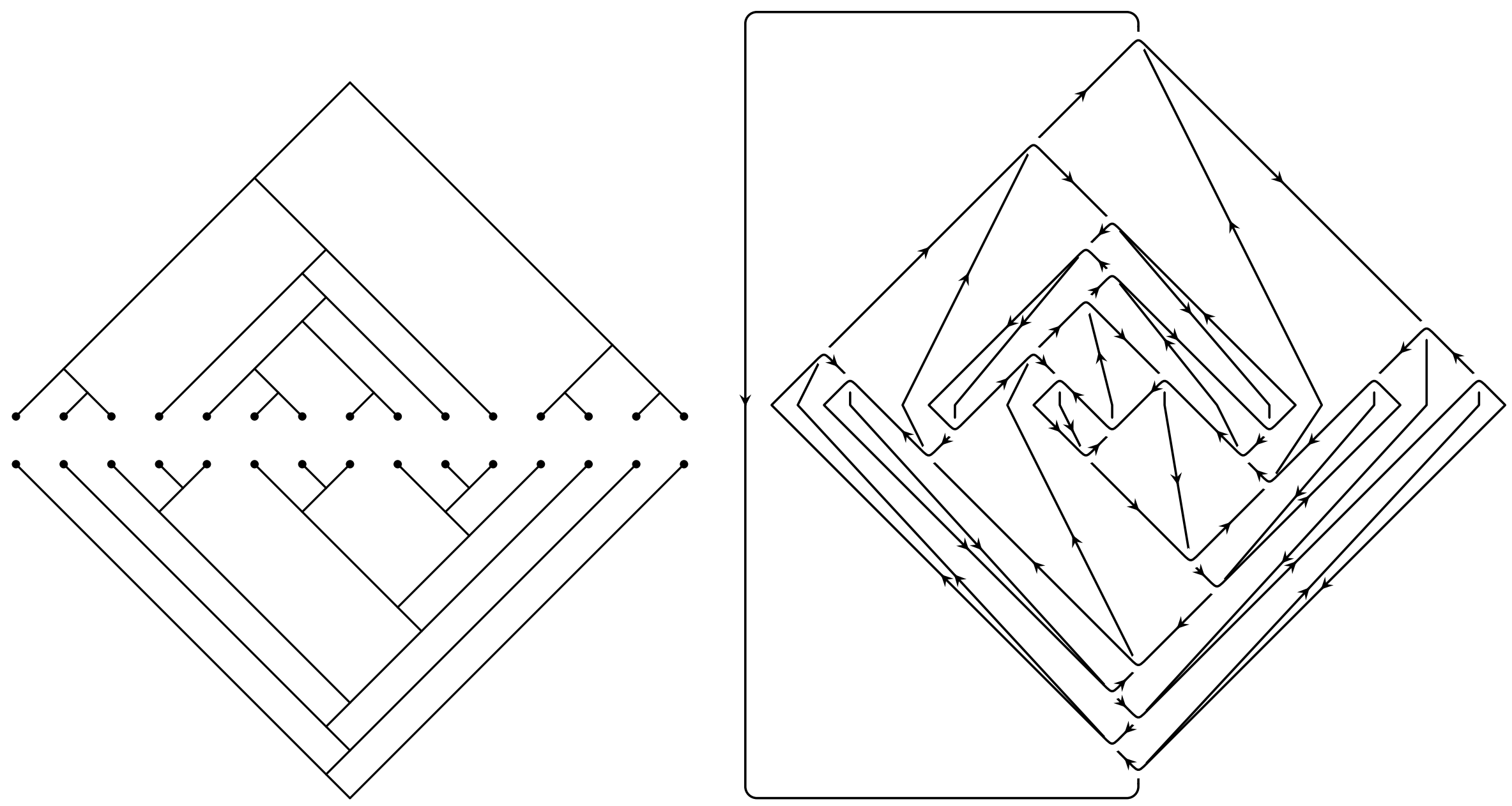}}
\caption{A diagram of $\vF$ with 15 leaves of the Kinoshita--Terasaka knot $11n_{42}$ (orientation
$K^{mr}$).}
\label{fig:mutant}
\end{figure}

\begin{theorem}\label{thm:mutants}
The Conway knot $11n_{34}$ and the Kinoshita--Terasaka knot $11n_{42}$ have different $\vF$-index, in every
orientation: $\ind_{\vF}(11n_{42})=15$ and $\ind_{\vF}(11n_{34})\ge 19$.
\end{theorem}
\begin{proof}
A diagram of $11n_{42}$ (in fact of $11n_{42}^{mr}$) with 15 leaves is shown in Figure~\ref{fig:mutant}. The
common HOMFLY-PT class of the two knots does not occur in $\vF$ with fewer than 15 leaves; it contains 39
diagrams with 15 leaves and 2\,541 with 17. Each of these diagrams represents a knot of Seifert genus $2$, as
computed via knot Floer homology, whereas $11n_{34}$ has Seifert genus $3$; hence none of them represents
$11n_{34}$, with any orientation. Since knots only occur for odd $n$ (Theorem~\ref{thm:oriented}\ref{it:parity}),
$\ind_{\vF}(11n_{34})\ge 19$ for every orientation and its mirror image.
\end{proof}
The $F$-, $F_3$- and $\vFt$-indices of the two knots coincide ($12$, $8$ and $9$).

More generally, Stoimenow tabulated the mutant knots up to 16 crossings~\cite{Stoimenow2010}. Up to 13
crossings his tables contain 865 groups of mutants (789 pairs, 43 triples, 32 groups of four and one of six
knots; counting the pairs within the groups, 1\,125 pairs); the 16 pairs with 11 crossings are also listed
in~\cite{MMMRS2015}. We converted his KnotScape numbers into KnotInfo names and checked every conversion with
SnapPy, from the Dowker--Thistlethwaite codes of his tables. Independently, the pairs of distinct knots with at
most 13 crossings having the same HOMFLY-PT polynomial (up to mirror image) and the same hyperbolic volume are
exactly these 1\,125 pairs.

Of these mutant pairs, 458 are distinguished by at least one index: 122 by the $F$-index, 176 by the
$F_3$-index, 97 by the $\vF$-index and 158 by the $\vFt$-index (Table~\ref{tab:mutants}). With 11 crossings,
12 of the 16 pairs are distinguished; for instance $11n_{35}$ and $11n_{43}$ have $F$-index $11$ and $12$, and
$11a_{24}$ and $11a_{26}$ have $\vFt$-index $9$ and $11$. The pairs $11n_{40}$/$11n_{46}$, $11n_{41}$/$11n_{47}$,
$11n_{71}$/$11n_{75}$ and $11n_{76}$/$11n_{78}$ are not distinguished by the computed values.

\subsection{The Borromean rings}\label{sec:borromean}
Jones observed that the $F$-index of the Borromean rings is at most $20$, and ``surely considerably
less''~\cite[Section~5.3]{Jones2017}, and posed the problem of computing their $F_3$- and
$F_2$-indices~\cite[Section~7]{Jones2019}. There the size of a tree diagram is the number of internal vertices of
each tree, that is $n-1$ in our normalization. For $F$ the question was answered by Sweeney~\cite{Sweeney2021},
who found diagrams of the Borromean rings with 9 leaves; we confirm his answer and answer the question for $F_3$.

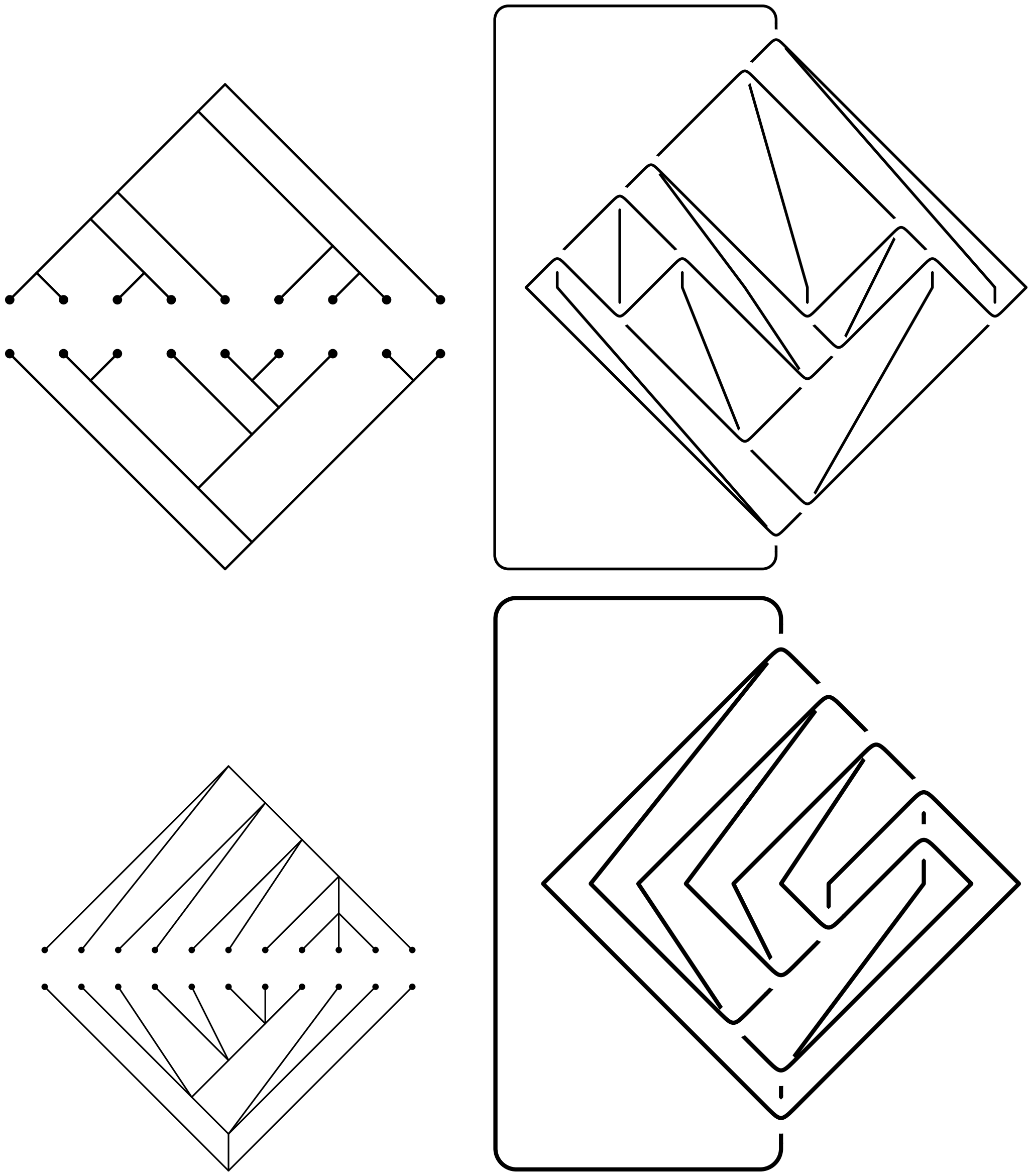
\begin{figure}[tbp]
\centering
\begin{tabular}{cc}
\resizebox{0.38\linewidth}{!}{\begin{tikzpicture}[scale=0.32, line cap=round]
  \draw (0.5,1) -- (0,0.5);
  \draw (0.5,1) -- (1,0.5);
  \draw (2.5,1) -- (2,0.5);
  \draw (2.5,1) -- (3,0.5);
  \draw (1.5,2) -- (0.5,1);
  \draw (1.5,2) -- (2.5,1);
  \draw (2,2.5) -- (1.5,2);
  \draw (2,2.5) -- (4,0.5);
  \draw (6.5,1) -- (6,0.5);
  \draw (6.5,1) -- (7,0.5);
  \draw (6,1.5) -- (5,0.5);
  \draw (6,1.5) -- (6.5,1);
  \draw (3.5,4) -- (2,2.5);
  \draw (3.5,4) -- (6,1.5);
  \draw (4,4.5) -- (3.5,4);
  \draw (4,4.5) -- (8,0.5);
  \draw (1.5,-1) -- (1,-0.5);
  \draw (1.5,-1) -- (2,-0.5);
  \draw (4.5,-1) -- (4,-0.5);
  \draw (4.5,-1) -- (5,-0.5);
  \draw (5,-1.5) -- (4.5,-1);
  \draw (5,-1.5) -- (6,-0.5);
  \draw (4.5,-2) -- (3,-0.5);
  \draw (4.5,-2) -- (5,-1.5);
  \draw (3.5,-3) -- (1.5,-1);
  \draw (3.5,-3) -- (4.5,-2);
  \draw (7.5,-1) -- (7,-0.5);
  \draw (7.5,-1) -- (8,-0.5);
  \draw (4.5,-4) -- (3.5,-3);
  \draw (4.5,-4) -- (7.5,-1);
  \draw (4,-4.5) -- (0,-0.5);
  \draw (4,-4.5) -- (4.5,-4);
  \fill (0,0.5) circle (2.5pt) (0,-0.5) circle (2.5pt);
  \fill (1,0.5) circle (2.5pt) (1,-0.5) circle (2.5pt);
  \fill (2,0.5) circle (2.5pt) (2,-0.5) circle (2.5pt);
  \fill (3,0.5) circle (2.5pt) (3,-0.5) circle (2.5pt);
  \fill (4,0.5) circle (2.5pt) (4,-0.5) circle (2.5pt);
  \fill (5,0.5) circle (2.5pt) (5,-0.5) circle (2.5pt);
  \fill (6,0.5) circle (2.5pt) (6,-0.5) circle (2.5pt);
  \fill (7,0.5) circle (2.5pt) (7,-0.5) circle (2.5pt);
  \fill (8,0.5) circle (2.5pt) (8,-0.5) circle (2.5pt);
\end{tikzpicture}} &
\resizebox{0.45\linewidth}{!}{\begin{tikzpicture}[scale=0.32, line cap=round, line join=round, >=stealth, thick]
  \draw[rounded corners=4pt] (7.682,7.682) -- (7.318,7.318);
  \draw[rounded corners=4pt] (8.296,7.661) -- (14.704,0.339) .. controls (15,0) .. (15,-0.45) -- (15,-0.55);
  \draw[rounded corners=4pt] (8.318,7.682) -- (15.682,0.318) .. controls (16,0) .. (15.682,-0.318) -- (15.318,-0.682);
  \draw[rounded corners=4pt] (6.682,6.682) -- (4.318,4.318);
  \draw[rounded corners=4pt] (7.124,6.567) -- (8.876,0.433) .. controls (9,0) .. (9,-0.45) -- (9,-0.55);
  \draw[rounded corners=4pt] (7.318,6.682) -- (11.682,2.318);
  \draw[rounded corners=4pt] (3.682,3.682) -- (3.318,3.318);
  \draw[rounded corners=4pt] (4.27,3.64) -- (6.73,0.36) .. controls (7,0) .. (7.25,-0.374) -- (8.75,-2.626);
  \draw[rounded corners=4pt] (4.318,3.682) -- (7.682,0.318) .. controls (8,0) .. (8.318,-0.318) -- (8.682,-0.682);
  \draw[rounded corners=4pt] (2.682,2.682) -- (1.318,1.318);
  \draw[rounded corners=4pt] (3,2.55) -- (3,0.45) .. controls (3,0) .. (3,-0.45) -- (3,-0.55);
  \draw[rounded corners=4pt] (3.318,2.682) -- (4.682,1.318);
  \draw[rounded corners=4pt] (0.682,0.682) -- (0.318,0.318) .. controls (0,0) .. (0.318,-0.318) -- (7.682,-7.682);
  \draw[rounded corners=4pt] (1,0.55) -- (1,0.45) .. controls (1,0) .. (1.296,-0.339) -- (7.704,-7.661);
  \draw[rounded corners=4pt] (1.318,0.682) -- (1.682,0.318) .. controls (2,0) .. (2.318,-0.318) -- (2.682,-0.682);
  \draw[rounded corners=4pt] (4.682,0.682) -- (4.318,0.318) .. controls (4,0) .. (3.682,-0.318) -- (3.318,-0.682);
  \draw[rounded corners=4pt] (5,0.55) -- (5,0.45) .. controls (5,0) .. (5.167,-0.418) -- (6.833,-4.582);
  \draw[rounded corners=4pt] (5.318,0.682) -- (5.682,0.318) .. controls (6,0) .. (6.318,-0.318) -- (8.682,-2.682);
  \draw[rounded corners=4pt] (11.682,1.682) -- (10.318,0.318) .. controls (10,0) .. (9.682,-0.318) -- (9.318,-0.682);
  \draw[rounded corners=4pt] (11.799,1.598) -- (11.201,0.402) .. controls (11,0) .. (10.799,-0.402) -- (10.201,-1.598);
  \draw[rounded corners=4pt] (12.318,1.682) -- (12.682,1.318);
  \draw[rounded corners=4pt] (12.682,0.682) -- (12.318,0.318) .. controls (12,0) .. (11.682,-0.318) -- (10.318,-1.682);
  \draw[rounded corners=4pt] (13,0.55) -- (13,0.45) .. controls (13,0) .. (12.777,-0.391) -- (9.223,-6.609);
  \draw[rounded corners=4pt] (13.318,0.682) -- (13.682,0.318) .. controls (14,0) .. (14.318,-0.318) -- (14.682,-0.682);
  \draw[rounded corners=4pt] (8.318,-7.682) -- (8.682,-7.318);
  \draw[rounded corners=4pt] (8.682,-6.682) -- (7.318,-5.318);
  \draw[rounded corners=4pt] (9.318,-6.682) -- (14.682,-1.318);
  \draw[rounded corners=4pt] (6.682,-4.682) -- (3.318,-1.318);
  \draw[rounded corners=4pt] (7.318,-4.682) -- (8.682,-3.318);
  \draw[rounded corners=4pt] (9.318,-2.682) -- (9.682,-2.318);
  \draw[rounded corners=4pt] (9.682,-1.682) -- (9.318,-1.318);
  \draw[rounded corners=4pt] (8,8.45) -- (8,9) -- (-1,9) -- (-1,-9) -- (8,-9) -- (8,-8.45);
  \draw (8,8.45) .. controls (8,8) .. (8.296,7.661);
  \draw (7.318,7.318) .. controls (7,7) .. (7.124,6.567);
  \draw (4.318,4.318) .. controls (4,4) .. (4.27,3.64);
  \draw (3.318,3.318) .. controls (3,3) .. (3,2.55);
  \draw (1.318,1.318) .. controls (1,1) .. (1,0.55);
  \draw (4.682,1.318) .. controls (5,1) .. (5,0.55);
  \draw (11.682,2.318) .. controls (12,2) .. (11.799,1.598);
  \draw (12.682,1.318) .. controls (13,1) .. (13,0.55);
  \draw (8,-8.45) .. controls (8,-8) .. (7.704,-7.661);
  \draw (8.682,-7.318) .. controls (9,-7) .. (9.223,-6.609);
  \draw (7.318,-5.318) .. controls (7,-5) .. (6.833,-4.582);
  \draw (3.318,-1.318) .. controls (3,-1) .. (3,-0.55);
  \draw (8.682,-3.318) .. controls (9,-3) .. (8.75,-2.626);
  \draw (9.682,-2.318) .. controls (10,-2) .. (10.201,-1.598);
  \draw (9.318,-1.318) .. controls (9,-1) .. (9,-0.55);
  \draw (14.682,-1.318) .. controls (15,-1) .. (15,-0.55);
  \draw[white, line width=6pt, line cap=butt] (7.682,7.682) .. controls (8,8) .. (8.318,7.682);
  \draw (7.682,7.682) .. controls (8,8) .. (8.318,7.682);
  \draw[white, line width=6pt, line cap=butt] (6.682,6.682) .. controls (7,7) .. (7.318,6.682);
  \draw (6.682,6.682) .. controls (7,7) .. (7.318,6.682);
  \draw[white, line width=6pt, line cap=butt] (3.682,3.682) .. controls (4,4) .. (4.318,3.682);
  \draw (3.682,3.682) .. controls (4,4) .. (4.318,3.682);
  \draw[white, line width=6pt, line cap=butt] (2.682,2.682) .. controls (3,3) .. (3.318,2.682);
  \draw (2.682,2.682) .. controls (3,3) .. (3.318,2.682);
  \draw[white, line width=6pt, line cap=butt] (0.682,0.682) .. controls (1,1) .. (1.318,0.682);
  \draw (0.682,0.682) .. controls (1,1) .. (1.318,0.682);
  \draw[white, line width=6pt, line cap=butt] (4.682,0.682) .. controls (5,1) .. (5.318,0.682);
  \draw (4.682,0.682) .. controls (5,1) .. (5.318,0.682);
  \draw[white, line width=6pt, line cap=butt] (11.682,1.682) .. controls (12,2) .. (12.318,1.682);
  \draw (11.682,1.682) .. controls (12,2) .. (12.318,1.682);
  \draw[white, line width=6pt, line cap=butt] (12.682,0.682) .. controls (13,1) .. (13.318,0.682);
  \draw (12.682,0.682) .. controls (13,1) .. (13.318,0.682);
  \draw[white, line width=6pt, line cap=butt] (7.682,-7.682) .. controls (8,-8) .. (8.318,-7.682);
  \draw (7.682,-7.682) .. controls (8,-8) .. (8.318,-7.682);
  \draw[white, line width=6pt, line cap=butt] (8.682,-6.682) .. controls (9,-7) .. (9.318,-6.682);
  \draw (8.682,-6.682) .. controls (9,-7) .. (9.318,-6.682);
  \draw[white, line width=6pt, line cap=butt] (6.682,-4.682) .. controls (7,-5) .. (7.318,-4.682);
  \draw (6.682,-4.682) .. controls (7,-5) .. (7.318,-4.682);
  \draw[white, line width=6pt, line cap=butt] (2.682,-0.682) .. controls (3,-1) .. (3.318,-0.682);
  \draw (2.682,-0.682) .. controls (3,-1) .. (3.318,-0.682);
  \draw[white, line width=6pt, line cap=butt] (8.682,-2.682) .. controls (9,-3) .. (9.318,-2.682);
  \draw (8.682,-2.682) .. controls (9,-3) .. (9.318,-2.682);
  \draw[white, line width=6pt, line cap=butt] (9.682,-1.682) .. controls (10,-2) .. (10.318,-1.682);
  \draw (9.682,-1.682) .. controls (10,-2) .. (10.318,-1.682);
  \draw[white, line width=6pt, line cap=butt] (8.682,-0.682) .. controls (9,-1) .. (9.318,-0.682);
  \draw (8.682,-0.682) .. controls (9,-1) .. (9.318,-0.682);
  \draw[white, line width=6pt, line cap=butt] (14.682,-0.682) .. controls (15,-1) .. (15.318,-0.682);
  \draw (14.682,-0.682) .. controls (15,-1) .. (15.318,-0.682);
\end{tikzpicture}} \\[2mm]
\resizebox{0.32\linewidth}{!}{\begin{tikzpicture}[scale=0.32, line cap=round]
  \draw (5,5.5) -- (0,0.5);
  \draw (5,5.5) -- (1,0.5);
  \draw (6,4.5) -- (2,0.5);
  \draw (6,4.5) -- (3,0.5);
  \draw (7,3.5) -- (4,0.5);
  \draw (7,3.5) -- (5,0.5);
  \draw (8,2.5) -- (6,0.5);
  \draw (8,1.5) -- (7,0.5);
  \draw (8,1.5) -- (8,0.5);
  \draw (8,1.5) -- (9,0.5);
  \draw (8,2.5) -- (8,1.5);
  \draw (8,2.5) -- (10,0.5);
  \draw (7,3.5) -- (8,2.5);
  \draw (6,4.5) -- (7,3.5);
  \draw (5,5.5) -- (6,4.5);
  \draw (5,-5.5) -- (0,-0.5);
  \draw (4,-3.5) -- (1,-0.5);
  \draw (4,-3.5) -- (2,-0.5);
  \draw (5,-2.5) -- (3,-0.5);
  \draw (5,-2.5) -- (4,-0.5);
  \draw (6,-1.5) -- (5,-0.5);
  \draw (6,-1.5) -- (6,-0.5);
  \draw (6,-1.5) -- (7,-0.5);
  \draw (5,-2.5) -- (6,-1.5);
  \draw (4,-3.5) -- (5,-2.5);
  \draw (5,-4.5) -- (4,-3.5);
  \draw (5,-4.5) -- (8,-0.5);
  \draw (5,-4.5) -- (9,-0.5);
  \draw (5,-5.5) -- (5,-4.5);
  \draw (5,-5.5) -- (10,-0.5);
  \fill (0,0.5) circle (2.5pt) (0,-0.5) circle (2.5pt);
  \fill (1,0.5) circle (2.5pt) (1,-0.5) circle (2.5pt);
  \fill (2,0.5) circle (2.5pt) (2,-0.5) circle (2.5pt);
  \fill (3,0.5) circle (2.5pt) (3,-0.5) circle (2.5pt);
  \fill (4,0.5) circle (2.5pt) (4,-0.5) circle (2.5pt);
  \fill (5,0.5) circle (2.5pt) (5,-0.5) circle (2.5pt);
  \fill (6,0.5) circle (2.5pt) (6,-0.5) circle (2.5pt);
  \fill (7,0.5) circle (2.5pt) (7,-0.5) circle (2.5pt);
  \fill (8,0.5) circle (2.5pt) (8,-0.5) circle (2.5pt);
  \fill (9,0.5) circle (2.5pt) (9,-0.5) circle (2.5pt);
  \fill (10,0.5) circle (2.5pt) (10,-0.5) circle (2.5pt);
\end{tikzpicture}} &
\resizebox{0.45\linewidth}{!}{\begin{tikzpicture}[scale=0.32, line cap=round, line join=round, >=stealth, thick]
  \draw[rounded corners=4pt] (4.682,4.682) -- (0.318,0.318) .. controls (0,0) .. (0.318,-0.318) -- (4.682,-4.682);
  \draw[rounded corners=4pt] (4.719,4.649) -- (1.281,0.351) .. controls (1,0) .. (1.318,-0.318) -- (3.682,-2.682);
  \draw[rounded corners=4pt] (5.318,4.682) -- (5.682,4.318);
  \draw[rounded corners=4pt] (5.682,3.682) -- (2.318,0.318) .. controls (2,0) .. (2.25,-0.374) -- (3.75,-2.626);
  \draw[rounded corners=4pt] (5.73,3.64) -- (3.27,0.36) .. controls (3,0) .. (3.318,-0.318) -- (4.682,-1.682);
  \draw[rounded corners=4pt] (6.318,3.682) -- (6.682,3.318);
  \draw[rounded corners=4pt] (6.682,2.682) -- (4.318,0.318) .. controls (4,0) .. (4.201,-0.402) -- (4.799,-1.598);
  \draw[rounded corners=4pt] (6.75,2.626) -- (5.25,0.374) .. controls (5,0) .. (5.318,-0.318) -- (5.682,-0.682);
  \draw[rounded corners=4pt] (7.318,2.682) -- (7.682,2.318);
  \draw[rounded corners=4pt] (7.682,1.682) -- (6.318,0.318) .. controls (6,0) .. (6,-0.45) -- (6,-0.55);
  \draw[rounded corners=4pt] (8,1.55) -- (8,1.45);
  \draw[rounded corners=4pt] (8.318,1.682) -- (9.682,0.318) .. controls (10,0) .. (9.682,-0.318) -- (5.318,-4.682);
  \draw[rounded corners=4pt] (7.682,0.682) -- (7.318,0.318) .. controls (7,0) .. (6.682,-0.318) -- (6.318,-0.682);
  \draw[rounded corners=4pt] (8,0.55) -- (8,0.45) .. controls (8,0) .. (7.73,-0.36) -- (5.27,-3.64);
  \draw[rounded corners=4pt] (8.318,0.682) -- (8.682,0.318) .. controls (9,0) .. (8.682,-0.318) -- (5.318,-3.682);
  \draw[rounded corners=4pt] (5,-4.55) -- (5,-4.45);
  \draw[rounded corners=4pt] (4.682,-3.682) -- (4.318,-3.318);
  \draw[rounded corners=4pt] (4.318,-2.682) -- (4.682,-2.318);
  \draw[rounded corners=4pt] (5.318,-1.682) -- (5.682,-1.318);
  \draw[rounded corners=4pt] (5,5.45) -- (5,6) -- (-1,6) -- (-1,-6) -- (5,-6) -- (5,-5.45);
  \draw (5,5.45) .. controls (5,5) .. (4.719,4.649);
  \draw (5.682,4.318) .. controls (6,4) .. (5.73,3.64);
  \draw (6.682,3.318) .. controls (7,3) .. (6.75,2.626);
  \draw (7.682,2.318) .. controls (8,2) .. (8,1.55);
  \draw (8,1.45) .. controls (8,1) .. (8,0.55);
  \draw (5,-5.45) .. controls (5,-5) .. (5,-4.55);
  \draw (5,-4.45) .. controls (5,-4) .. (5.27,-3.64);
  \draw (4.318,-3.318) .. controls (4,-3) .. (3.75,-2.626);
  \draw (4.682,-2.318) .. controls (5,-2) .. (4.799,-1.598);
  \draw (5.682,-1.318) .. controls (6,-1) .. (6,-0.55);
  \draw[white, line width=6pt, line cap=butt] (4.682,4.682) .. controls (5,5) .. (5.318,4.682);
  \draw (4.682,4.682) .. controls (5,5) .. (5.318,4.682);
  \draw[white, line width=6pt, line cap=butt] (5.682,3.682) .. controls (6,4) .. (6.318,3.682);
  \draw (5.682,3.682) .. controls (6,4) .. (6.318,3.682);
  \draw[white, line width=6pt, line cap=butt] (6.682,2.682) .. controls (7,3) .. (7.318,2.682);
  \draw (6.682,2.682) .. controls (7,3) .. (7.318,2.682);
  \draw[white, line width=6pt, line cap=butt] (7.682,1.682) .. controls (8,2) .. (8.318,1.682);
  \draw (7.682,1.682) .. controls (8,2) .. (8.318,1.682);
  \draw[white, line width=6pt, line cap=butt] (7.682,0.682) .. controls (8,1) .. (8.318,0.682);
  \draw (7.682,0.682) .. controls (8,1) .. (8.318,0.682);
  \draw[white, line width=6pt, line cap=butt] (4.682,-4.682) .. controls (5,-5) .. (5.318,-4.682);
  \draw (4.682,-4.682) .. controls (5,-5) .. (5.318,-4.682);
  \draw[white, line width=6pt, line cap=butt] (4.682,-3.682) .. controls (5,-4) .. (5.318,-3.682);
  \draw (4.682,-3.682) .. controls (5,-4) .. (5.318,-3.682);
  \draw[white, line width=6pt, line cap=butt] (3.682,-2.682) .. controls (4,-3) .. (4.318,-2.682);
  \draw (3.682,-2.682) .. controls (4,-3) .. (4.318,-2.682);
  \draw[white, line width=6pt, line cap=butt] (4.682,-1.682) .. controls (5,-2) .. (5.318,-1.682);
  \draw (4.682,-1.682) .. controls (5,-2) .. (5.318,-1.682);
  \draw[white, line width=6pt, line cap=butt] (5.682,-0.682) .. controls (6,-1) .. (6.318,-0.682);
  \draw (5.682,-0.682) .. controls (6,-1) .. (6.318,-0.682);
\end{tikzpicture}}
\end{tabular}
\caption{The Borromean rings. Top: a binary tree diagram with 9 leaves (Sweeney). Bottom: a ternary tree
diagram with $n=6$.}
\label{fig:borromean}
\end{figure}

\begin{theorem}\label{thm:borromean}
The Borromean rings have $F$-index $9$ and $F_3$-index $6$ (in the normalization of~\cite{Jones2019}, $8$ and $5$).
\end{theorem}
\begin{proof}
For each $n$ all the tree diagrams of size $n$ are enumerated, without any reduction. A diagram can give the
Borromean rings only if:
\begin{enumerate}
\item its link has $3$ components;
\item all the linking numbers vanish;
\item its HOMFLY-PT polynomial is that of the Borromean rings up to mirror image (with vanishing linking
  numbers it does not depend on the orientations of the components).
\end{enumerate}
No diagram of $F$ with $n\le 8$ satisfies these three conditions (with $n=8$, 52\,918 of the 184\,041 diagrams
have $3$ components and 41\,872 of them have vanishing linking numbers, but none has the right polynomial),
and none of $F_3$ with $n\le 5$. The diagrams of Figure~\ref{fig:borromean} give the Borromean rings: the
verified canonical retriangulations of their complements and of that of $L6a4$ (Section~\ref{sec:certified-tables})
have a combinatorial isomorphism taking meridians to meridians on every cusp, which gives a homeomorphism of the
complements extending to one of $S^3$ that takes one link to the other. The same certified test shows that all the
diagrams of these sizes satisfying the three conditions above give the Borromean rings.
\end{proof}
In $F$ the Borromean rings are realized by 64 diagrams with 9 leaves, all found by
Sweeney~\cite{Sweeney2021}, who answered Jones' question for $F$; in $F_3$ by 14 diagrams with $n=6$.

\subsection{Connected sums}\label{sec:sums}
By Theorems~\ref{thm:unoriented}\ref{it:sum} and~\ref{thm:oriented}\ref{it:sum-or}, grafting minimal diagrams
gives
\[
\ind(K_1\#K_2)\le\ind(K_1)+\ind(K_2)-1
\]
for all four indices. Table~\ref{tab:sums} shows that this inequality is often strict.

\begin{table}[H]
\centering
\begin{tabular}{lcccc}
\toprule
knot & $F$ & $F_3$ & $\vF$ & $\vFt$ \\
\midrule
$3_1\#3_1$ (granny knot) & $\mathbf{8}$ (9) & $\mathbf{5}$ (7) & $17$ (17) & $\mathbf{5}$ (9) \\
$3_1\#3_1^m$ (square knot) & $9$ (9) & $\mathbf{6}$ (7) & $\mathbf{13}$ (17) & $\mathbf{7}$ (9) \\
$3_1\#4_1$ & $\mathbf{9}$ (11) & $\mathbf{6}$ (8) & $\mathbf{13}$ (15) & $\mathbf{7}$ (9) \\
\bottomrule
\end{tabular}
\caption{The four indices of three connected sums; in brackets, $\ind(K_1)+\ind(K_2)-1$. Bold: the index is
smaller.}
\label{tab:sums}
\end{table}

The values are computed as follows (script \texttt{connected\_sums.py}). The enumeration skips the non-prime tree
diagrams, which are grafts, but this loses nothing below the bound: by the uniqueness of the prime decomposition, a
graft realizing $K_1\#K_2$ either has two pieces realizing $K_1$ and $K_2$, and then its size is at least
$\ind(K_1)+\ind(K_2)-1$, or has a trivial piece, and then the other piece is a smaller diagram of $K_1\#K_2$. Hence
the first $n$ at which the HOMFLY-PT class of $K_1\#K_2$ (the product of the polynomials of the factors) occurs is a
lower bound. At that $n$ the diagrams of the class are connected sums of two diagrams with $3$ crossings and
determinant $3$ (trefoils) or with $4$ crossings and determinant $5$ (figure-eight knots), the signatures of the
two trefoils telling apart the granny knot and the square knot. All the factors are invertible, so the orientations
play no role. Remarkably, the granny knot has $F_3$-index $5$, only one more than the trefoil, and the same
$\vFt$-index as the trefoil, $5$: a single diagram with $n=5$ realizes it in both groups.

\subsection{Comparison with classical invariants}\label{sec:bounds}
Table~\ref{tab:bounds} compares the indices with the braid index $b$ (KnotInfo, knots with at most 12
crossings) and with the crossing number $c$.

\begin{table}[H]
\centering
{\footnotesize\begin{tabular}{lccp{3.6cm}ccc}
\toprule
index & knots & $\min(\mathrm{ind}-b)$ & attained by (ind/$b$) & $\max(\mathrm{ind}-c)$ & $\min(2(\mathrm{ind}-1)-c)$ & $c=2(\mathrm{ind}-1)$ \\
\midrule
$F$ & 2149 & 3 & $3_{1}$ (5/2) & 3 & 4 & 0 \\
$F_3$ & 2977 & 1 & $6_{1}$ (5/4), $10_{13}$ (7/6), $10_{35}$ (7/6), \dots\ (6) & 1 & 0 & 105 \\
$\vec F$ & 1047 & 4 & $4_{1}$ (7/3) & 14 & 8 & 0 \\
$\vec F_3$ & 2977 & 2 & $4_{1}$ (5/3), $5_{2}$ (5/3), $8_{12}$ (7/5), \dots\ (9) & 2 & 0 & 3 \\
\bottomrule
\end{tabular}
}
\caption{Indices versus braid index $b$ and crossing number $c$, over the knots with known index. The last
column counts the knots attaining $c=2(\ind-1)$.}
\label{tab:bounds}
\end{table}

The braid index is strictly smaller than each index for every non-trivial knot, although
Theorem~\ref{thm:unoriented}\ref{it:grid} only covers $F$ and $\vF$; the Golan--Sapir bound $\ind_F\le 12c$ is far from
sharp ($\ind_F\le c+3$ on the data). Since a tree diagram of size $n$ gives a link diagram with
exactly $2(n-1)$ crossings, the equality $c=2(\ind-1)$ holds if and only if some tree diagram gives a minimal
crossing diagram of $K$ (such a tree diagram has then size $\ind(K)$). In our data this happens for $F_3$ and
$\vFt$ only:

\noindent\textbf{$F_3$} (105 knots): $8_{19}$, $8_{20}$, $10_{124}$, $10_{128}$, $10_{132}$, $10_{133}$, $10_{134}$, $10_{136}$, $10_{137}$, $10_{138}$, $10_{140}$, $10_{141}$, $10_{147}$, $10_{150}$, $10_{165}$, $12n_{1}$, $12n_{2}$, $12n_{3}$, $12n_{6}$, $12n_{10}$, $12n_{20}$, $12n_{38}$, $12n_{39}$, $12n_{40}$, $12n_{45}$, $12n_{68}$, $12n_{71}$, $12n_{73}$, $12n_{80}$, $12n_{84}$, $12n_{97}$, $12n_{119}$, $12n_{120}$, $12n_{121}$, $12n_{130}$, $12n_{142}$, $12n_{146}$, $12n_{150}$, $12n_{152}$, $12n_{153}$, $12n_{154}$, $12n_{156}$, $12n_{157}$, $12n_{159}$, $12n_{160}$, $12n_{162}$, $12n_{172}$, $12n_{176}$, $12n_{182}$, $12n_{188}$, $12n_{196}$, $12n_{198}$, $12n_{218}$, $12n_{271}$, $12n_{275}$, $12n_{286}$, $12n_{288}$, $12n_{293}$, $12n_{300}$, $12n_{301}$, $12n_{310}$, $12n_{312}$, $12n_{323}$, $12n_{326}$, $12n_{327}$, $12n_{333}$, $12n_{338}$, $12n_{342}$, $12n_{358}$, $12n_{362}$, $12n_{368}$, $12n_{370}$, $12n_{379}$, $12n_{381}$, $12n_{388}$, $12n_{399}$, $12n_{407}$, $12n_{408}$, $12n_{420}$, $12n_{427}$, $12n_{428}$, $12n_{438}$, $12n_{443}$, $12n_{444}$, $12n_{449}$, $12n_{451}$, $12n_{452}$, $12n_{459}$, $12n_{462}$, $12n_{485}$, $12n_{491}$, $12n_{502}$, $12n_{512}$, $12n_{568}$, $12n_{590}$, $12n_{591}$, $12n_{592}$, $12n_{597}$, $12n_{603}$, $12n_{625}$, $12n_{634}$, $12n_{705}$, $12n_{706}$, $12n_{742}$, $12n_{838}$.

\smallskip\noindent\textbf{$\vec F_3$} (3 knots): $12n_{3}$, $12n_{20}$, $12n_{119}$.

\medskip\noindent The genus bound of Theorem~\ref{thm:oriented}\ref{it:genus}:

\medskip
\noindent\textbf{$\vec F$}: 2191 knots with known index and genus; $\min\big(\tfrac{\mathrm{ind}-1}{2}-g\big) = 2$, attained by 3 knots (with at most 10 crossings: $4_{1}$, $8_{18}$, $10_{123}$).

\smallskip
\noindent\textbf{$\vec F_3$}: 12711 knots with known index and genus; $\min\big(\tfrac{\mathrm{ind}-1}{2}-g\big) = 0$, attained by 2571 knots (with at most 10 crossings: $5_{1}$, $7_{1}$, $8_{2}$, $8_{5}$, $8_{7}$, $8_{9}$, $8_{10}$, $8_{17}$, $8_{19}$, $9_{1}$, $9_{6}$, $9_{9}$, $9_{11}$, $9_{16}$, $9_{20}$, $9_{24}$, $9_{28}$, $9_{43}$, $10_{2}$, $10_{5}$, $10_{9}$, $10_{17}$, $10_{46}$, $10_{47}$, $10_{48}$, $10_{49}$, $10_{62}$, $10_{64}$, $10_{66}$, $10_{79}$, $10_{80}$, $10_{82}$, $10_{85}$, $10_{91}$, $10_{94}$, $10_{99}$, $10_{100}$, $10_{104}$, $10_{106}$, $10_{109}$, $10_{118}$, $10_{124}$, $10_{125}$, $10_{128}$, $10_{134}$, $10_{139}$, $10_{141}$, $10_{149}$, $10_{150}$, $10_{151}$, $10_{152}$, $10_{153}$, $10_{154}$, $10_{156}$, $10_{161}$).

\medskip\noindent\emph{Oriented and unoriented indices.} By Theorem~\ref{thm:oriented}\ref{it:inequalities},
$\ind_F\le\ind_{\vF}$ and $\ind_{F_3}\le\ind_{\vFt}$. Since the oriented indices of a knot are odd, equality
requires an odd unoriented index; for instance the trefoil has $F_3$-index $4$ and $\vFt$-index $5$, while its
$F$-index is $5$ and its $\vF$-index $9$.
In $\vF$ equality is rare: among the 1\,604 knots for which both indices are known it holds only for 7 knots,
among them the figure-eight knot (index $7$), $8_{18}$ ($11$) and $10_{45}$ ($13$), always in both orientations;
usually $\ind_{\vF}$ exceeds $\ind_F$ by $2$ to $5$. In $\vFt$ it is frequent: it holds for 1\,236 of the
12\,141 knots for which both indices are known (851 of them in both orientations), starting with $4_1$, $5_1$,
$5_2$ (index $5$) and eight knots with 8 crossings (index $7$); for $8_{17}$ it holds only for the orientation
$K$, since $\ind_{F_3}(8_{17})=\ind_{\vFt}(K)=7$ and $\ind_{\vFt}(K^r)=9$. The difference
$\ind_{\vFt}-\ind_{F_3}$ never exceeds $4$.

\section{Open problems}\label{sec:problems}

Problem~\ref{prob:invariants} remains open: the oriented indices detect non-invertibility in many cases, but
not always (997 non-invertible knots have equal $\vF$-index in both orientations, among them $8_{17}$,
$9_{32}$ and $9_{33}$). Our data suggest the following questions.

\begin{question}
Determine the $\vF$-index of the 50 prime knots with at most 10 crossings for which it is still unknown (13 with
9 crossings and 37 with 10, all with $\vF$-index at least $21$; Section~\ref{sec:targeted}). Is it $23$ for
$9_7$, $9_{18}$, $9_{23}$, $10_6$ and $10_{14}$, where the descent stops?
\end{question}

\begin{question}
Which refinements of the oriented indices (for instance the number of minimal diagrams realizing $K$ and
$K^r$, or invariants of these diagrams) detect non-invertibility when $\ind(K)=\ind(K^r)$?
\end{question}
\begin{remark}
By Theorem~\ref{thm:oriented}\ref{it:inverse}, the map $g\mapsto g^{-1}$ is a size-preserving bijection between the
diagrams representing $L$ and those representing $L^{mr}$, both in $\vF$ and in $\vFt$. Applied to $L=K^m$, it
gives a size-preserving bijection between the diagrams of $K^m$ and those of $K^r$. If $K$ is positive amphicheiral,
that is $K=K^m$, then $K$ and $K^r$ have the same index and, for every size, the same number of diagrams; in
particular they have the same number of minimal diagrams. Hence, for the non-invertible positive amphicheiral knots,
a refinement as in Question~1 must use features of the diagrams which are not preserved by exchanging the two trees.
No such obstruction exists for chiral and negative amphicheiral knots.
\end{remark}

\begin{question}
Is $\operatorname{braid}(K)<\ind_{F_3}(K)$ for every non-trivial knot $K$? (On the data,
$\ind_{F_3}(K)-\operatorname{braid}(K)\ge 1$, with equality for six knots.)
\end{question}

\begin{question}
Characterize the knots for which the bounds of Theorem~\ref{thm:unoriented}\ref{it:crossings} and
Theorem~\ref{thm:oriented}\ref{it:genus} are attained: the knots with $c(K)=2(\ind_{F_3}(K)-1)$, i.e.\ those
having a minimal crossing diagram coming from $F_3$ (105 knots in our data), and the knots with
$g(K)=(\ind_{\vFt}(K)-1)/2$, i.e.\ those for which the checkerboard surface of a minimal diagram of $\vFt$ is a
minimal genus Seifert surface (2\,571 knots). Neither bound is attained for $F$ and $\vF$ in our data.
\end{question}

\begin{question}
Which mutant pairs are distinguished by the Thompson indices? Is there a mutation that always preserves
them, for instance the mutations of pretzel knots?
\end{question}

\appendix
\section{Comparison with Sweeney's table}\label{app:sweeney}

D.~Sweeney enumerated the tree diagrams of $F$ and of $\vF$ with few leaves~\cite{Sweeney2021}.
Table~\ref{tab:sweeney-range} compares the ranges of his computation and of ours.
\begin{itemize}
\item \emph{$F$.} His enumeration is complete up to 10 leaves; with 11 leaves it covers the diagrams with one
  component, and only a sample of them is published. We used the same reductions as Sweeney (prime diagrams, one
  representative for the four reflections). He gives the $F$-index of the prime knots with at most 10
  crossings, as exact values up to 11 leaves and ``$\ge 12$'' otherwise.
\item \emph{$\vF$.} His files contain all the diagrams of $\vF$ with at most 11 leaves, grouped by HOMFLY-PT
  polynomial but without the names of the knots, and no values of the $\vF$-index. Up to 11 leaves they agree
  with our enumeration: the same 11 HOMFLY-PT classes of knots, each appearing first at the same $n$.
\end{itemize}

\begin{table}[H]
\centering
\begin{tabular}{lll}
\toprule
 & Sweeney~\cite{Sweeney2021} & this paper \\
\midrule
$F$: diagrams & all up to 10 leaves; knots with 11 leaves & all up to 13 leaves \\
$F$: knots & $\le 10$ crossings, values up to 11 & $\le 13$ crossings (12\,965 knots) \\
$\vF$: diagrams & all up to 11 leaves (unlabelled) & all up to 17 leaves (19: polynomials) \\
$\vF$: knots & -- & $\le 13$ crossings, both orientations \\
$F_3$, $\vFt$ & -- & $n\le 9$ and $n\le 11$ \\
\bottomrule
\end{tabular}
\caption{Ranges of Sweeney's computation and of ours.}
\label{tab:sweeney-range}
\end{table}

He also determined the $F$-index of the Borromean rings, $9$, answering Jones' question for $F$; we confirm it
(Section~\ref{sec:borromean}). We determined its $F_3$-index.

\begin{figure}[ht]
\centering
\begin{tabular}{cc}
\resizebox{0.28\linewidth}{!}{\begin{tikzpicture}[scale=0.32, line cap=round]
  \draw (2.5,1) -- (2,0.5);
  \draw (2.5,1) -- (3,0.5);
  \draw (3,1.5) -- (2.5,1);
  \draw (3,1.5) -- (4,0.5);
  \draw (2.5,2) -- (1,0.5);
  \draw (2.5,2) -- (3,1.5);
  \draw (3,2.5) -- (2.5,2);
  \draw (3,2.5) -- (5,0.5);
  \draw (7.5,1) -- (7,0.5);
  \draw (7.5,1) -- (8,0.5);
  \draw (7,1.5) -- (6,0.5);
  \draw (7,1.5) -- (7.5,1);
  \draw (7.5,2) -- (7,1.5);
  \draw (7.5,2) -- (9,0.5);
  \draw (5,4.5) -- (3,2.5);
  \draw (5,4.5) -- (7.5,2);
  \draw (4.5,5) -- (0,0.5);
  \draw (4.5,5) -- (5,4.5);
  \draw (0.5,-1) -- (0,-0.5);
  \draw (0.5,-1) -- (1,-0.5);
  \draw (1,-1.5) -- (0.5,-1);
  \draw (1,-1.5) -- (2,-0.5);
  \draw (5.5,-1) -- (5,-0.5);
  \draw (5.5,-1) -- (6,-0.5);
  \draw (5,-1.5) -- (4,-0.5);
  \draw (5,-1.5) -- (5.5,-1);
  \draw (5.5,-2) -- (5,-1.5);
  \draw (5.5,-2) -- (7,-0.5);
  \draw (5,-2.5) -- (3,-0.5);
  \draw (5,-2.5) -- (5.5,-2);
  \draw (5.5,-3) -- (5,-2.5);
  \draw (5.5,-3) -- (8,-0.5);
  \draw (6,-3.5) -- (5.5,-3);
  \draw (6,-3.5) -- (9,-0.5);
  \draw (4.5,-5) -- (1,-1.5);
  \draw (4.5,-5) -- (6,-3.5);
  \fill (0,0.5) circle (2.5pt) (0,-0.5) circle (2.5pt);
  \fill (1,0.5) circle (2.5pt) (1,-0.5) circle (2.5pt);
  \fill (2,0.5) circle (2.5pt) (2,-0.5) circle (2.5pt);
  \fill (3,0.5) circle (2.5pt) (3,-0.5) circle (2.5pt);
  \fill (4,0.5) circle (2.5pt) (4,-0.5) circle (2.5pt);
  \fill (5,0.5) circle (2.5pt) (5,-0.5) circle (2.5pt);
  \fill (6,0.5) circle (2.5pt) (6,-0.5) circle (2.5pt);
  \fill (7,0.5) circle (2.5pt) (7,-0.5) circle (2.5pt);
  \fill (8,0.5) circle (2.5pt) (8,-0.5) circle (2.5pt);
  \fill (9,0.5) circle (2.5pt) (9,-0.5) circle (2.5pt);
\end{tikzpicture}} &
\resizebox{0.35\linewidth}{!}{\begin{tikzpicture}[scale=0.32, line cap=round, line join=round, >=stealth, thick]
  \draw[rounded corners=4pt] (8.682,8.682) -- (0.318,0.318) .. controls (0,0) .. (0.318,-0.318) -- (0.682,-0.682);
  \draw[rounded corners=4pt] (8.701,8.664) -- (1.299,0.336) .. controls (1,0) .. (1,-0.45) -- (1,-0.55);
  \draw[rounded corners=4pt] (9.318,8.682) -- (9.682,8.318);
  \draw[rounded corners=4pt] (9.682,7.682) -- (6.318,4.318);
  \draw[rounded corners=4pt] (10.056,7.553) -- (10.944,0.447) .. controls (11,0) .. (11,-0.45) -- (11,-0.55);
  \draw[rounded corners=4pt] (10.318,7.682) -- (14.682,3.318);
  \draw[rounded corners=4pt] (5.682,3.682) -- (5.318,3.318);
  \draw[rounded corners=4pt] (6.27,3.64) -- (8.73,0.36) .. controls (9,0) .. (9.201,-0.402) -- (9.799,-1.598);
  \draw[rounded corners=4pt] (6.318,3.682) -- (9.682,0.318) .. controls (10,0) .. (10.318,-0.318) -- (10.682,-0.682);
  \draw[rounded corners=4pt] (4.682,2.682) -- (2.318,0.318) .. controls (2,0) .. (1.682,-0.318) -- (1.318,-0.682);
  \draw[rounded corners=4pt] (4.75,2.626) -- (3.25,0.374) .. controls (3,0) .. (2.799,-0.402) -- (2.201,-1.598);
  \draw[rounded corners=4pt] (5.318,2.682) -- (5.682,2.318);
  \draw[rounded corners=4pt] (5.682,1.682) -- (5.318,1.318);
  \draw[rounded corners=4pt] (6.201,1.598) -- (6.799,0.402) .. controls (7,0) .. (7.27,-0.36) -- (9.73,-3.64);
  \draw[rounded corners=4pt] (6.318,1.682) -- (7.682,0.318) .. controls (8,0) .. (8.318,-0.318) -- (9.682,-1.682);
  \draw[rounded corners=4pt] (4.682,0.682) -- (4.318,0.318) .. controls (4,0) .. (3.682,-0.318) -- (2.318,-1.682);
  \draw[rounded corners=4pt] (5,0.55) -- (5,0.45) .. controls (5,0) .. (5.183,-0.411) -- (8.817,-8.589);
  \draw[rounded corners=4pt] (5.318,0.682) -- (5.682,0.318) .. controls (6,0) .. (6.318,-0.318) -- (9.682,-3.682);
  \draw[rounded corners=4pt] (14.682,2.682) -- (14.318,2.318);
  \draw[rounded corners=4pt] (15.25,2.626) -- (16.75,0.374) .. controls (17,0) .. (16.712,-0.346) -- (12.288,-5.654);
  \draw[rounded corners=4pt] (15.318,2.682) -- (17.682,0.318) .. controls (18,0) .. (17.682,-0.318) -- (12.318,-5.682);
  \draw[rounded corners=4pt] (13.682,1.682) -- (12.318,0.318) .. controls (12,0) .. (11.682,-0.318) -- (11.318,-0.682);
  \draw[rounded corners=4pt] (13.799,1.598) -- (13.201,0.402) .. controls (13,0) .. (12.75,-0.374) -- (11.25,-2.626);
  \draw[rounded corners=4pt] (14.318,1.682) -- (14.682,1.318);
  \draw[rounded corners=4pt] (14.682,0.682) -- (14.318,0.318) .. controls (14,0) .. (13.682,-0.318) -- (11.318,-2.682);
  \draw[rounded corners=4pt] (15,0.55) -- (15,0.45) .. controls (15,0) .. (14.719,-0.351) -- (11.281,-4.649);
  \draw[rounded corners=4pt] (15.318,0.682) -- (15.682,0.318) .. controls (16,0) .. (15.682,-0.318) -- (11.318,-4.682);
  \draw[rounded corners=4pt] (8.682,-8.682) -- (2.318,-2.318);
  \draw[rounded corners=4pt] (9.318,-8.682) -- (11.682,-6.318);
  \draw[rounded corners=4pt] (1.682,-1.682) -- (1.318,-1.318);
  \draw[rounded corners=4pt] (11.682,-5.682) -- (11.318,-5.318);
  \draw[rounded corners=4pt] (10.682,-4.682) -- (10.318,-4.318);
  \draw[rounded corners=4pt] (10.318,-3.682) -- (10.682,-3.318);
  \draw[rounded corners=4pt] (10.682,-2.682) -- (10.318,-2.318);
  \draw[rounded corners=4pt] (10.318,-1.682) -- (10.682,-1.318);
  \draw[rounded corners=4pt] (9,9.45) -- (9,10) -- (-1,10) -- (-1,-10) -- (9,-10) -- (9,-9.45);
  \draw (9,9.45) .. controls (9,9) .. (8.701,8.664);
  \draw (9.682,8.318) .. controls (10,8) .. (10.056,7.553);
  \draw (6.318,4.318) .. controls (6,4) .. (6.27,3.64);
  \draw (5.318,3.318) .. controls (5,3) .. (4.75,2.626);
  \draw (5.682,2.318) .. controls (6,2) .. (6.201,1.598);
  \draw (5.318,1.318) .. controls (5,1) .. (5,0.55);
  \draw (14.682,3.318) .. controls (15,3) .. (15.25,2.626);
  \draw (14.318,2.318) .. controls (14,2) .. (13.799,1.598);
  \draw (14.682,1.318) .. controls (15,1) .. (15,0.55);
  \draw (9,-9.45) .. controls (9,-9) .. (8.817,-8.589);
  \draw (2.318,-2.318) .. controls (2,-2) .. (2.201,-1.598);
  \draw (1.318,-1.318) .. controls (1,-1) .. (1,-0.55);
  \draw (11.682,-6.318) .. controls (12,-6) .. (12.288,-5.654);
  \draw (11.318,-5.318) .. controls (11,-5) .. (11.281,-4.649);
  \draw (10.318,-4.318) .. controls (10,-4) .. (9.73,-3.64);
  \draw (10.682,-3.318) .. controls (11,-3) .. (11.25,-2.626);
  \draw (10.318,-2.318) .. controls (10,-2) .. (9.799,-1.598);
  \draw (10.682,-1.318) .. controls (11,-1) .. (11,-0.55);
  \draw[white, line width=6pt, line cap=butt] (8.682,8.682) .. controls (9,9) .. (9.318,8.682);
  \draw (8.682,8.682) .. controls (9,9) .. (9.318,8.682);
  \draw[white, line width=6pt, line cap=butt] (9.682,7.682) .. controls (10,8) .. (10.318,7.682);
  \draw (9.682,7.682) .. controls (10,8) .. (10.318,7.682);
  \draw[white, line width=6pt, line cap=butt] (5.682,3.682) .. controls (6,4) .. (6.318,3.682);
  \draw (5.682,3.682) .. controls (6,4) .. (6.318,3.682);
  \draw[white, line width=6pt, line cap=butt] (4.682,2.682) .. controls (5,3) .. (5.318,2.682);
  \draw (4.682,2.682) .. controls (5,3) .. (5.318,2.682);
  \draw[white, line width=6pt, line cap=butt] (5.682,1.682) .. controls (6,2) .. (6.318,1.682);
  \draw (5.682,1.682) .. controls (6,2) .. (6.318,1.682);
  \draw[white, line width=6pt, line cap=butt] (4.682,0.682) .. controls (5,1) .. (5.318,0.682);
  \draw (4.682,0.682) .. controls (5,1) .. (5.318,0.682);
  \draw[white, line width=6pt, line cap=butt] (14.682,2.682) .. controls (15,3) .. (15.318,2.682);
  \draw (14.682,2.682) .. controls (15,3) .. (15.318,2.682);
  \draw[white, line width=6pt, line cap=butt] (13.682,1.682) .. controls (14,2) .. (14.318,1.682);
  \draw (13.682,1.682) .. controls (14,2) .. (14.318,1.682);
  \draw[white, line width=6pt, line cap=butt] (14.682,0.682) .. controls (15,1) .. (15.318,0.682);
  \draw (14.682,0.682) .. controls (15,1) .. (15.318,0.682);
  \draw[white, line width=6pt, line cap=butt] (8.682,-8.682) .. controls (9,-9) .. (9.318,-8.682);
  \draw (8.682,-8.682) .. controls (9,-9) .. (9.318,-8.682);
  \draw[white, line width=6pt, line cap=butt] (1.682,-1.682) .. controls (2,-2) .. (2.318,-1.682);
  \draw (1.682,-1.682) .. controls (2,-2) .. (2.318,-1.682);
  \draw[white, line width=6pt, line cap=butt] (0.682,-0.682) .. controls (1,-1) .. (1.318,-0.682);
  \draw (0.682,-0.682) .. controls (1,-1) .. (1.318,-0.682);
  \draw[white, line width=6pt, line cap=butt] (11.682,-5.682) .. controls (12,-6) .. (12.318,-5.682);
  \draw (11.682,-5.682) .. controls (12,-6) .. (12.318,-5.682);
  \draw[white, line width=6pt, line cap=butt] (10.682,-4.682) .. controls (11,-5) .. (11.318,-4.682);
  \draw (10.682,-4.682) .. controls (11,-5) .. (11.318,-4.682);
  \draw[white, line width=6pt, line cap=butt] (9.682,-3.682) .. controls (10,-4) .. (10.318,-3.682);
  \draw (9.682,-3.682) .. controls (10,-4) .. (10.318,-3.682);
  \draw[white, line width=6pt, line cap=butt] (10.682,-2.682) .. controls (11,-3) .. (11.318,-2.682);
  \draw (10.682,-2.682) .. controls (11,-3) .. (11.318,-2.682);
  \draw[white, line width=6pt, line cap=butt] (9.682,-1.682) .. controls (10,-2) .. (10.318,-1.682);
  \draw (9.682,-1.682) .. controls (10,-2) .. (10.318,-1.682);
  \draw[white, line width=6pt, line cap=butt] (10.682,-0.682) .. controls (11,-1) .. (11.318,-0.682);
  \draw (10.682,-0.682) .. controls (11,-1) .. (11.318,-0.682);
\end{tikzpicture}}
\end{tabular}
\caption{$10_{156}$ with 10 leaves.}
\label{fig:sw-10156}
\end{figure}
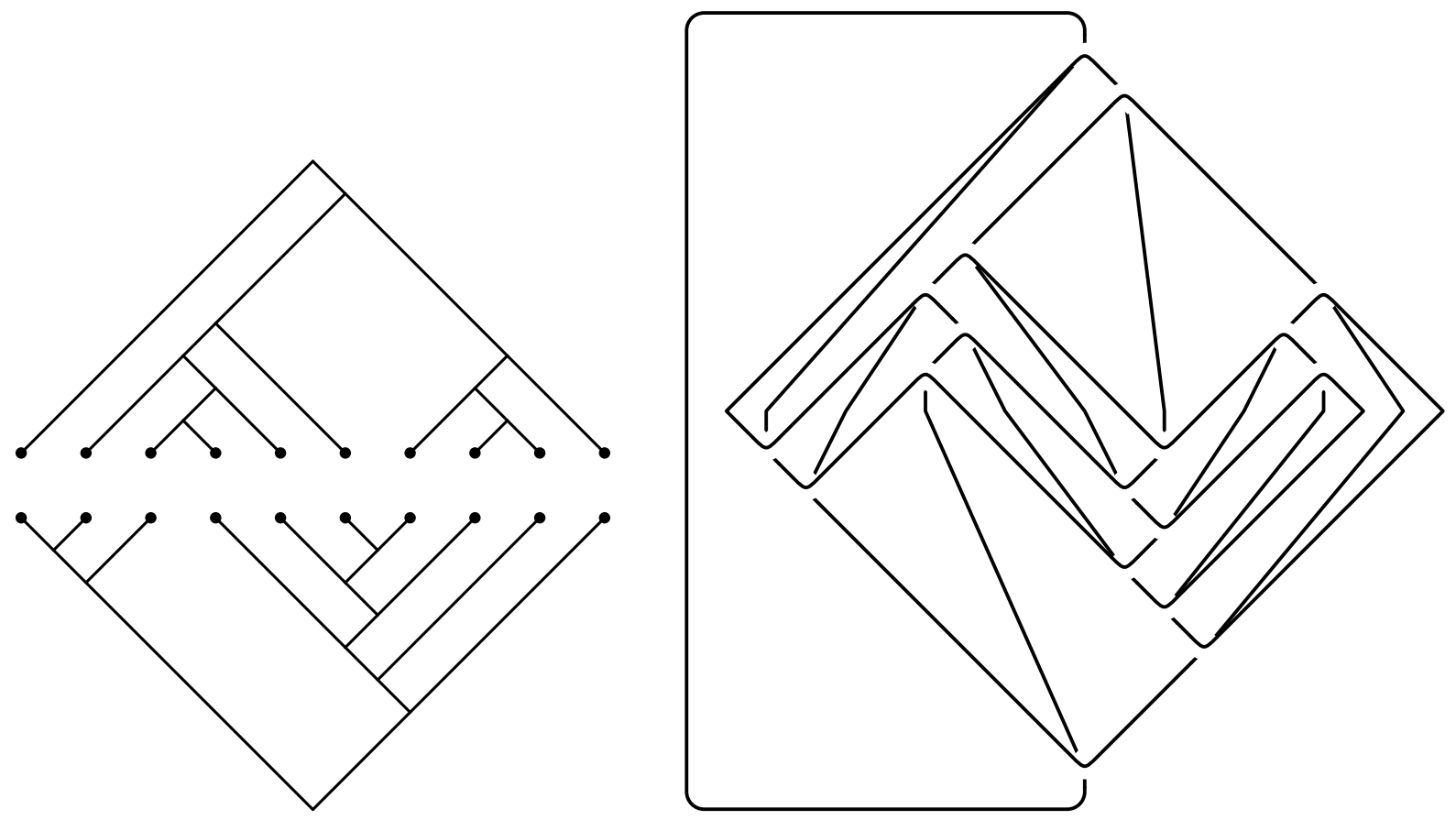

For the $F$-index, of the 249 prime knots with at most 10 crossings, 245 have the same value in his table and in
ours, including all the 87 entries ``$\ge 12$''. For two knots our value differs from his, and for two knots his
table gives no value: their HOMFLY-PT class is in his data, but the polynomial alone cannot separate them from
$5_1$ and $8_{16}$ (Table~\ref{tab:sweeney}). The corresponding diagrams are in
Figures~\ref{fig:sw-1045}--\ref{fig:sw-10156}.

\begin{table}[ht]
\centering
\begin{tabular}{lp{3.8cm}cp{7.6cm}}
\toprule
knot & Sweeney & $\ind_F$ & reason \\
\midrule
$10_{45}$ & 11 & \textbf{13} & his 11-leaf diagram is a different knot with the same HOMFLY-PT polynomial
(volume $16.8704$ versus $13.7161$); no diagram with 11 or 12 leaves gives $10_{45}$ \\
$10_{163}$ & 11 & \textbf{10} & a 10-leaf diagram exists (its complement is K10n35, i.e.\ $10_{163}$ in
KnotInfo); it is also in his data, in an unlabelled class \\
\midrule
$10_{132}$ & not listed (its class, shared with $5_1$, occurs at 8) & \textbf{9} & no 8-leaf diagram of the class is $10_{132}$; a 9-leaf
diagram is \\
$10_{156}$ & not listed (its class, shared with $8_{16}$, occurs at 10) & \textbf{10} & a 10-leaf diagram exists; the class does not occur with
fewer leaves \\
\bottomrule
\end{tabular}
\caption{Differences from Sweeney's table of $F$-indices. Top: values that differ. Bottom: values missing in
his table.}
\label{tab:sweeney}
\end{table}

\begin{figure}[ht]
\centering
\begin{tabular}{cc}
\resizebox{0.32\linewidth}{!}{\begin{tikzpicture}[scale=0.32, line cap=round]
  \draw (4.5,1) -- (4,0.5);
  \draw (4.5,1) -- (5,0.5);
  \draw (4,1.5) -- (3,0.5);
  \draw (4,1.5) -- (4.5,1);
  \draw (7.5,1) -- (7,0.5);
  \draw (7.5,1) -- (8,0.5);
  \draw (9.5,1) -- (9,0.5);
  \draw (9.5,1) -- (10,0.5);
  \draw (8.5,2) -- (7.5,1);
  \draw (8.5,2) -- (9.5,1);
  \draw (9,2.5) -- (8.5,2);
  \draw (9,2.5) -- (11,0.5);
  \draw (8.5,3) -- (6,0.5);
  \draw (8.5,3) -- (9,2.5);
  \draw (9,3.5) -- (8.5,3);
  \draw (9,3.5) -- (12,0.5);
  \draw (7.5,5) -- (4,1.5);
  \draw (7.5,5) -- (9,3.5);
  \draw (7,5.5) -- (2,0.5);
  \draw (7,5.5) -- (7.5,5);
  \draw (6.5,6) -- (1,0.5);
  \draw (6.5,6) -- (7,5.5);
  \draw (6,6.5) -- (0,0.5);
  \draw (6,6.5) -- (6.5,6);
  \draw (0.5,-1) -- (0,-0.5);
  \draw (0.5,-1) -- (1,-0.5);
  \draw (2.5,-1) -- (2,-0.5);
  \draw (2.5,-1) -- (3,-0.5);
  \draw (3,-1.5) -- (2.5,-1);
  \draw (3,-1.5) -- (4,-0.5);
  \draw (2,-2.5) -- (0.5,-1);
  \draw (2,-2.5) -- (3,-1.5);
  \draw (5.5,-1) -- (5,-0.5);
  \draw (5.5,-1) -- (6,-0.5);
  \draw (3,-3.5) -- (2,-2.5);
  \draw (3,-3.5) -- (5.5,-1);
  \draw (8.5,-1) -- (8,-0.5);
  \draw (8.5,-1) -- (9,-0.5);
  \draw (8,-1.5) -- (7,-0.5);
  \draw (8,-1.5) -- (8.5,-1);
  \draw (4.5,-5) -- (3,-3.5);
  \draw (4.5,-5) -- (8,-1.5);
  \draw (11.5,-1) -- (11,-0.5);
  \draw (11.5,-1) -- (12,-0.5);
  \draw (11,-1.5) -- (10,-0.5);
  \draw (11,-1.5) -- (11.5,-1);
  \draw (6,-6.5) -- (4.5,-5);
  \draw (6,-6.5) -- (11,-1.5);
  \fill (0,0.5) circle (2.5pt) (0,-0.5) circle (2.5pt);
  \fill (1,0.5) circle (2.5pt) (1,-0.5) circle (2.5pt);
  \fill (2,0.5) circle (2.5pt) (2,-0.5) circle (2.5pt);
  \fill (3,0.5) circle (2.5pt) (3,-0.5) circle (2.5pt);
  \fill (4,0.5) circle (2.5pt) (4,-0.5) circle (2.5pt);
  \fill (5,0.5) circle (2.5pt) (5,-0.5) circle (2.5pt);
  \fill (6,0.5) circle (2.5pt) (6,-0.5) circle (2.5pt);
  \fill (7,0.5) circle (2.5pt) (7,-0.5) circle (2.5pt);
  \fill (8,0.5) circle (2.5pt) (8,-0.5) circle (2.5pt);
  \fill (9,0.5) circle (2.5pt) (9,-0.5) circle (2.5pt);
  \fill (10,0.5) circle (2.5pt) (10,-0.5) circle (2.5pt);
  \fill (11,0.5) circle (2.5pt) (11,-0.5) circle (2.5pt);
  \fill (12,0.5) circle (2.5pt) (12,-0.5) circle (2.5pt);
\end{tikzpicture}} &
\resizebox{0.35\linewidth}{!}{\begin{tikzpicture}[scale=0.32, line cap=round, line join=round, >=stealth, thick]
  \draw[rounded corners=4pt] (11.682,11.682) -- (0.318,0.318) .. controls (0,0) .. (0.318,-0.318) -- (0.682,-0.682);
  \draw[rounded corners=4pt] (11.696,11.668) -- (1.304,0.332) .. controls (1,0) .. (1,-0.45) -- (1,-0.55);
  \draw[rounded corners=4pt] (12.318,11.682) -- (12.682,11.318);
  \draw[rounded corners=4pt] (12.682,10.682) -- (2.318,0.318) .. controls (2,0) .. (1.682,-0.318) -- (1.318,-0.682);
  \draw[rounded corners=4pt] (12.697,10.667) -- (3.303,0.333) .. controls (3,0) .. (3.109,-0.437) -- (3.891,-3.563);
  \draw[rounded corners=4pt] (13.318,10.682) -- (13.682,10.318);
  \draw[rounded corners=4pt] (13.682,9.682) -- (4.318,0.318) .. controls (4,0) .. (4.318,-0.318) -- (4.682,-0.682);
  \draw[rounded corners=4pt] (13.699,9.666) -- (5.301,0.334) .. controls (5,0) .. (5,-0.45) -- (5,-0.55);
  \draw[rounded corners=4pt] (14.318,9.682) -- (14.682,9.318);
  \draw[rounded corners=4pt] (14.682,8.682) -- (8.318,2.318);
  \draw[rounded corners=4pt] (14.817,8.589) -- (11.183,0.411) .. controls (11,0) .. (11,-0.45) -- (11,-0.55);
  \draw[rounded corners=4pt] (15.318,8.682) -- (17.682,6.318);
  \draw[rounded corners=4pt] (7.682,1.682) -- (6.318,0.318) .. controls (6,0) .. (5.682,-0.318) -- (5.318,-0.682);
  \draw[rounded corners=4pt] (7.799,1.598) -- (7.201,0.402) .. controls (7,0) .. (6.799,-0.402) -- (6.201,-1.598);
  \draw[rounded corners=4pt] (8.318,1.682) -- (8.682,1.318);
  \draw[rounded corners=4pt] (8.682,0.682) -- (8.318,0.318) .. controls (8,0) .. (7.682,-0.318) -- (6.318,-1.682);
  \draw[rounded corners=4pt] (9,0.55) -- (9,0.45) .. controls (9,0) .. (8.799,-0.402) -- (6.201,-5.598);
  \draw[rounded corners=4pt] (9.318,0.682) -- (9.682,0.318) .. controls (10,0) .. (10.318,-0.318) -- (10.682,-0.682);
  \draw[rounded corners=4pt] (17.682,5.682) -- (17.318,5.318);
  \draw[rounded corners=4pt] (18.288,5.654) -- (22.712,0.346) .. controls (23,0) .. (23,-0.45) -- (23,-0.55);
  \draw[rounded corners=4pt] (18.318,5.682) -- (23.682,0.318) .. controls (24,0) .. (23.682,-0.318) -- (23.318,-0.682);
  \draw[rounded corners=4pt] (16.682,4.682) -- (12.318,0.318) .. controls (12,0) .. (11.682,-0.318) -- (11.318,-0.682);
  \draw[rounded corners=4pt] (16.719,4.649) -- (13.281,0.351) .. controls (13,0) .. (12.817,-0.411) -- (9.183,-8.589);
  \draw[rounded corners=4pt] (17.318,4.682) -- (17.682,4.318);
  \draw[rounded corners=4pt] (17.682,3.682) -- (17.318,3.318);
  \draw[rounded corners=4pt] (18.27,3.64) -- (20.73,0.36) .. controls (21,0) .. (21.201,-0.402) -- (21.799,-1.598);
  \draw[rounded corners=4pt] (18.318,3.682) -- (21.682,0.318) .. controls (22,0) .. (22.318,-0.318) -- (22.682,-0.682);
  \draw[rounded corners=4pt] (16.682,2.682) -- (15.318,1.318);
  \draw[rounded corners=4pt] (17,2.55) -- (17,0.45) .. controls (17,0) .. (17,-0.45) -- (17,-0.55);
  \draw[rounded corners=4pt] (17.318,2.682) -- (18.682,1.318);
  \draw[rounded corners=4pt] (14.682,0.682) -- (14.318,0.318) .. controls (14,0) .. (14.318,-0.318) -- (15.682,-1.682);
  \draw[rounded corners=4pt] (15,0.55) -- (15,0.45) .. controls (15,0) .. (15.201,-0.402) -- (15.799,-1.598);
  \draw[rounded corners=4pt] (15.318,0.682) -- (15.682,0.318) .. controls (16,0) .. (16.318,-0.318) -- (16.682,-0.682);
  \draw[rounded corners=4pt] (18.682,0.682) -- (18.318,0.318) .. controls (18,0) .. (17.682,-0.318) -- (17.318,-0.682);
  \draw[rounded corners=4pt] (19,0.55) -- (19,0.45) .. controls (19,0) .. (18.773,-0.389) -- (12.227,-11.611);
  \draw[rounded corners=4pt] (19.318,0.682) -- (19.682,0.318) .. controls (20,0) .. (20.318,-0.318) -- (21.682,-1.682);
  \draw[rounded corners=4pt] (11.682,-11.682) -- (9.318,-9.318);
  \draw[rounded corners=4pt] (12.318,-11.682) -- (21.682,-2.318);
  \draw[rounded corners=4pt] (8.682,-8.682) -- (6.318,-6.318);
  \draw[rounded corners=4pt] (9.318,-8.682) -- (15.682,-2.318);
  \draw[rounded corners=4pt] (5.682,-5.682) -- (4.318,-4.318);
  \draw[rounded corners=4pt] (6.318,-5.682) -- (10.682,-1.318);
  \draw[rounded corners=4pt] (3.682,-3.682) -- (1.318,-1.318);
  \draw[rounded corners=4pt] (4.318,-3.682) -- (5.682,-2.318);
  \draw[rounded corners=4pt] (5.682,-1.682) -- (5.318,-1.318);
  \draw[rounded corners=4pt] (16.318,-1.682) -- (16.682,-1.318);
  \draw[rounded corners=4pt] (22.318,-1.682) -- (22.682,-1.318);
  \draw[rounded corners=4pt] (12,12.45) -- (12,13) -- (-1,13) -- (-1,-13) -- (12,-13) -- (12,-12.45);
  \draw (12,12.45) .. controls (12,12) .. (11.696,11.668);
  \draw (12.682,11.318) .. controls (13,11) .. (12.697,10.667);
  \draw (13.682,10.318) .. controls (14,10) .. (13.699,9.666);
  \draw (14.682,9.318) .. controls (15,9) .. (14.817,8.589);
  \draw (8.318,2.318) .. controls (8,2) .. (7.799,1.598);
  \draw (8.682,1.318) .. controls (9,1) .. (9,0.55);
  \draw (17.682,6.318) .. controls (18,6) .. (18.288,5.654);
  \draw (17.318,5.318) .. controls (17,5) .. (16.719,4.649);
  \draw (17.682,4.318) .. controls (18,4) .. (18.27,3.64);
  \draw (17.318,3.318) .. controls (17,3) .. (17,2.55);
  \draw (15.318,1.318) .. controls (15,1) .. (15,0.55);
  \draw (18.682,1.318) .. controls (19,1) .. (19,0.55);
  \draw (12,-12.45) .. controls (12,-12) .. (12.227,-11.611);
  \draw (9.318,-9.318) .. controls (9,-9) .. (9.183,-8.589);
  \draw (6.318,-6.318) .. controls (6,-6) .. (6.201,-5.598);
  \draw (4.318,-4.318) .. controls (4,-4) .. (3.891,-3.563);
  \draw (1.318,-1.318) .. controls (1,-1) .. (1,-0.55);
  \draw (5.682,-2.318) .. controls (6,-2) .. (6.201,-1.598);
  \draw (5.318,-1.318) .. controls (5,-1) .. (5,-0.55);
  \draw (10.682,-1.318) .. controls (11,-1) .. (11,-0.55);
  \draw (15.682,-2.318) .. controls (16,-2) .. (15.799,-1.598);
  \draw (16.682,-1.318) .. controls (17,-1) .. (17,-0.55);
  \draw (21.682,-2.318) .. controls (22,-2) .. (21.799,-1.598);
  \draw (22.682,-1.318) .. controls (23,-1) .. (23,-0.55);
  \draw[white, line width=6pt, line cap=butt] (11.682,11.682) .. controls (12,12) .. (12.318,11.682);
  \draw (11.682,11.682) .. controls (12,12) .. (12.318,11.682);
  \draw[white, line width=6pt, line cap=butt] (12.682,10.682) .. controls (13,11) .. (13.318,10.682);
  \draw (12.682,10.682) .. controls (13,11) .. (13.318,10.682);
  \draw[white, line width=6pt, line cap=butt] (13.682,9.682) .. controls (14,10) .. (14.318,9.682);
  \draw (13.682,9.682) .. controls (14,10) .. (14.318,9.682);
  \draw[white, line width=6pt, line cap=butt] (14.682,8.682) .. controls (15,9) .. (15.318,8.682);
  \draw (14.682,8.682) .. controls (15,9) .. (15.318,8.682);
  \draw[white, line width=6pt, line cap=butt] (7.682,1.682) .. controls (8,2) .. (8.318,1.682);
  \draw (7.682,1.682) .. controls (8,2) .. (8.318,1.682);
  \draw[white, line width=6pt, line cap=butt] (8.682,0.682) .. controls (9,1) .. (9.318,0.682);
  \draw (8.682,0.682) .. controls (9,1) .. (9.318,0.682);
  \draw[white, line width=6pt, line cap=butt] (17.682,5.682) .. controls (18,6) .. (18.318,5.682);
  \draw (17.682,5.682) .. controls (18,6) .. (18.318,5.682);
  \draw[white, line width=6pt, line cap=butt] (16.682,4.682) .. controls (17,5) .. (17.318,4.682);
  \draw (16.682,4.682) .. controls (17,5) .. (17.318,4.682);
  \draw[white, line width=6pt, line cap=butt] (17.682,3.682) .. controls (18,4) .. (18.318,3.682);
  \draw (17.682,3.682) .. controls (18,4) .. (18.318,3.682);
  \draw[white, line width=6pt, line cap=butt] (16.682,2.682) .. controls (17,3) .. (17.318,2.682);
  \draw (16.682,2.682) .. controls (17,3) .. (17.318,2.682);
  \draw[white, line width=6pt, line cap=butt] (14.682,0.682) .. controls (15,1) .. (15.318,0.682);
  \draw (14.682,0.682) .. controls (15,1) .. (15.318,0.682);
  \draw[white, line width=6pt, line cap=butt] (18.682,0.682) .. controls (19,1) .. (19.318,0.682);
  \draw (18.682,0.682) .. controls (19,1) .. (19.318,0.682);
  \draw[white, line width=6pt, line cap=butt] (11.682,-11.682) .. controls (12,-12) .. (12.318,-11.682);
  \draw (11.682,-11.682) .. controls (12,-12) .. (12.318,-11.682);
  \draw[white, line width=6pt, line cap=butt] (8.682,-8.682) .. controls (9,-9) .. (9.318,-8.682);
  \draw (8.682,-8.682) .. controls (9,-9) .. (9.318,-8.682);
  \draw[white, line width=6pt, line cap=butt] (5.682,-5.682) .. controls (6,-6) .. (6.318,-5.682);
  \draw (5.682,-5.682) .. controls (6,-6) .. (6.318,-5.682);
  \draw[white, line width=6pt, line cap=butt] (3.682,-3.682) .. controls (4,-4) .. (4.318,-3.682);
  \draw (3.682,-3.682) .. controls (4,-4) .. (4.318,-3.682);
  \draw[white, line width=6pt, line cap=butt] (0.682,-0.682) .. controls (1,-1) .. (1.318,-0.682);
  \draw (0.682,-0.682) .. controls (1,-1) .. (1.318,-0.682);
  \draw[white, line width=6pt, line cap=butt] (5.682,-1.682) .. controls (6,-2) .. (6.318,-1.682);
  \draw (5.682,-1.682) .. controls (6,-2) .. (6.318,-1.682);
  \draw[white, line width=6pt, line cap=butt] (4.682,-0.682) .. controls (5,-1) .. (5.318,-0.682);
  \draw (4.682,-0.682) .. controls (5,-1) .. (5.318,-0.682);
  \draw[white, line width=6pt, line cap=butt] (10.682,-0.682) .. controls (11,-1) .. (11.318,-0.682);
  \draw (10.682,-0.682) .. controls (11,-1) .. (11.318,-0.682);
  \draw[white, line width=6pt, line cap=butt] (15.682,-1.682) .. controls (16,-2) .. (16.318,-1.682);
  \draw (15.682,-1.682) .. controls (16,-2) .. (16.318,-1.682);
  \draw[white, line width=6pt, line cap=butt] (16.682,-0.682) .. controls (17,-1) .. (17.318,-0.682);
  \draw (16.682,-0.682) .. controls (17,-1) .. (17.318,-0.682);
  \draw[white, line width=6pt, line cap=butt] (21.682,-1.682) .. controls (22,-2) .. (22.318,-1.682);
  \draw (21.682,-1.682) .. controls (22,-2) .. (22.318,-1.682);
  \draw[white, line width=6pt, line cap=butt] (22.682,-0.682) .. controls (23,-1) .. (23.318,-0.682);
  \draw (22.682,-0.682) .. controls (23,-1) .. (23.318,-0.682);
\end{tikzpicture}} \\[1mm]
\resizebox{0.28\linewidth}{!}{\begin{tikzpicture}[scale=0.32, line cap=round]
  \draw (1.5,1) -- (1,0.5);
  \draw (1.5,1) -- (2,0.5);
  \draw (2,1.5) -- (1.5,1);
  \draw (2,1.5) -- (3,0.5);
  \draw (2.5,2) -- (2,1.5);
  \draw (2.5,2) -- (4,0.5);
  \draw (2,2.5) -- (0,0.5);
  \draw (2,2.5) -- (2.5,2);
  \draw (7.5,1) -- (7,0.5);
  \draw (7.5,1) -- (8,0.5);
  \draw (7,1.5) -- (6,0.5);
  \draw (7,1.5) -- (7.5,1);
  \draw (6.5,2) -- (5,0.5);
  \draw (6.5,2) -- (7,1.5);
  \draw (7,2.5) -- (6.5,2);
  \draw (7,2.5) -- (9,0.5);
  \draw (4.5,5) -- (2,2.5);
  \draw (4.5,5) -- (7,2.5);
  \draw (5,5.5) -- (4.5,5);
  \draw (5,5.5) -- (10,0.5);
  \draw (0.5,-1) -- (0,-0.5);
  \draw (0.5,-1) -- (1,-0.5);
  \draw (4.5,-1) -- (4,-0.5);
  \draw (4.5,-1) -- (5,-0.5);
  \draw (4,-1.5) -- (3,-0.5);
  \draw (4,-1.5) -- (4.5,-1);
  \draw (4.5,-2) -- (4,-1.5);
  \draw (4.5,-2) -- (6,-0.5);
  \draw (5,-2.5) -- (4.5,-2);
  \draw (5,-2.5) -- (7,-0.5);
  \draw (4.5,-3) -- (2,-0.5);
  \draw (4.5,-3) -- (5,-2.5);
  \draw (3.5,-4) -- (0.5,-1);
  \draw (3.5,-4) -- (4.5,-3);
  \draw (9.5,-1) -- (9,-0.5);
  \draw (9.5,-1) -- (10,-0.5);
  \draw (9,-1.5) -- (8,-0.5);
  \draw (9,-1.5) -- (9.5,-1);
  \draw (5,-5.5) -- (3.5,-4);
  \draw (5,-5.5) -- (9,-1.5);
  \fill (0,0.5) circle (2.5pt) (0,-0.5) circle (2.5pt);
  \fill (1,0.5) circle (2.5pt) (1,-0.5) circle (2.5pt);
  \fill (2,0.5) circle (2.5pt) (2,-0.5) circle (2.5pt);
  \fill (3,0.5) circle (2.5pt) (3,-0.5) circle (2.5pt);
  \fill (4,0.5) circle (2.5pt) (4,-0.5) circle (2.5pt);
  \fill (5,0.5) circle (2.5pt) (5,-0.5) circle (2.5pt);
  \fill (6,0.5) circle (2.5pt) (6,-0.5) circle (2.5pt);
  \fill (7,0.5) circle (2.5pt) (7,-0.5) circle (2.5pt);
  \fill (8,0.5) circle (2.5pt) (8,-0.5) circle (2.5pt);
  \fill (9,0.5) circle (2.5pt) (9,-0.5) circle (2.5pt);
  \fill (10,0.5) circle (2.5pt) (10,-0.5) circle (2.5pt);
\end{tikzpicture}} &
\resizebox{0.35\linewidth}{!}{\begin{tikzpicture}[scale=0.32, line cap=round, line join=round, >=stealth, thick]
  \draw[rounded corners=4pt] (9.682,9.682) -- (9.318,9.318);
  \draw[rounded corners=4pt] (10.301,9.666) -- (18.699,0.334) .. controls (19,0) .. (19,-0.45) -- (19,-0.55);
  \draw[rounded corners=4pt] (10.318,9.682) -- (19.682,0.318) .. controls (20,0) .. (19.682,-0.318) -- (19.318,-0.682);
  \draw[rounded corners=4pt] (8.682,8.682) -- (4.318,4.318);
  \draw[rounded corners=4pt] (9,8.55) -- (9,0.45) .. controls (9,0) .. (9,-0.45) -- (9,-0.55);
  \draw[rounded corners=4pt] (9.318,8.682) -- (13.682,4.318);
  \draw[rounded corners=4pt] (3.682,3.682) -- (0.318,0.318) .. controls (0,0) .. (0.318,-0.318) -- (0.682,-0.682);
  \draw[rounded corners=4pt] (3.73,3.64) -- (1.27,0.36) .. controls (1,0) .. (1,-0.45) -- (1,-0.55);
  \draw[rounded corners=4pt] (4.318,3.682) -- (4.682,3.318);
  \draw[rounded corners=4pt] (4.682,2.682) -- (4.318,2.318);
  \draw[rounded corners=4pt] (5.25,2.626) -- (6.75,0.374) .. controls (7,0) .. (7.201,-0.402) -- (7.799,-1.598);
  \draw[rounded corners=4pt] (5.318,2.682) -- (7.682,0.318) .. controls (8,0) .. (8.318,-0.318) -- (8.682,-0.682);
  \draw[rounded corners=4pt] (3.682,1.682) -- (3.318,1.318);
  \draw[rounded corners=4pt] (4.201,1.598) -- (4.799,0.402) .. controls (5,0) .. (5.281,-0.351) -- (8.719,-4.649);
  \draw[rounded corners=4pt] (4.318,1.682) -- (5.682,0.318) .. controls (6,0) .. (6.318,-0.318) -- (7.682,-1.682);
  \draw[rounded corners=4pt] (2.682,0.682) -- (2.318,0.318) .. controls (2,0) .. (1.682,-0.318) -- (1.318,-0.682);
  \draw[rounded corners=4pt] (3,0.55) -- (3,0.45) .. controls (3,0) .. (3.223,-0.391) -- (6.777,-6.609);
  \draw[rounded corners=4pt] (3.318,0.682) -- (3.682,0.318) .. controls (4,0) .. (4.318,-0.318) -- (8.682,-4.682);
  \draw[rounded corners=4pt] (13.682,3.682) -- (13.318,3.318);
  \draw[rounded corners=4pt] (14.27,3.64) -- (16.73,0.36) .. controls (17,0) .. (17.201,-0.402) -- (17.799,-1.598);
  \draw[rounded corners=4pt] (14.318,3.682) -- (17.682,0.318) .. controls (18,0) .. (18.318,-0.318) -- (18.682,-0.682);
  \draw[rounded corners=4pt] (12.682,2.682) -- (10.318,0.318) .. controls (10,0) .. (9.682,-0.318) -- (9.318,-0.682);
  \draw[rounded corners=4pt] (12.75,2.626) -- (11.25,0.374) .. controls (11,0) .. (10.75,-0.374) -- (9.25,-2.626);
  \draw[rounded corners=4pt] (13.318,2.682) -- (13.682,2.318);
  \draw[rounded corners=4pt] (13.682,1.682) -- (12.318,0.318) .. controls (12,0) .. (11.682,-0.318) -- (9.318,-2.682);
  \draw[rounded corners=4pt] (13.799,1.598) -- (13.201,0.402) .. controls (13,0) .. (12.73,-0.36) -- (10.27,-3.64);
  \draw[rounded corners=4pt] (14.318,1.682) -- (14.682,1.318);
  \draw[rounded corners=4pt] (14.682,0.682) -- (14.318,0.318) .. controls (14,0) .. (13.682,-0.318) -- (10.318,-3.682);
  \draw[rounded corners=4pt] (15,0.55) -- (15,0.45) .. controls (15,0) .. (14.799,-0.402) -- (10.201,-9.598);
  \draw[rounded corners=4pt] (15.318,0.682) -- (15.682,0.318) .. controls (16,0) .. (16.318,-0.318) -- (17.682,-1.682);
  \draw[rounded corners=4pt] (9.682,-9.682) -- (7.318,-7.318);
  \draw[rounded corners=4pt] (10.318,-9.682) -- (17.682,-2.318);
  \draw[rounded corners=4pt] (6.682,-6.682) -- (1.318,-1.318);
  \draw[rounded corners=4pt] (7.318,-6.682) -- (8.682,-5.318);
  \draw[rounded corners=4pt] (9.318,-4.682) -- (9.682,-4.318);
  \draw[rounded corners=4pt] (9.682,-3.682) -- (9.318,-3.318);
  \draw[rounded corners=4pt] (8.682,-2.682) -- (8.318,-2.318);
  \draw[rounded corners=4pt] (8.318,-1.682) -- (8.682,-1.318);
  \draw[rounded corners=4pt] (18.318,-1.682) -- (18.682,-1.318);
  \draw[rounded corners=4pt] (10,10.45) -- (10,11) -- (-1,11) -- (-1,-11) -- (10,-11) -- (10,-10.45);
  \draw (10,10.45) .. controls (10,10) .. (10.301,9.666);
  \draw (9.318,9.318) .. controls (9,9) .. (9,8.55);
  \draw (4.318,4.318) .. controls (4,4) .. (3.73,3.64);
  \draw (4.682,3.318) .. controls (5,3) .. (5.25,2.626);
  \draw (4.318,2.318) .. controls (4,2) .. (4.201,1.598);
  \draw (3.318,1.318) .. controls (3,1) .. (3,0.55);
  \draw (13.682,4.318) .. controls (14,4) .. (14.27,3.64);
  \draw (13.318,3.318) .. controls (13,3) .. (12.75,2.626);
  \draw (13.682,2.318) .. controls (14,2) .. (13.799,1.598);
  \draw (14.682,1.318) .. controls (15,1) .. (15,0.55);
  \draw (10,-10.45) .. controls (10,-10) .. (10.201,-9.598);
  \draw (7.318,-7.318) .. controls (7,-7) .. (6.777,-6.609);
  \draw (1.318,-1.318) .. controls (1,-1) .. (1,-0.55);
  \draw (8.682,-5.318) .. controls (9,-5) .. (8.719,-4.649);
  \draw (9.682,-4.318) .. controls (10,-4) .. (10.27,-3.64);
  \draw (9.318,-3.318) .. controls (9,-3) .. (9.25,-2.626);
  \draw (8.318,-2.318) .. controls (8,-2) .. (7.799,-1.598);
  \draw (8.682,-1.318) .. controls (9,-1) .. (9,-0.55);
  \draw (17.682,-2.318) .. controls (18,-2) .. (17.799,-1.598);
  \draw (18.682,-1.318) .. controls (19,-1) .. (19,-0.55);
  \draw[white, line width=6pt, line cap=butt] (9.682,9.682) .. controls (10,10) .. (10.318,9.682);
  \draw (9.682,9.682) .. controls (10,10) .. (10.318,9.682);
  \draw[white, line width=6pt, line cap=butt] (8.682,8.682) .. controls (9,9) .. (9.318,8.682);
  \draw (8.682,8.682) .. controls (9,9) .. (9.318,8.682);
  \draw[white, line width=6pt, line cap=butt] (3.682,3.682) .. controls (4,4) .. (4.318,3.682);
  \draw (3.682,3.682) .. controls (4,4) .. (4.318,3.682);
  \draw[white, line width=6pt, line cap=butt] (4.682,2.682) .. controls (5,3) .. (5.318,2.682);
  \draw (4.682,2.682) .. controls (5,3) .. (5.318,2.682);
  \draw[white, line width=6pt, line cap=butt] (3.682,1.682) .. controls (4,2) .. (4.318,1.682);
  \draw (3.682,1.682) .. controls (4,2) .. (4.318,1.682);
  \draw[white, line width=6pt, line cap=butt] (2.682,0.682) .. controls (3,1) .. (3.318,0.682);
  \draw (2.682,0.682) .. controls (3,1) .. (3.318,0.682);
  \draw[white, line width=6pt, line cap=butt] (13.682,3.682) .. controls (14,4) .. (14.318,3.682);
  \draw (13.682,3.682) .. controls (14,4) .. (14.318,3.682);
  \draw[white, line width=6pt, line cap=butt] (12.682,2.682) .. controls (13,3) .. (13.318,2.682);
  \draw (12.682,2.682) .. controls (13,3) .. (13.318,2.682);
  \draw[white, line width=6pt, line cap=butt] (13.682,1.682) .. controls (14,2) .. (14.318,1.682);
  \draw (13.682,1.682) .. controls (14,2) .. (14.318,1.682);
  \draw[white, line width=6pt, line cap=butt] (14.682,0.682) .. controls (15,1) .. (15.318,0.682);
  \draw (14.682,0.682) .. controls (15,1) .. (15.318,0.682);
  \draw[white, line width=6pt, line cap=butt] (9.682,-9.682) .. controls (10,-10) .. (10.318,-9.682);
  \draw (9.682,-9.682) .. controls (10,-10) .. (10.318,-9.682);
  \draw[white, line width=6pt, line cap=butt] (6.682,-6.682) .. controls (7,-7) .. (7.318,-6.682);
  \draw (6.682,-6.682) .. controls (7,-7) .. (7.318,-6.682);
  \draw[white, line width=6pt, line cap=butt] (0.682,-0.682) .. controls (1,-1) .. (1.318,-0.682);
  \draw (0.682,-0.682) .. controls (1,-1) .. (1.318,-0.682);
  \draw[white, line width=6pt, line cap=butt] (8.682,-4.682) .. controls (9,-5) .. (9.318,-4.682);
  \draw (8.682,-4.682) .. controls (9,-5) .. (9.318,-4.682);
  \draw[white, line width=6pt, line cap=butt] (9.682,-3.682) .. controls (10,-4) .. (10.318,-3.682);
  \draw (9.682,-3.682) .. controls (10,-4) .. (10.318,-3.682);
  \draw[white, line width=6pt, line cap=butt] (8.682,-2.682) .. controls (9,-3) .. (9.318,-2.682);
  \draw (8.682,-2.682) .. controls (9,-3) .. (9.318,-2.682);
  \draw[white, line width=6pt, line cap=butt] (7.682,-1.682) .. controls (8,-2) .. (8.318,-1.682);
  \draw (7.682,-1.682) .. controls (8,-2) .. (8.318,-1.682);
  \draw[white, line width=6pt, line cap=butt] (8.682,-0.682) .. controls (9,-1) .. (9.318,-0.682);
  \draw (8.682,-0.682) .. controls (9,-1) .. (9.318,-0.682);
  \draw[white, line width=6pt, line cap=butt] (17.682,-1.682) .. controls (18,-2) .. (18.318,-1.682);
  \draw (17.682,-1.682) .. controls (18,-2) .. (18.318,-1.682);
  \draw[white, line width=6pt, line cap=butt] (18.682,-0.682) .. controls (19,-1) .. (19.318,-0.682);
  \draw (18.682,-0.682) .. controls (19,-1) .. (19.318,-0.682);
\end{tikzpicture}}
\end{tabular}
\caption{Top: $10_{45}$ with 13 leaves. Bottom: Sweeney's 11-leaf diagram, which has the HOMFLY-PT
polynomial of $10_{45}$ but is a different knot.}
\label{fig:sw-1045}
\end{figure}
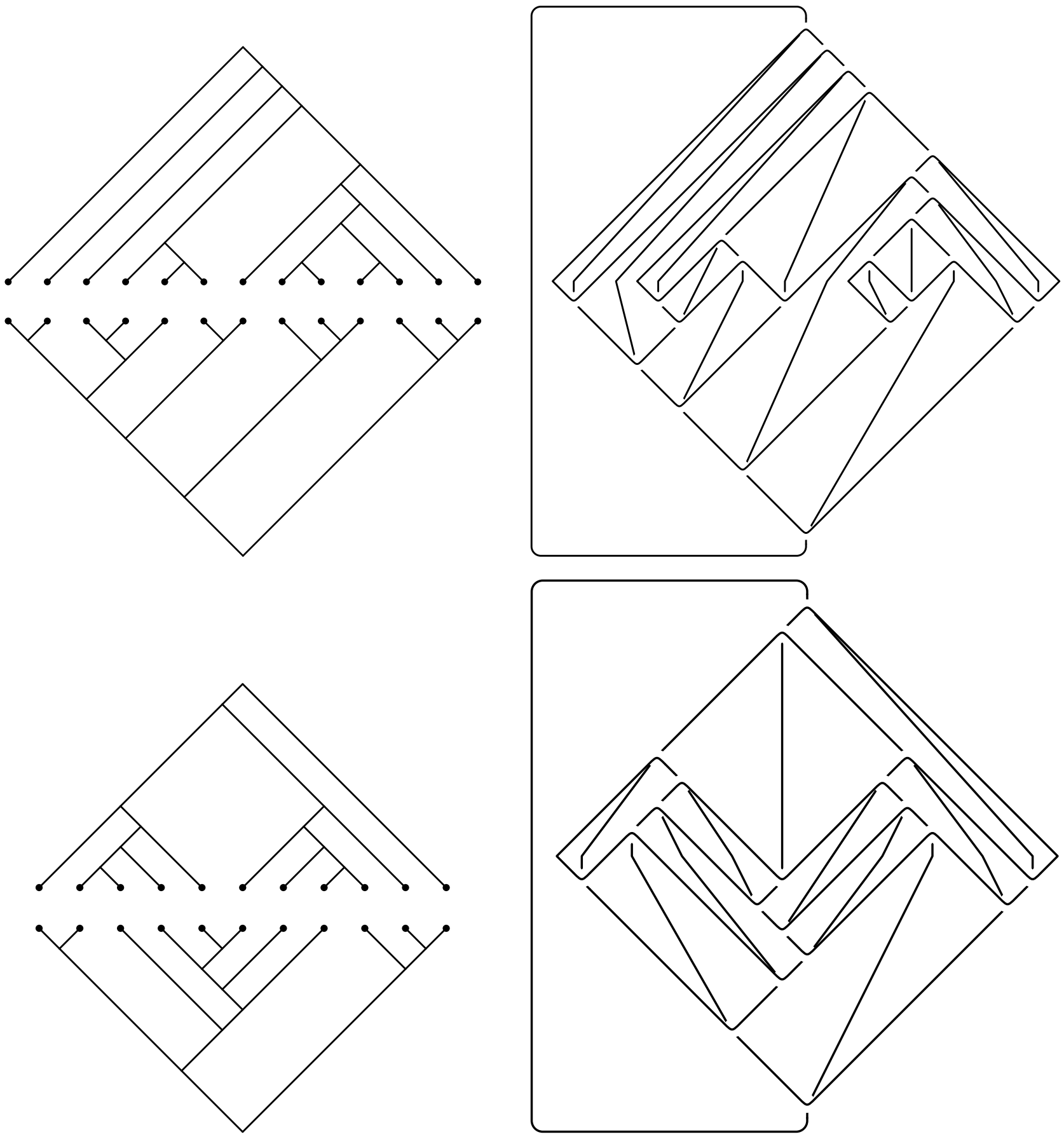

\begin{figure}[ht]
\centering
\begin{tabular}{cc}
\resizebox{0.28\linewidth}{!}{\begin{tikzpicture}[scale=0.32, line cap=round]
  \draw (2.5,1) -- (2,0.5);
  \draw (2.5,1) -- (3,0.5);
  \draw (3,1.5) -- (2.5,1);
  \draw (3,1.5) -- (4,0.5);
  \draw (2.5,2) -- (1,0.5);
  \draw (2.5,2) -- (3,1.5);
  \draw (5.5,1) -- (5,0.5);
  \draw (5.5,1) -- (6,0.5);
  \draw (3.5,3) -- (2.5,2);
  \draw (3.5,3) -- (5.5,1);
  \draw (7.5,1) -- (7,0.5);
  \draw (7.5,1) -- (8,0.5);
  \draw (4.5,4) -- (3.5,3);
  \draw (4.5,4) -- (7.5,1);
  \draw (5,4.5) -- (4.5,4);
  \draw (5,4.5) -- (9,0.5);
  \draw (4.5,5) -- (0,0.5);
  \draw (4.5,5) -- (5,4.5);
  \draw (0.5,-1) -- (0,-0.5);
  \draw (0.5,-1) -- (1,-0.5);
  \draw (1,-1.5) -- (0.5,-1);
  \draw (1,-1.5) -- (2,-0.5);
  \draw (4.5,-1) -- (4,-0.5);
  \draw (4.5,-1) -- (5,-0.5);
  \draw (4,-1.5) -- (3,-0.5);
  \draw (4,-1.5) -- (4.5,-1);
  \draw (6.5,-1) -- (6,-0.5);
  \draw (6.5,-1) -- (7,-0.5);
  \draw (5,-2.5) -- (4,-1.5);
  \draw (5,-2.5) -- (6.5,-1);
  \draw (8.5,-1) -- (8,-0.5);
  \draw (8.5,-1) -- (9,-0.5);
  \draw (6,-3.5) -- (5,-2.5);
  \draw (6,-3.5) -- (8.5,-1);
  \draw (4.5,-5) -- (1,-1.5);
  \draw (4.5,-5) -- (6,-3.5);
  \fill (0,0.5) circle (2.5pt) (0,-0.5) circle (2.5pt);
  \fill (1,0.5) circle (2.5pt) (1,-0.5) circle (2.5pt);
  \fill (2,0.5) circle (2.5pt) (2,-0.5) circle (2.5pt);
  \fill (3,0.5) circle (2.5pt) (3,-0.5) circle (2.5pt);
  \fill (4,0.5) circle (2.5pt) (4,-0.5) circle (2.5pt);
  \fill (5,0.5) circle (2.5pt) (5,-0.5) circle (2.5pt);
  \fill (6,0.5) circle (2.5pt) (6,-0.5) circle (2.5pt);
  \fill (7,0.5) circle (2.5pt) (7,-0.5) circle (2.5pt);
  \fill (8,0.5) circle (2.5pt) (8,-0.5) circle (2.5pt);
  \fill (9,0.5) circle (2.5pt) (9,-0.5) circle (2.5pt);
\end{tikzpicture}} &
\resizebox{0.35\linewidth}{!}{\begin{tikzpicture}[scale=0.32, line cap=round, line join=round, >=stealth, thick]
  \draw[rounded corners=4pt] (8.682,8.682) -- (0.318,0.318) .. controls (0,0) .. (0.318,-0.318) -- (0.682,-0.682);
  \draw[rounded corners=4pt] (8.701,8.664) -- (1.299,0.336) .. controls (1,0) .. (1,-0.45) -- (1,-0.55);
  \draw[rounded corners=4pt] (9.318,8.682) -- (9.682,8.318);
  \draw[rounded corners=4pt] (9.682,7.682) -- (9.318,7.318);
  \draw[rounded corners=4pt] (10.296,7.661) -- (16.704,0.339) .. controls (17,0) .. (17,-0.45) -- (17,-0.55);
  \draw[rounded corners=4pt] (10.318,7.682) -- (17.682,0.318) .. controls (18,0) .. (17.682,-0.318) -- (17.318,-0.682);
  \draw[rounded corners=4pt] (8.682,6.682) -- (7.318,5.318);
  \draw[rounded corners=4pt] (9.223,6.609) -- (12.777,0.391) .. controls (13,0) .. (13,-0.45) -- (13,-0.55);
  \draw[rounded corners=4pt] (9.318,6.682) -- (14.682,1.318);
  \draw[rounded corners=4pt] (6.682,4.682) -- (5.318,3.318);
  \draw[rounded corners=4pt] (7.167,4.582) -- (8.833,0.418) .. controls (9,0) .. (9,-0.45) -- (9,-0.55);
  \draw[rounded corners=4pt] (7.318,4.682) -- (10.682,1.318);
  \draw[rounded corners=4pt] (4.682,2.682) -- (2.318,0.318) .. controls (2,0) .. (1.682,-0.318) -- (1.318,-0.682);
  \draw[rounded corners=4pt] (4.75,2.626) -- (3.25,0.374) .. controls (3,0) .. (2.799,-0.402) -- (2.201,-1.598);
  \draw[rounded corners=4pt] (5.318,2.682) -- (5.682,2.318);
  \draw[rounded corners=4pt] (5.682,1.682) -- (5.318,1.318);
  \draw[rounded corners=4pt] (6.201,1.598) -- (6.799,0.402) .. controls (7,0) .. (7.201,-0.402) -- (7.799,-1.598);
  \draw[rounded corners=4pt] (6.318,1.682) -- (7.682,0.318) .. controls (8,0) .. (8.318,-0.318) -- (8.682,-0.682);
  \draw[rounded corners=4pt] (4.682,0.682) -- (4.318,0.318) .. controls (4,0) .. (3.682,-0.318) -- (2.318,-1.682);
  \draw[rounded corners=4pt] (5,0.55) -- (5,0.45) .. controls (5,0) .. (5.183,-0.411) -- (8.817,-8.589);
  \draw[rounded corners=4pt] (5.318,0.682) -- (5.682,0.318) .. controls (6,0) .. (6.318,-0.318) -- (7.682,-1.682);
  \draw[rounded corners=4pt] (10.682,0.682) -- (10.318,0.318) .. controls (10,0) .. (9.682,-0.318) -- (9.318,-0.682);
  \draw[rounded corners=4pt] (11,0.55) -- (11,0.45) .. controls (11,0) .. (10.891,-0.437) -- (10.109,-3.563);
  \draw[rounded corners=4pt] (11.318,0.682) -- (11.682,0.318) .. controls (12,0) .. (12.318,-0.318) -- (12.682,-0.682);
  \draw[rounded corners=4pt] (14.682,0.682) -- (14.318,0.318) .. controls (14,0) .. (13.682,-0.318) -- (13.318,-0.682);
  \draw[rounded corners=4pt] (15,0.55) -- (15,0.45) .. controls (15,0) .. (14.799,-0.402) -- (12.201,-5.598);
  \draw[rounded corners=4pt] (15.318,0.682) -- (15.682,0.318) .. controls (16,0) .. (16.318,-0.318) -- (16.682,-0.682);
  \draw[rounded corners=4pt] (8.682,-8.682) -- (2.318,-2.318);
  \draw[rounded corners=4pt] (9.318,-8.682) -- (11.682,-6.318);
  \draw[rounded corners=4pt] (1.682,-1.682) -- (1.318,-1.318);
  \draw[rounded corners=4pt] (11.682,-5.682) -- (10.318,-4.318);
  \draw[rounded corners=4pt] (12.318,-5.682) -- (16.682,-1.318);
  \draw[rounded corners=4pt] (9.682,-3.682) -- (8.318,-2.318);
  \draw[rounded corners=4pt] (10.318,-3.682) -- (12.682,-1.318);
  \draw[rounded corners=4pt] (8.318,-1.682) -- (8.682,-1.318);
  \draw[rounded corners=4pt] (9,9.45) -- (9,10) -- (-1,10) -- (-1,-10) -- (9,-10) -- (9,-9.45);
  \draw (9,9.45) .. controls (9,9) .. (8.701,8.664);
  \draw (9.682,8.318) .. controls (10,8) .. (10.296,7.661);
  \draw (9.318,7.318) .. controls (9,7) .. (9.223,6.609);
  \draw (7.318,5.318) .. controls (7,5) .. (7.167,4.582);
  \draw (5.318,3.318) .. controls (5,3) .. (4.75,2.626);
  \draw (5.682,2.318) .. controls (6,2) .. (6.201,1.598);
  \draw (5.318,1.318) .. controls (5,1) .. (5,0.55);
  \draw (10.682,1.318) .. controls (11,1) .. (11,0.55);
  \draw (14.682,1.318) .. controls (15,1) .. (15,0.55);
  \draw (9,-9.45) .. controls (9,-9) .. (8.817,-8.589);
  \draw (2.318,-2.318) .. controls (2,-2) .. (2.201,-1.598);
  \draw (1.318,-1.318) .. controls (1,-1) .. (1,-0.55);
  \draw (11.682,-6.318) .. controls (12,-6) .. (12.201,-5.598);
  \draw (10.318,-4.318) .. controls (10,-4) .. (10.109,-3.563);
  \draw (8.318,-2.318) .. controls (8,-2) .. (7.799,-1.598);
  \draw (8.682,-1.318) .. controls (9,-1) .. (9,-0.55);
  \draw (12.682,-1.318) .. controls (13,-1) .. (13,-0.55);
  \draw (16.682,-1.318) .. controls (17,-1) .. (17,-0.55);
  \draw[white, line width=6pt, line cap=butt] (8.682,8.682) .. controls (9,9) .. (9.318,8.682);
  \draw (8.682,8.682) .. controls (9,9) .. (9.318,8.682);
  \draw[white, line width=6pt, line cap=butt] (9.682,7.682) .. controls (10,8) .. (10.318,7.682);
  \draw (9.682,7.682) .. controls (10,8) .. (10.318,7.682);
  \draw[white, line width=6pt, line cap=butt] (8.682,6.682) .. controls (9,7) .. (9.318,6.682);
  \draw (8.682,6.682) .. controls (9,7) .. (9.318,6.682);
  \draw[white, line width=6pt, line cap=butt] (6.682,4.682) .. controls (7,5) .. (7.318,4.682);
  \draw (6.682,4.682) .. controls (7,5) .. (7.318,4.682);
  \draw[white, line width=6pt, line cap=butt] (4.682,2.682) .. controls (5,3) .. (5.318,2.682);
  \draw (4.682,2.682) .. controls (5,3) .. (5.318,2.682);
  \draw[white, line width=6pt, line cap=butt] (5.682,1.682) .. controls (6,2) .. (6.318,1.682);
  \draw (5.682,1.682) .. controls (6,2) .. (6.318,1.682);
  \draw[white, line width=6pt, line cap=butt] (4.682,0.682) .. controls (5,1) .. (5.318,0.682);
  \draw (4.682,0.682) .. controls (5,1) .. (5.318,0.682);
  \draw[white, line width=6pt, line cap=butt] (10.682,0.682) .. controls (11,1) .. (11.318,0.682);
  \draw (10.682,0.682) .. controls (11,1) .. (11.318,0.682);
  \draw[white, line width=6pt, line cap=butt] (14.682,0.682) .. controls (15,1) .. (15.318,0.682);
  \draw (14.682,0.682) .. controls (15,1) .. (15.318,0.682);
  \draw[white, line width=6pt, line cap=butt] (8.682,-8.682) .. controls (9,-9) .. (9.318,-8.682);
  \draw (8.682,-8.682) .. controls (9,-9) .. (9.318,-8.682);
  \draw[white, line width=6pt, line cap=butt] (1.682,-1.682) .. controls (2,-2) .. (2.318,-1.682);
  \draw (1.682,-1.682) .. controls (2,-2) .. (2.318,-1.682);
  \draw[white, line width=6pt, line cap=butt] (0.682,-0.682) .. controls (1,-1) .. (1.318,-0.682);
  \draw (0.682,-0.682) .. controls (1,-1) .. (1.318,-0.682);
  \draw[white, line width=6pt, line cap=butt] (11.682,-5.682) .. controls (12,-6) .. (12.318,-5.682);
  \draw (11.682,-5.682) .. controls (12,-6) .. (12.318,-5.682);
  \draw[white, line width=6pt, line cap=butt] (9.682,-3.682) .. controls (10,-4) .. (10.318,-3.682);
  \draw (9.682,-3.682) .. controls (10,-4) .. (10.318,-3.682);
  \draw[white, line width=6pt, line cap=butt] (7.682,-1.682) .. controls (8,-2) .. (8.318,-1.682);
  \draw (7.682,-1.682) .. controls (8,-2) .. (8.318,-1.682);
  \draw[white, line width=6pt, line cap=butt] (8.682,-0.682) .. controls (9,-1) .. (9.318,-0.682);
  \draw (8.682,-0.682) .. controls (9,-1) .. (9.318,-0.682);
  \draw[white, line width=6pt, line cap=butt] (12.682,-0.682) .. controls (13,-1) .. (13.318,-0.682);
  \draw (12.682,-0.682) .. controls (13,-1) .. (13.318,-0.682);
  \draw[white, line width=6pt, line cap=butt] (16.682,-0.682) .. controls (17,-1) .. (17.318,-0.682);
  \draw (16.682,-0.682) .. controls (17,-1) .. (17.318,-0.682);
\end{tikzpicture}}
\end{tabular}
\caption{$10_{163}$ with 10 leaves.}
\label{fig:sw-10163}
\end{figure}
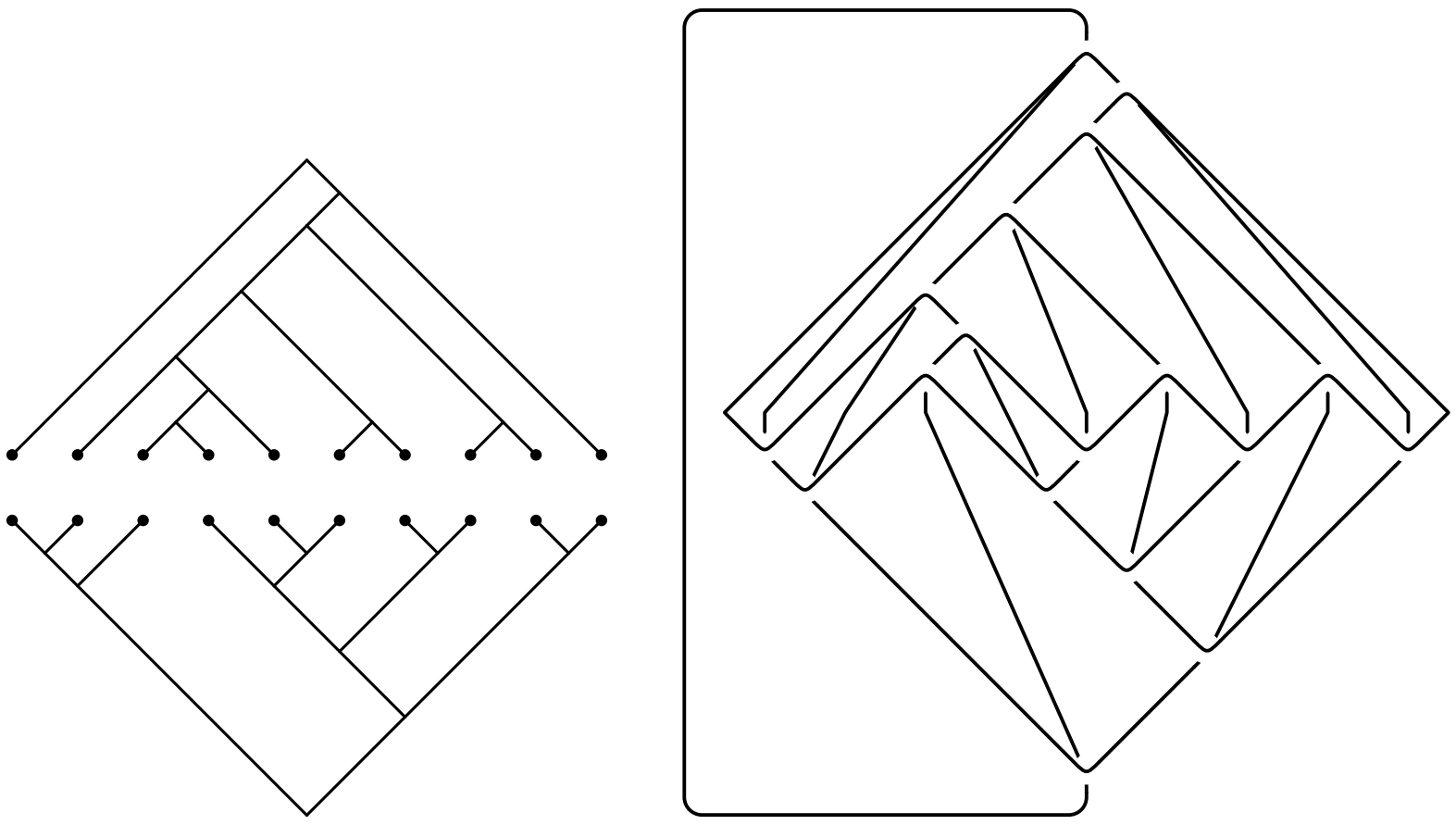

\begin{figure}[ht]
\centering
\begin{tabular}{cc}
\resizebox{0.27\linewidth}{!}{\begin{tikzpicture}[scale=0.32, line cap=round]
  \draw (1.5,1) -- (1,0.5);
  \draw (1.5,1) -- (2,0.5);
  \draw (6.5,1) -- (6,0.5);
  \draw (6.5,1) -- (7,0.5);
  \draw (6,1.5) -- (5,0.5);
  \draw (6,1.5) -- (6.5,1);
  \draw (5.5,2) -- (4,0.5);
  \draw (5.5,2) -- (6,1.5);
  \draw (5,2.5) -- (3,0.5);
  \draw (5,2.5) -- (5.5,2);
  \draw (5.5,3) -- (5,2.5);
  \draw (5.5,3) -- (8,0.5);
  \draw (4.5,4) -- (1.5,1);
  \draw (4.5,4) -- (5.5,3);
  \draw (4,4.5) -- (0,0.5);
  \draw (4,4.5) -- (4.5,4);
  \draw (0.5,-1) -- (0,-0.5);
  \draw (0.5,-1) -- (1,-0.5);
  \draw (2.5,-1) -- (2,-0.5);
  \draw (2.5,-1) -- (3,-0.5);
  \draw (3,-1.5) -- (2.5,-1);
  \draw (3,-1.5) -- (4,-0.5);
  \draw (3.5,-2) -- (3,-1.5);
  \draw (3.5,-2) -- (5,-0.5);
  \draw (4,-2.5) -- (3.5,-2);
  \draw (4,-2.5) -- (6,-0.5);
  \draw (7.5,-1) -- (7,-0.5);
  \draw (7.5,-1) -- (8,-0.5);
  \draw (5,-3.5) -- (4,-2.5);
  \draw (5,-3.5) -- (7.5,-1);
  \draw (4,-4.5) -- (0.5,-1);
  \draw (4,-4.5) -- (5,-3.5);
  \fill (0,0.5) circle (2.5pt) (0,-0.5) circle (2.5pt);
  \fill (1,0.5) circle (2.5pt) (1,-0.5) circle (2.5pt);
  \fill (2,0.5) circle (2.5pt) (2,-0.5) circle (2.5pt);
  \fill (3,0.5) circle (2.5pt) (3,-0.5) circle (2.5pt);
  \fill (4,0.5) circle (2.5pt) (4,-0.5) circle (2.5pt);
  \fill (5,0.5) circle (2.5pt) (5,-0.5) circle (2.5pt);
  \fill (6,0.5) circle (2.5pt) (6,-0.5) circle (2.5pt);
  \fill (7,0.5) circle (2.5pt) (7,-0.5) circle (2.5pt);
  \fill (8,0.5) circle (2.5pt) (8,-0.5) circle (2.5pt);
\end{tikzpicture}} &
\resizebox{0.35\linewidth}{!}{\begin{tikzpicture}[scale=0.32, line cap=round, line join=round, >=stealth, thick]
  \draw[rounded corners=4pt] (7.682,7.682) -- (0.318,0.318) .. controls (0,0) .. (0.318,-0.318) -- (0.682,-0.682);
  \draw[rounded corners=4pt] (7.704,7.661) -- (1.296,0.339) .. controls (1,0) .. (1,-0.45) -- (1,-0.55);
  \draw[rounded corners=4pt] (8.318,7.682) -- (8.682,7.318);
  \draw[rounded corners=4pt] (8.682,6.682) -- (3.318,1.318);
  \draw[rounded corners=4pt] (8.777,6.609) -- (5.223,0.391) .. controls (5,0) .. (5,-0.45) -- (5,-0.55);
  \draw[rounded corners=4pt] (9.318,6.682) -- (10.682,5.318);
  \draw[rounded corners=4pt] (2.682,0.682) -- (2.318,0.318) .. controls (2,0) .. (1.682,-0.318) -- (1.318,-0.682);
  \draw[rounded corners=4pt] (3,0.55) -- (3,0.45) .. controls (3,0) .. (3.238,-0.382) -- (7.762,-7.618);
  \draw[rounded corners=4pt] (3.318,0.682) -- (3.682,0.318) .. controls (4,0) .. (4.318,-0.318) -- (4.682,-0.682);
  \draw[rounded corners=4pt] (10.682,4.682) -- (10.318,4.318);
  \draw[rounded corners=4pt] (11.281,4.649) -- (14.719,0.351) .. controls (15,0) .. (15,-0.45) -- (15,-0.55);
  \draw[rounded corners=4pt] (11.318,4.682) -- (15.682,0.318) .. controls (16,0) .. (15.682,-0.318) -- (15.318,-0.682);
  \draw[rounded corners=4pt] (9.682,3.682) -- (6.318,0.318) .. controls (6,0) .. (5.682,-0.318) -- (5.318,-0.682);
  \draw[rounded corners=4pt] (9.73,3.64) -- (7.27,0.36) .. controls (7,0) .. (6.799,-0.402) -- (6.201,-1.598);
  \draw[rounded corners=4pt] (10.318,3.682) -- (10.682,3.318);
  \draw[rounded corners=4pt] (10.682,2.682) -- (8.318,0.318) .. controls (8,0) .. (7.682,-0.318) -- (6.318,-1.682);
  \draw[rounded corners=4pt] (10.75,2.626) -- (9.25,0.374) .. controls (9,0) .. (8.75,-0.374) -- (7.25,-2.626);
  \draw[rounded corners=4pt] (11.318,2.682) -- (11.682,2.318);
  \draw[rounded corners=4pt] (11.682,1.682) -- (10.318,0.318) .. controls (10,0) .. (9.682,-0.318) -- (7.318,-2.682);
  \draw[rounded corners=4pt] (11.799,1.598) -- (11.201,0.402) .. controls (11,0) .. (10.73,-0.36) -- (8.27,-3.64);
  \draw[rounded corners=4pt] (12.318,1.682) -- (12.682,1.318);
  \draw[rounded corners=4pt] (12.682,0.682) -- (12.318,0.318) .. controls (12,0) .. (11.682,-0.318) -- (8.318,-3.682);
  \draw[rounded corners=4pt] (13,0.55) -- (13,0.45) .. controls (13,0) .. (12.799,-0.402) -- (10.201,-5.598);
  \draw[rounded corners=4pt] (13.318,0.682) -- (13.682,0.318) .. controls (14,0) .. (14.318,-0.318) -- (14.682,-0.682);
  \draw[rounded corners=4pt] (7.682,-7.682) -- (1.318,-1.318);
  \draw[rounded corners=4pt] (8.318,-7.682) -- (9.682,-6.318);
  \draw[rounded corners=4pt] (9.682,-5.682) -- (8.318,-4.318);
  \draw[rounded corners=4pt] (10.318,-5.682) -- (14.682,-1.318);
  \draw[rounded corners=4pt] (7.682,-3.682) -- (7.318,-3.318);
  \draw[rounded corners=4pt] (6.682,-2.682) -- (6.318,-2.318);
  \draw[rounded corners=4pt] (5.682,-1.682) -- (5.318,-1.318);
  \draw[rounded corners=4pt] (8,8.45) -- (8,9) -- (-1,9) -- (-1,-9) -- (8,-9) -- (8,-8.45);
  \draw (8,8.45) .. controls (8,8) .. (7.704,7.661);
  \draw (8.682,7.318) .. controls (9,7) .. (8.777,6.609);
  \draw (3.318,1.318) .. controls (3,1) .. (3,0.55);
  \draw (10.682,5.318) .. controls (11,5) .. (11.281,4.649);
  \draw (10.318,4.318) .. controls (10,4) .. (9.73,3.64);
  \draw (10.682,3.318) .. controls (11,3) .. (10.75,2.626);
  \draw (11.682,2.318) .. controls (12,2) .. (11.799,1.598);
  \draw (12.682,1.318) .. controls (13,1) .. (13,0.55);
  \draw (8,-8.45) .. controls (8,-8) .. (7.762,-7.618);
  \draw (1.318,-1.318) .. controls (1,-1) .. (1,-0.55);
  \draw (9.682,-6.318) .. controls (10,-6) .. (10.201,-5.598);
  \draw (8.318,-4.318) .. controls (8,-4) .. (8.27,-3.64);
  \draw (7.318,-3.318) .. controls (7,-3) .. (7.25,-2.626);
  \draw (6.318,-2.318) .. controls (6,-2) .. (6.201,-1.598);
  \draw (5.318,-1.318) .. controls (5,-1) .. (5,-0.55);
  \draw (14.682,-1.318) .. controls (15,-1) .. (15,-0.55);
  \draw[white, line width=6pt, line cap=butt] (7.682,7.682) .. controls (8,8) .. (8.318,7.682);
  \draw (7.682,7.682) .. controls (8,8) .. (8.318,7.682);
  \draw[white, line width=6pt, line cap=butt] (8.682,6.682) .. controls (9,7) .. (9.318,6.682);
  \draw (8.682,6.682) .. controls (9,7) .. (9.318,6.682);
  \draw[white, line width=6pt, line cap=butt] (2.682,0.682) .. controls (3,1) .. (3.318,0.682);
  \draw (2.682,0.682) .. controls (3,1) .. (3.318,0.682);
  \draw[white, line width=6pt, line cap=butt] (10.682,4.682) .. controls (11,5) .. (11.318,4.682);
  \draw (10.682,4.682) .. controls (11,5) .. (11.318,4.682);
  \draw[white, line width=6pt, line cap=butt] (9.682,3.682) .. controls (10,4) .. (10.318,3.682);
  \draw (9.682,3.682) .. controls (10,4) .. (10.318,3.682);
  \draw[white, line width=6pt, line cap=butt] (10.682,2.682) .. controls (11,3) .. (11.318,2.682);
  \draw (10.682,2.682) .. controls (11,3) .. (11.318,2.682);
  \draw[white, line width=6pt, line cap=butt] (11.682,1.682) .. controls (12,2) .. (12.318,1.682);
  \draw (11.682,1.682) .. controls (12,2) .. (12.318,1.682);
  \draw[white, line width=6pt, line cap=butt] (12.682,0.682) .. controls (13,1) .. (13.318,0.682);
  \draw (12.682,0.682) .. controls (13,1) .. (13.318,0.682);
  \draw[white, line width=6pt, line cap=butt] (7.682,-7.682) .. controls (8,-8) .. (8.318,-7.682);
  \draw (7.682,-7.682) .. controls (8,-8) .. (8.318,-7.682);
  \draw[white, line width=6pt, line cap=butt] (0.682,-0.682) .. controls (1,-1) .. (1.318,-0.682);
  \draw (0.682,-0.682) .. controls (1,-1) .. (1.318,-0.682);
  \draw[white, line width=6pt, line cap=butt] (9.682,-5.682) .. controls (10,-6) .. (10.318,-5.682);
  \draw (9.682,-5.682) .. controls (10,-6) .. (10.318,-5.682);
  \draw[white, line width=6pt, line cap=butt] (7.682,-3.682) .. controls (8,-4) .. (8.318,-3.682);
  \draw (7.682,-3.682) .. controls (8,-4) .. (8.318,-3.682);
  \draw[white, line width=6pt, line cap=butt] (6.682,-2.682) .. controls (7,-3) .. (7.318,-2.682);
  \draw (6.682,-2.682) .. controls (7,-3) .. (7.318,-2.682);
  \draw[white, line width=6pt, line cap=butt] (5.682,-1.682) .. controls (6,-2) .. (6.318,-1.682);
  \draw (5.682,-1.682) .. controls (6,-2) .. (6.318,-1.682);
  \draw[white, line width=6pt, line cap=butt] (4.682,-0.682) .. controls (5,-1) .. (5.318,-0.682);
  \draw (4.682,-0.682) .. controls (5,-1) .. (5.318,-0.682);
  \draw[white, line width=6pt, line cap=butt] (14.682,-0.682) .. controls (15,-1) .. (15.318,-0.682);
  \draw (14.682,-0.682) .. controls (15,-1) .. (15.318,-0.682);
\end{tikzpicture}}
\end{tabular}
\caption{$10_{132}$ with 9 leaves.}
\label{fig:sw-10132}
\end{figure}
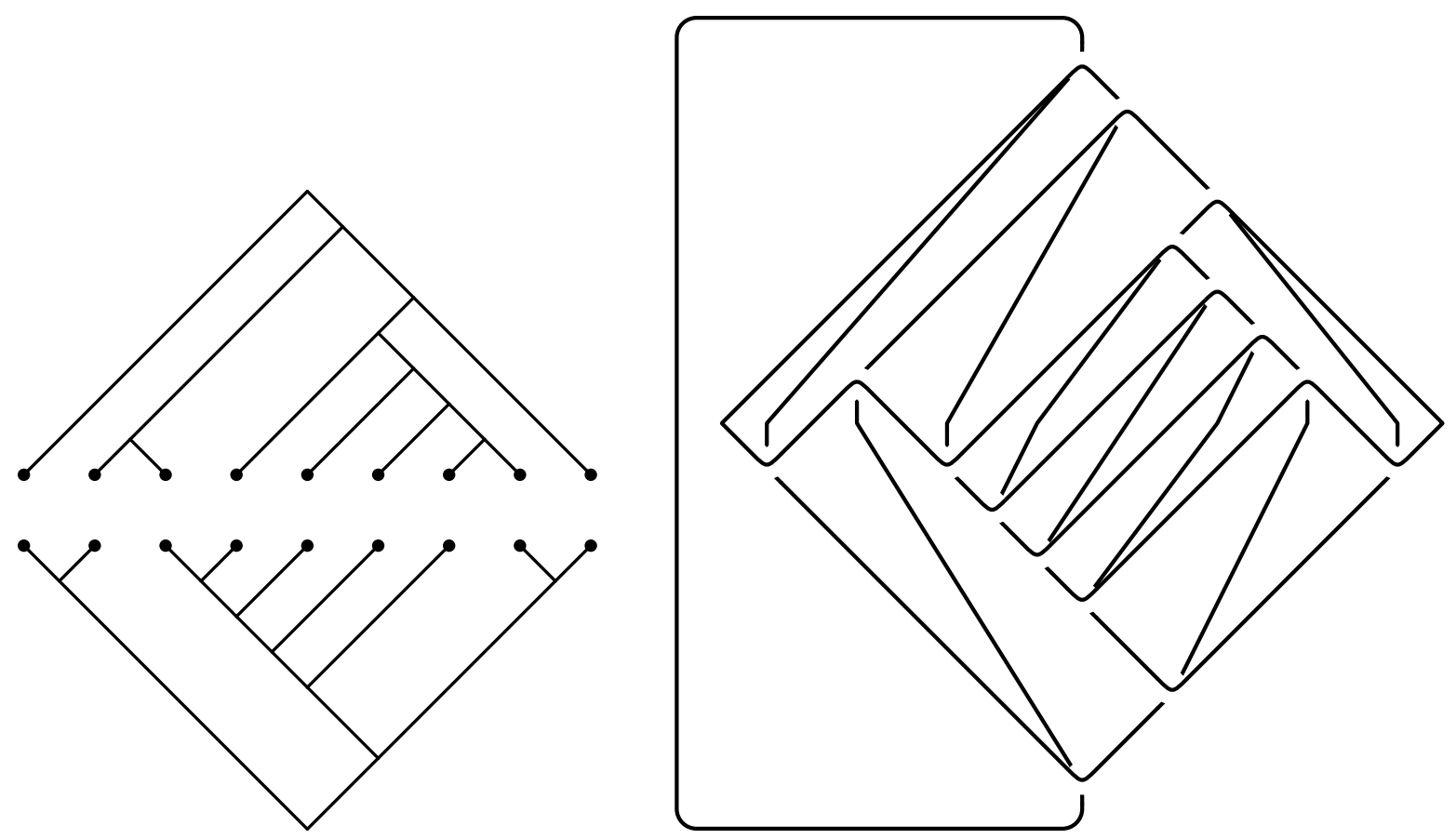

\clearpage
\section{Tables}\label{app:tables}

This appendix collects the long tables: the four indices of the knots with at most 10 crossings
(Table~\ref{tab:small}), the non-invertible knots whose two orientations have different $\vF$- or
$\vFt$-index (Tables~\ref{tab:oriented-Fvec} and~\ref{tab:oriented-F3vec}), and the mutant pairs
distinguished by the indices (Table~\ref{tab:mutants}). The complete tables are in the repository~\cite{Repository} and in the data archive~\cite{Data}.

{\footnotesize\blocktables
\begin{longtable}{lcccc@{\hspace{0.6em}}|@{\hspace{0.6em}}lcccc@{\hspace{0.6em}}|@{\hspace{0.6em}}lcccc}
\toprule
knot & $F$ & $F_3$ & $\vec F$ & $\vec F_3$ & knot & $F$ & $F_3$ & $\vec F$ & $\vec F_3$ & knot & $F$ & $F_3$ & $\vec F$ & $\vec F_3$ \\
\midrule
\endfirsthead
\toprule
knot & $F$ & $F_3$ & $\vec F$ & $\vec F_3$ & knot & $F$ & $F_3$ & $\vec F$ & $\vec F_3$ & knot & $F$ & $F_3$ & $\vec F$ & $\vec F_3$ \\
\midrule
\endhead
\bottomrule
\endfoot
\bottomrule
\caption{The four indices of the prime knots with at most 10 crossings. For the oriented groups and a non-invertible knot with different values, the entry is index$(K)$\,/\,index$(K^r)$, $K$ being the KnotInfo orientation; otherwise it is the index of every orientation. $\geq m$: exact lower bound; for the oriented groups the values and the lower bounds are odd, since knots only occur for odd $n$. $19^{*}$: exact value, given by Theorem~\ref{thm:oriented}\ref{it:reverse} or by a diagram with 19 leaves found by a targeted search (Section~\ref{sec:targeted}). $19^{\dagger}$: non-invertible knot with a diagram with 19 leaves of one orientation (found by the targeted enumeration of Section~\ref{sec:targeted}); the other one has index $19$, $21$ or $23$ (Theorem~\ref{thm:oriented}\ref{it:reverse}).}\label{tab:small}
\endlastfoot
$3_{1}$ & 5 & 4 & 9 & 5 & $4_{1}$ & 7 & 5 & 7 & 5 & $5_{1}$ & 8 & 5 & 15 & 5 \\
$5_{2}$ & 7 & 5 & 13 & 5 & $6_{1}$ & 8 & 5 & 11 & 7 & $6_{2}$ & 8 & 5 & 11 & 7 \\
$6_{3}$ & 9 & 6 & 11 & 7 & $7_{1}$ & 9 & 6 & 21 & 7 & $7_{2}$ & 10 & 6 & 17 & 7 \\
$7_{3}$ & 9 & 6 & 19 & 7 & $7_{4}$ & 8 & 6 & 17 & 7 & $7_{5}$ & 9 & 6 & 19 & 7 \\
$7_{6}$ & 9 & 6 & 13 & 7 & $7_{7}$ & 9 & 6 & 11 & 7 & $8_{1}$ & 10 & 7 & 13 & 9 \\
$8_{2}$ & 10 & 7 & 17 & 7 & $8_{3}$ & 11 & 7 & 15 & 9 & $8_{4}$ & 10 & 7 & 13 & 7 \\
$8_{5}$ & 10 & 7 & 17 & 7 & $8_{6}$ & 9 & 6 & 15 & 7 & $8_{7}$ & 10 & 6 & 15 & 7 \\
$8_{8}$ & 10 & 6 & 13 & 7 & $8_{9}$ & 9 & 7 & 15 & 7 & $8_{10}$ & 10 & 6 & 15 & 7 \\
\end{longtable}
}

{\footnotesize\blocktables\begin{longtable}{lcc@{\hspace{0.6em}}|@{\hspace{0.6em}}lcc@{\hspace{0.6em}}|@{\hspace{0.6em}}lcc@{\hspace{0.6em}}|@{\hspace{0.6em}}lcc}
\toprule
knot & $K$ & $K^r$ & knot & $K$ & $K^r$ & knot & $K$ & $K^r$ & knot & $K$ & $K^r$ \\
\midrule
\endfirsthead
\toprule
knot & $K$ & $K^r$ & knot & $K$ & $K^r$ & knot & $K$ & $K^r$ & knot & $K$ & $K^r$ \\
\midrule
\endhead
\bottomrule
\endfoot
\bottomrule
\caption{Non-invertible knots with at most 13 crossings whose two orientations have different $\vec F$-index, with all tree diagrams with $n\le 17$ (311 knots: 85 with both orientations realized, 226 with one orientation $>17$). $K$ is the orientation of KnotInfo; index$(K)$ = index$(K^{mr})$ and index$(K^r)$ = index$(K^m)$; all the values are odd, since knots only occur for odd $n$. $19^{*}$: exact value, given by Theorem~\ref{thm:oriented}\ref{it:reverse} or by a diagram of size $19$ found by a targeted search (Section~\ref{sec:targeted}).}\label{tab:oriented-Fvec}
\endlastfoot
$10_{88}$ & 15 & 13 & $10_{93}$ & 15 & 17 & $10_{95}$ & 17 & $19^{*}$ & $10_{102}$ & 15 & 17 \\
$11a_{26}$ & 17 & 15 & $11a_{30}$ & 17 & 15 & $11a_{76}$ & 15 & 17 & $11a_{103}$ & 17 & 15 \\
$11a_{125}$ & 17 & $19^{*}$ & $11a_{256}$ & 15 & 17 & $11a_{270}$ & 15 & 17 & $11a_{287}$ & 17 & 15 \\
$11a_{316}$ & $19^{*}$ & 17 & $11n_{67}$ & 15 & 17 & $11n_{80}$ & $19^{*}$ & 17 & $11n_{132}$ & 15 & 17 \\
$11n_{177}$ & 15 & 17 & $12a_{1}$ & 17 & 15 & $12a_{8}$ & $19^{*}$ & 17 & $12a_{57}$ & 17 & $19^{*}$ \\
$12a_{131}$ & 17 & $19^{*}$ & $12a_{202}$ & 17 & 15 & $12a_{214}$ & 15 & 17 & $12a_{262}$ & 15 & 17 \\
$12a_{286}$ & 17 & 15 & $12a_{301}$ & 17 & $19^{*}$ & $12a_{328}$ & $19^{*}$ & 17 & $12a_{342}$ & 17 & 15 \\
$12a_{361}$ & 17 & 15 & $12a_{387}$ & 15 & 17 & $12a_{388}$ & 17 & 15 & $12a_{452}$ & $19^{*}$ & 17 \\
\end{longtable}
}

{\small\blocktables\begin{longtable}{lcc@{\hspace{0.6em}}|@{\hspace{0.6em}}lcc@{\hspace{0.6em}}|@{\hspace{0.6em}}lcc@{\hspace{0.6em}}|@{\hspace{0.6em}}lcc}
\toprule
knot & $K$ & $K^r$ & knot & $K$ & $K^r$ & knot & $K$ & $K^r$ & knot & $K$ & $K^r$ \\
\midrule
\endfirsthead
\toprule
knot & $K$ & $K^r$ & knot & $K$ & $K^r$ & knot & $K$ & $K^r$ & knot & $K$ & $K^r$ \\
\midrule
\endhead
\bottomrule
\endfoot
\bottomrule
\caption{Non-invertible knots with at most 13 crossings whose two orientations have different $\vec F_3$-index, with all tree diagrams with $n\le 11$ (946 knots: 710 with both orientations realized, 236 with one orientation $>11$). $K$ is the orientation of KnotInfo; index$(K)$ = index$(K^{mr})$ and index$(K^r)$ = index$(K^m)$; all the values are odd, since knots only occur for odd $n$. $13^{*}$: exact value, given by Theorem~\ref{thm:oriented}\ref{it:reverse3} or by a diagram of size $13$ found by a targeted search (Section~\ref{sec:targeted}).}\label{tab:oriented-F3vec}
\endlastfoot
$8_{17}$ & 7 & 9 & $10_{147}$ & 7 & 9 & $10_{149}$ & 7 & 9 & $10_{150}$ & 7 & 9 \\
$10_{153}$ & 7 & 9 & $11a_{15}$ & 11 & 9 & $11a_{19}$ & 11 & 9 & $11a_{69}$ & 11 & 9 \\
$11a_{71}$ & 11 & 9 & $11a_{80}$ & 9 & 11 & $11a_{112}$ & 9 & 11 & $11a_{113}$ & 9 & 11 \\
$11a_{115}$ & 9 & 11 & $11a_{141}$ & 9 & 11 & $11a_{163}$ & 9 & 11 & $11a_{167}$ & 9 & 11 \\
$11a_{218}$ & 11 & 9 & $11a_{232}$ & 11 & 9 & $11a_{233}$ & 9 & 11 & $11a_{269}$ & 11 & 9 \\
$11a_{273}$ & 11 & 9 & $11a_{274}$ & 11 & 9 & $11a_{284}$ & 9 & 11 & $11a_{286}$ & 9 & 11 \\
$11a_{300}$ & 9 & 11 & $11a_{326}$ & 11 & 9 & $11a_{330}$ & 9 & 11 & $11a_{331}$ & 11 & 9 \\
$11n_{8}$ & 9 & 7 & $11n_{24}$ & 9 & 7 & $11n_{50}$ & 7 & 9 & $11n_{184}$ & 11 & 9 \\
\end{longtable}
}

{\scriptsize\blocktables\setlength{\tabcolsep}{1.6pt}\begin{longtable}{lcccc@{\hspace{0.6em}}|@{\hspace{0.6em}}lcccc}
\toprule
pair & $F$ & $F_3$ & $\vec F$ & $\vec F_3$ & pair & $F$ & $F_3$ & $\vec F$ & $\vec F_3$ \\
\midrule
\endfirsthead
\toprule
pair & $F$ & $F_3$ & $\vec F$ & $\vec F_3$ & pair & $F$ & $F_3$ & $\vec F$ & $\vec F_3$ \\
\midrule
\endhead
\bottomrule
\endfoot
\bottomrule
\caption{Pairs of mutant knots with at most 13 crossings (Stoimenow~\cite{Stoimenow2010}; 1125 pairs in all) that are distinguished by at least one index (458 pairs). Entries: index of the first knot / index of the second; in bold the indices that distinguish them. For the oriented groups the index is $\min(\mathrm{index}(K), \mathrm{index}(K^r))$, and the lower bounds are odd, since knots only occur for odd $n$.}\label{tab:mutants}
\endlastfoot
$11a_{19}$, $11a_{25}$ & 13/13 & 8/8 & \textbf{17/$\geq 19$} & 9/9 & $11a_{24}$, $11a_{26}$ & 13/13 & 8/8 & \textbf{$\geq 19$/15} & \textbf{9/11} \\
$11a_{44}$, $11a_{47}$ & 13/13 & 8/8 & \textbf{$\geq 19$/17} & \textbf{9/11} & $11a_{57}$, $11a_{231}$ & 12/12 & \textbf{8/9} & $\geq 19$/$\geq 19$ & \textbf{9/11} \\
$11a_{251}$, $11a_{253}$ & 13/13 & \textbf{8/9} & \textbf{$\geq 19$/17} & 11/11 & $11a_{252}$, $11a_{254}$ & 12/12 & 8/8 & \textbf{17/$\geq 19$} & 11/11 \\
$11n_{34}$, $11n_{42}$ & 12/12 & 8/8 & \textbf{$\geq 19$/15} & 9/9 & $11n_{35}$, $11n_{43}$ & \textbf{11/12} & 8/8 & $\geq 21$/$\geq 21$ & 9/9 \\
$11n_{36}$, $11n_{44}$ & 12/12 & 8/8 & \textbf{$\geq 19$/17} & 9/9 & $11n_{39}$, $11n_{45}$ & 13/13 & 8/8 & \textbf{15/$\geq 19$} & 9/9 \\
$11n_{73}$, $11n_{74}$ & \textbf{12/13} & 8/8 & \textbf{$\geq 19$/17} & 9/9 & $11n_{151}$, $11n_{152}$ & 13/13 & 8/8 & \textbf{17/$\geq 19$} & 9/9 \\
$12a_{29}$, $12a_{113}$ & 13/13 & 9/9 & $\geq 19$/$\geq 19$ & \textbf{9/11} & $12a_{30}$, $12a_{157}$ & $\geq 14$/$\geq 14$ & \textbf{8/9} & $\geq 19$/$\geq 19$ & 11/11 \\
$12a_{33}$, $12a_{157}$ & $\geq 14$/$\geq 14$ & \textbf{8/9} & $\geq 19$/$\geq 19$ & 11/11 & $12a_{36}$, $12a_{694}$ & $\geq 14$/$\geq 14$ & 9/9 & $\geq 21$/$\geq 21$ & \textbf{9/11} \\
\end{longtable}
}

\bibliographystyle{amsplain}
\bibliography{refs}

\end{document}